\input amstex
\documentstyle{amsppt}
\NoBlackBoxes
\hcorrection {-0.2cm}
\magnification=1200
\baselineskip=16pt
\font\we=cmb10 at 14.4truept
\topmatter
\title\nofrills 
\centerline{\we Relative Bott-Chern Secondary Characteristic
Classes} 
\endtitle
\author\centerline  {\bf Lin Weng}
\endauthor
\NoRunningHeads
\endtopmatter
\TagsOnRight
\centerline {Department of Mathematics, Graduate School of
Science, Osaka University,}
\centerline {Toyonaka, Osaka 560, Japan}
\vskip 0.2cm
\centerline  {(October 8, 1998)}
\vskip 1.0cm
\noindent
{\bf Abstract.} In this paper, we introduce six axioms
for relative Bott-Chern secondary characteristic classes and
prove the uniqueness and existence theorem for them. Such a
work provides us a natural way  to understand and hence to
prove the   arithmetic Grothendieck-Riemann-Roch theorem.
\vskip 1.5cm
\noindent
{\bf Content}
\vskip 0.20cm
\noindent
Introduction
\vskip 0.20cm
\noindent
1. Classical Bott-Chern secondary characteristic classes
\vskip 0.20cm
\noindent
2. Axioms of relative Bott-Chern secondary characteristic
classes
\vskip 0.20cm
\noindent
3. Deformation to the normal cone
\vskip 0.20cm
\noindent
4. Uniqueness   of relative Bott-Chern
secondary characteristic classes
\vskip 0.20cm
\noindent
5. Some intermediate results 
\vskip 0.20cm
\noindent
6. Proof of the uniqueness of relative Bott-Chern
secondary characteristic classes
\vskip 0.20cm
\noindent
7. Existence of relative Bott-Chern
secondary characteristic classes
\vskip 0.20cm
\noindent
Acknowledgement
\vskip 0.20cm
\noindent
References
\vfill\eject
\centerline {\bf Introduction}
\vskip 0.50cm
\noindent 
This  paper comes partly as a result of our renewed attemp
to simplify and clarify our previous work on the so-called
relative Bott-Chern secondary characteristic classes and an
arithmetic  Grothendieck-Riemann-Roch theorem for local complete
intersection morpohisms done in 1991 (see e.g., [We1,2]).  Our
(initial and final) aim is to expose the natural structure of the
arithmetic Grothendieck-Riemann-Roch theorem 
via a theory for  relative Bott-Chern secondary characteristic
classes.
\vskip 0.20cm
Yet, the works done in [We1,2] are far from being perfect. To make
things worse, as what we found  recently, the original
uniqueness assertion for the relative Bott-Chern secondary
characteristic classes stated in [We1,2] is not entirely
correct. Fortunately, this mistake is by no means a fatal one. It
is our main task here to correct this mistake and to fix some
other flaws (appeared in [We1,2]) so as to make the theory of
relative Bott-Chern secondary characteristic classes and our
arithmetic Grothendieck-Riemann-Roch theorem for local complete
intersection morphisms revive. In this way, we then hope that the
reader can  understand the arithmetic
Grothendieck-Riemann-Roch theorem  at its own disposal.
\vskip 0.20cm
This paper is organized as follows.  We begin with a
review of a theory of the classical Bott-Chern
secondary characteristic classes in Chapter 1.  Then, in Chapter
2, we introduce five of six key axioms for our relative
Bott-Chern secondary characteristic classes, which consist of the
downstairs rule, the projection rule, the functorial rule, the
uniqueness rule with respect to hermitian vector bundles, and the
uniqueness rule with respect to metrized morphisms. The final
axiom, the deformation to the normal cone rule, for  relative
Bott-Chern secondary characteristic classes is introduced in
Chapter 3. After this,  we state   two
uniqueness theorems for the relative Bott-Chern secondary
characteristic classes in Chapter 4.  In Chapter 5, we
prove a few intermediate results as a preparation to
prove the uniqueness theorems. In Chapter 6, we complete the
proof of the uniqueness theorems for relative Bott-Chern
secondary  characteristic classes. Finally, in Chapter 7, we
prove a weak version of the existence theorem for  relative
Bott-Chern secondary characteristic classes by an effective
construction, which is sufficient and necessary for the
application to the arithmetic Grothendieck-Riemann-Roch theorem.
\vfill\eject
\noindent
{\bf 1. Classical Bott-Chern secondary characteristic
classes}
\vskip 0.20cm
\noindent
A) Axioms for classical Bott-Chern classes
\vskip 0.20cm
\noindent
B) Uniqueness of classical Bott-Chern  classes
\vskip 0.20cm
\noindent
C) Existence of  classical Bott-Chern  classes
\vskip 0.20cm
\noindent
Associated to a hermitian vector bundle $(E,\rho)$ on a complex
manifold $M$ are  Chern forms $\text{ch}(E,\rho)$. It is
well-known that the differential forms $\text{ch}(E,\rho)$ are
$d$-closed on $M$. Hence they define  de Rham cohomology classes
$[\text{ch}(E,\rho)]\in H^*(M,{\Bbb R})$. Moreover, the classes
$[\text{ch}(E,\rho)]$ do not depend on the choice of the metric
$\rho$.
\vskip 0.20cm
On the other hand, if we view $\text{ch}(E,\rho)$ as genuine
differential forms,  they do depend on the metric $\rho$.
Hence a natural question is to understand how  Chern forms
depend on  metrics. This problem is solved by Bott and Chern
in [BC], via the so-called classical Bott-Chern secondary
characteristic classes. We in this chapter review the accociated
theory.
\vskip 0.20cm
\noindent
{\bf A. Axioms for classical Bott-Chern classes}
\vskip 0.20cm
\noindent
(A.1) Let $(E,\rho)$ be  a hermitian vector bundle  of
rank $r$ on a complex manifold $M$. It is well-known that there
exists a unique connection $\nabla_{E,\rho}$ on $E$ which is
compatible with the complex structure $\bar\partial$ of $E$ and
preserves the metric $\rho$. The associated curvature is then
defined to be $\nabla_{E,\rho}^2$. 
\vskip 0.20cm
With respect to a local frame $s_U$ over an open trivialization
chart $U$ for $E$, the connection $\nabla_{E,\rho}$ (resp. the
curvature  $\nabla_{E,\rho}^2$) may be represented by an
$r\times r$-matrix $\omega_U$ (resp. $\Omega_U$) of differential
1-forms (resp. differential (1,1)-forms) on
$U$. Furthermore, $\Omega_U$ satisfies the {\it Bianchi
identity}
$$d\Omega=\Omega\wedge\omega-\omega\wedge\omega.\eqno(1.1)$$
Moreover, if $V$ is another open trivialization chart for $E$
with
$s_V$ a local frame of $E$ on $V$,
$g_{VU}:U\cap V\to \text{GL}(r,{\Bbb C})$, the associated
holomorphic transformation,  then  
on $U\cap V$, $s_U=s_V\cdot g_{VU}$, and
$$\Omega_U=g_{VU}^{-1}\cdot \Omega_V\cdot g_{VU}.\eqno(1.2)$$
Here $\Omega_V$ denotes the curvature forms of $(E,\rho)$ on $V$
with respect to the local frame $s_V$.
\vskip 0.20cm
\noindent
(A.2) To facilitate the ensuring discussion, we now recall a
result from algebra. Let
$B\subset {\Bbb R}$ be a subring, and let $\phi\in
B[[T_1,\dots,T_r]]$ be a symmetric power series. For every $k\geq
0$, denote by
$\phi_{[k]}$ the degree $k$ homogeneous component of 
$\phi$. Then there exists a unique polynomial map
$$\Phi_{[k]}: M_r({\Bbb C})\rightarrow {\Bbb C}$$ from the ring of
$r\times r$ complex matrices such that 
\vskip 0.20cm
\item {(1)} $\Phi_{[k]}$ is invariant under the conjugation 
of $\text{GL}(n,{\Bbb C})$;
\item {(2)}
$\Phi_{[k]}(\text{diag}(a_1,\dots,a_r))=\phi_{[k]}(a_1,\dots,
a_r).$
\vskip 0.20cm
Thus, more generally, for any $B$-algebra $A$, we have a
well-defined  map
$$\Phi\,=\,\oplus_{k\geq 0}\Phi_{[k]}: M_r(A)\rightarrow A.$$ 
This then implies that for a nilpotent subalgebra $I$ of $A$, 
we may also have a well-defined map
$$\Phi\,=\,\oplus_{k\geq 0}\Phi_{[k]}: M_r(I)\rightarrow A.$$
\vskip 0.20cm
\noindent
(A.3) We apply (A.1) and (A.2) as follows. With the same notation
as above, if
$(E,\rho)$ is a  hermitian vector bundle of rank
$r$ on a complex  manifold $M$, and  $\phi\in
B[[T_1,\dots,T_r]]$ is a symmetric power series, 
then on a local trivialization chart $U$ of $E$, define
differential forms $\phi(E,\rho;U)$ on $U$ by
$$\phi(E,\rho;U):=
\Phi({1\over{2\pi\sqrt{-1}}}\cdot \Omega_U).$$ Here we set
$A:=A(U):=\oplus_{p\geq 0}A^{p,p}(U)$  be  the space of all
$(p,p)$-differential forms on $U$ and $I=:\oplus_{p\geq
1}A^{p,p}(U)$.
\vskip 0.20cm
Then by (1.2), and A.2.1, $\{(U,\phi(E,\rho;U))\}$ defines global
differential forms on
$M$. Denote these resulting forms by $\phi(E,\rho)$
and call them the {\it characteristic forms} assocated to
$(E,\rho)$ with respect to $\phi$.
\vskip 0.20cm
\noindent
{\bf Theorem.} {\it With the same notation as
above, 
\vskip 0.20cm \noindent{(1)} $\phi(E,\rho)$ are $d$-closed
differential forms on
$M$, i.e. $d\big(\phi(E,\rho)\big)=0$.
\vskip 0.20cm \noindent{(2)} For any map 
 $f: N\rightarrow M$ of complex manifolds,
$f^*\big(\phi(E,\rho)\big)=\phi(f^*E, f^*\rho).$
\vskip 0.20cm \noindent{(3)} If $(E,\rho)$ is a hermitian line bundle $(L,\rho)$,
then $\phi(L,\rho)=\phi(c_1(L,\rho))$. Here 
$c_1(L,\rho)$ is the so-called first Chern form of $(L,\rho)$
defined to be the differential form
$dd^c\big(-\log\|s\|_\rho^2\big).$
\vskip 0.20cm \noindent{(4)} The de Rham cohomology classes $[\phi(E,\rho)]\in
H^*(M,{\Bbb R})$ of
$\phi(E,\rho)$ do not  depend on the choice of $\rho$, but the
forms
$\phi(E,\rho)$ themselves do  depend on $\rho$.}
\vskip 0.20cm
The proof of this theorem may be found in  standard
textbooks which contain the theory of characteristic forms. (In
fact, (1) is a direct consequence of the Bianchi
identity (1.1).) So we omit it.  Instead, we point out that
first,  (1), (2), (3) and (4) charactrize the characteristic
classes associated to $\phi$ uniquely; and secondly, (4) is the
starting point of our story.
\vskip 0.20cm
\noindent
(A.4) Now a natural question is   how  characteristic forms
depend on hermitian metrics. This
problem is solved by Bott and Chern via the so-called
the classical Bott-Chern secondary characteristic classes [BC].
\vskip 0.45cm
To state their fundamental result, following [GS2], we start with
the following axioms for the  classical Bott-Chern  secondary
characteristic classes, $\phi_{\text{BC}}(E., \rho.)$, with
values in $\tilde A(M):=
A(M)/(\text{Im}\partial+\text{Im}\bar\partial)$, associated  to
a symmetric power series $\phi$, a short exact sequence of
vector bundles $E.: 0\rightarrow E_1\rightarrow E_2\rightarrow 
E_3\rightarrow 0$ on $M$, and hermitian metrics $\rho_j$ on $E_j$
for  $j=1,2,3$. 
\vskip 0.20cm
\noindent
{\it Axiom 1.} ({\bf Downstairs Rule}) In ${A}(M)$,
$$dd^c\phi_{\text{BC}}(E.,\rho.)=\phi(E_2,\rho_2)\,-\,
\phi(E_1\oplus E_3,\rho_1\oplus \rho_3).$$
\vskip 0.20cm
\noindent{\it Axiom 2.} ({\bf Functorial Rule}) For any morphism
$f: N\rightarrow M$ of complex manifolds,
$$f^*\big(\phi_{\text{BC}}(E.,\rho.)\big)\,=\,\phi_{\text{BC}}(f^*E.,f^*\rho.).$$
\vskip 0.20cm
\noindent{\it Axiom 3.} ({\bf Uniqueness Rule}) 
If $(E.,\rho.)$ splits,
 i.e. $(E_2,\rho_2)=(E_1\oplus E_3,\rho_1\oplus
 \rho_3)$, then $$\phi_{\text{BC}}(E.,\rho.)=0.$$
\vskip 0.20cm
\noindent
(A.5) We have the following
\vskip 0.20cm
\noindent
{\bf Theorem.} ({\bf Existence and Uniqueness of Classical
Bott-Chern Secondary Characteristic  Classes})  ([BC] \& [GS2])
{\it For any symmetric power series $\phi$, there exists a unique
construction
$\phi_{\text{BC}}$ such that  associated to each   short exact
sequence of vector bundles on a  complex manifold $M$ $$E.:\qquad
0\rightarrow E_1\rightarrow E_2\rightarrow  E_3\rightarrow 0$$ 
and hermitian metrics
$\rho_j$ on $E_j$ for 
$j=1,\, 2,\,3$, is a unique element
$\phi_{\text{BC}}(E.,\rho.)\in \tilde{A}(M)$ which satisfies 
Axioms 1, 2, and 3.}
\vskip 0.20cm
\noindent
{\it Remark 1.1.} We will call the unique element
$\phi_{\text{BC}}(E.,\rho.)\in \tilde{A}(M)$ the {\it classical
Bott-Chern secondary characteristic classes}, or {\it classical
Bott-Chern  classes} or simply {\it Bott-Chern classes},
associated to $(E.,\rho.)$ and $\phi$. Useful examples are given
when $\phi$ corresponds to ch, the Chern characteristic
class, or td, the Todd characteristic class. In these cases,
we denote them by $\text{ch}_{\text{BC}}$ and
$\text{td}_{\text{BC}}$ respectively.
\vskip 0.20cm
\noindent
{\it Remark 1.2.} $\phi_{\text{BC}}(E.,\rho.)$'s are not
differential forms, rather  they are in $\tilde A$, the quotient
space of differential forms modulo the exact forms.
\vskip 0.20cm
\noindent
{\it Remark 1.3.} When $E_3=0$, so that $E_1=E_2=:E$, we write
$\phi_{\text{BC}}(E.,\rho.)$ simply as
$\phi_{\text{BC}}(E;\rho_2,\rho_1)$. Thus, by Axiom 1,  
$$dd^c\phi_{\text{BC}}(E;\rho_2,\rho_1)=\phi(E,\rho_2)-\phi(E,
\rho_1).$$ So $\phi_{\text{BC}}(E;\rho_2,\rho_1)$
measures how $\phi(E,\rho)$ depends on the metric $\rho$.
\vskip 0.20cm
\noindent
{\bf B. Uniqueness of classical Bott-Chern classes}
\vskip 0.20cm
\noindent
(B.1) For the time being, assume that there exist elements
$\phi_{\text{BC}}(E.,\rho.)\in\tilde A(M)$  satisfying Axioms 1,
2 and 3. Our aim in this section is to prove the uniqueness
part of Theorem A.5 by using a ${\Bbb P}^1$-deformation following
[GS2].
\vskip 0.20cm
Let $1_\infty$ be a section of the bundle ${\Cal O}_{{\Bbb
P}^1}(1)$, which vanishes at $\infty$ and has the
value $1$ at 0. So we have an inclusion
$\text{Id}_{E_1}\otimes 1_\infty: E_1\rightarrow E_1(1)$. Set
$$DE_1:=E_1(1):=E_1\otimes {\Cal O}_{{
\Bbb P}^1}(1),\ \  DE_2:=(E_2\oplus E_1(1))/E_1,\ \ 
DE_3:=E_3.$$
 Then we obtain an exact sequence of vector bundles
on $M\times {\Bbb P}^1$:
$$DE.:\ \ \ \ \ 0\rightarrow DE_1
\rightarrow  DE_2\rightarrow DE_3\rightarrow 0.$$
Easily we see that
the restriction of $DE.$ on $M\times\{0\}$ gives $E.$ on $M$ while
the restrction of $DE.$ on $M\times\{\infty\}$ gives the splitting
exact sequence $0\to E_1\to E_1\oplus E_3\to E_3\to 0$ on
$M$. (So we deform $E.$ to a split exact sequence.) By a partition
of unity, we may further choose hermitian metrics $D\rho_i$ on
$DE_i$,
$i=1,2,3$, such that the induced metrics on  $M\times\{0\}$
(resp. on
$M\times\{\infty\}$) via restrictions coincide with the original
metrics $\rho_i$ (resp. splits). Now define an element in $A(M)$
by setting
$$\phi_{\text{BC}}'(E.,\rho.):=\int_{{\Bbb
P}^1}[\log|z|^2]\cdot \Big(\phi(DE_2,D\rho_2)\Big).$$ By an
abuse of notation, we  use $\phi_{\text{BC}}'(E.,\rho.)$ to
denote its image in $\tilde A(M)$ as well.
\vskip 0.20cm
To prove the uniqueness, it suffices to justify the following
\vskip 0.20cm
\noindent
{\bf Claim.} {\it With the same notation as above,
\vskip 0.20cm
\noindent 
(1) $\phi_{\text{BC}}'(E.,\rho.)\in
\tilde A(M)$ does not depend on  $D\rho_2$;
\vskip 0.20cm
\noindent
(2) $\phi_{\text{BC}}'(E.,\rho.)=\phi_{\text{BC}}(E.,\rho.)$ in
$\tilde A(M)$.}
\vskip 0.20cm
\noindent
(B.2) To prove Claim B.1.1, we use a ${\Bbb P}^1$-deformation
again.  Suppose that $D\rho_2$ and $D{\rho}_2'$ are metrics on
$DE_2$ such that they induce the same metrics on $M\times \{0\}$
as well as on  $M\times \{\infty\}$.
Consider the product
$M\times {\Bbb P}^1\times {\Bbb P}^1$ with  points $(y,z,u)$. We
have the following natural maps:
$$M\times{\Bbb  P^1}\buildrel i_u^{12}\over \hookrightarrow M\times {\Bbb  P}^1\times {\Bbb  P}^1
\buildrel p_{12}\over \rightarrow M\times {\Bbb 
P}^1,\quad M\times{\Bbb  P^1}\buildrel i_z^{13}\over
\hookrightarrow M\times {\Bbb  P}^1\times {\Bbb  P}^1
\buildrel p_{13}\over \rightarrow M\times {\Bbb  P}^1,$$
defined by $$i_u^{12}(y,z):=(y,z,u), \ \
p_{12}(y,z,u):=(y,z),\quad i_z^{13}(y,u):=(y,z,u),\ \
p_{13}(y,z,u):=(y,u).$$ Also let $p_1:M\times {\Bbb 
P}^1\rightarrow  M$ be the  projection to the first factor. Then 
by a partition of unity,  we may find a metric
$\tau$ on the  bundle $p_{12}^*DE_2$ such that
\vskip 0.20cm
\noindent {(i)} $(i_0^{12})^*(p_{12}^* DE_2,\tau)\simeq
(DE_2,D\rho_2)$ and  $(i_\infty^{12})^*(p_{12}^* DE_2,\tau)\simeq
(DE_2,D\rho_2')$;
\vskip 0.20cm
\noindent {(ii)} $(i_0^{13})^*(p_{12}^*DE_2,\tau)\simeq
p_1^*(E_2,\rho_2)$ and
$(i_\infty^{13})^*(p_{12}^*DE_2,\tau)\simeq  p_1^*(E_1\oplus
E_3,\rho_1\oplus
\rho_3)$.
\vskip 0.20cm
Hence, 
$$\eqalign {~&\int_{{\Bbb 
P}^1}[\log |z|^2]\,\phi(DE_2,D\rho_2)-
\int_{{\Bbb  P}^1}[\log |z|^2]\,\phi(DE_2,D\rho_2')\cr
=&\int_{{\Bbb  P}^1}[\log
|z|^2]\,\Big(\phi(DE_2,D\rho_2)-\phi(DE_2,D\rho_2')\Big)\cr
=&\int_{{\Bbb  P}^1}[\log
|z|^2]\,\Big(\phi\big((i_0^{12})^*(p_{12}^*DE_2,\tau)\big)
-\phi\big((i_\infty^{12})^*(p_{12}^*DE_2,\tau)\big)\Big)\cr
=&\int_{{\Bbb  P}^1\times {\Bbb  P}^1}[\log |z|^2]\,[\log |u|^2]
\Big(d_ud^c_u\big(\phi(p_{12}^*DE_2,\tau)\big)\Big)\quad(\text{by\
Stokes'\ formula}).\cr}$$ But if we 
let $\partial=\partial_M+\partial_z+\partial_u$ and $\bar\partial=\bar\partial_M
+\bar\partial_z+\bar\partial_u$ be the differentials on $M\times {\Bbb  P}^1\times {\Bbb  P}^1$,
then by the fact that characteristic forms are $d$-closed, we have
$$\eqalign {\int_{{\Bbb  P}^1\times {\Bbb  P}^1}&[\log |z|^2]\,
[\log|u|^2]\Big(d_ud^c_u\big(\phi(p_{12}^*DE_2,\tau)\big)\Big)
\cr
=&\int_{{\Bbb  P}^1\times {\Bbb  P}^1}[\log |z|^2]\,[\log
|u|^2]\Big(d_zd^c_z\big(\phi(p_{12}^*DE_2,\tau)\big)\Big)
.\cr}$$ Thus, using 
Stokes' formula again, we have
$$\eqalign {~&\int_{{\Bbb  P}^1}[\log |z|^2]\,\phi(DE_2,D\rho_2)-
\int_{{\Bbb  P}^1}[\log |z|^2]\,\phi(DE_2,D\rho_2')\cr
=&\int_{{\Bbb  P}^1}[\log
|u|^2]\,\Big(\phi\big((i_0^{13})^*(p_{12}^*DE_2,\tau)\big)
-\phi\big((i_\infty^{13})^*(p_{12}^*DE_2,\tau)\big)\Big)\cr
=&\int_{{\Bbb  P}^1}[\log
|u|^2]\,\Big(p_1^*\big(\phi(E_2,\rho_2)-\phi(E_1\oplus
E_3,\rho_1\oplus\rho_3)\big)\Big)\cr 
=&\int_{{\Bbb  P}^1}[\log
|u|^2]\,\Big(p_1^*\big(\phi(E_2,\rho_2)-\phi(E_1\oplus
E_3,\rho_1\oplus\rho_3)\big)\Big)\cr =&0.\cr}$$ Here, in the last
step, we use the fact that $p_1^*(\phi(E_2,\rho_2)-\phi(E_1\oplus
E_3,\rho_1\oplus\rho_3))$ is a constant form with respect to
${\Bbb  P}^1$  and hence its integration with respect to
$\log |u|^2$ over ${\Bbb P}^1$ is identically zero. This proves
Claim B.1.1.
\vskip 0.20cm
\noindent
(B.3) Next we  prove  Claim B.1.2. By Axiom 1, 
 $$\eqalign {~&\int_{{\Bbb  P}^1}[\log |z|^2]
\,d_{M\times {\Bbb  P}^1}d^c_{M\times {\Bbb 
P}^1}\phi_{\text{BC}}(DE.,D\rho.)\cr =&\int_{{\Bbb  P}^1}[\log
|z|^2]\,\phi(DE_2,D\rho_2) -\int_{{\Bbb  P}^1}[\log
|z|^2]\,\phi(DE_1\oplus DE_3, D\rho_1\oplus D\rho_3)\cr
 =&\int_{{\Bbb  P}^1}[\log
|z|^2]\,\phi(DE_2,D\rho_2).\cr}\eqno(1.3)$$
Here, in the last step, we use the fact that $\phi(DE_1\oplus
DE_3, D\rho_1\oplus D\rho_3)$ remains  the same  if
we make a change from $z$ to $z^{-1}$. Therefore, in $\tilde
A(M)$,
$$\eqalign {~&\phi_{\text{BC}}'(E.,\rho.)\cr
=&\int_{{\Bbb 
P}^1}[\log |z|^2]\,\phi(DE_2,D\rho_2)\cr
=&\int_{{\Bbb  P}^1}[\log
|z|^2]\,d_{M\times {\Bbb  P}^1}d^c_{M\times {\Bbb 
P}^1}\phi_{\text{BC}}
(DE.,D\rho.)\quad(\text{by\ (1.3)})\cr 
=&i_0^*\phi_{\text{BC}}
(DE.,D\rho.)-i_\infty^*\phi_{\text{BC}}
(DE.,D\rho.)\cr 
=&\phi_{\text{BC}}(E.,\rho.)-i_\infty^*\phi_{\text{BC}}
(DE.,D\rho.)\quad(\text{by\ Axiom\ 2\ as}\
i_0^*(DE.,D\rho.)=(E.,\rho.))\cr
=&\phi_{\text{BC}}(E.,\rho.)-0\quad(\text{by\ Axiom\ 3\ as}\
i_\infty^*(DE.,D\rho.)\ \text{splits})\cr
=&\phi_{\text{BC}}(E.,\rho.).\cr}$$ 
This completes the proof of Claim B.1.2 and hence the
uniqueness of the classical Bott-Chern secondary characteristic
classes.
\vskip 0.20cm
\noindent
{\bf C. The existence of classical Bott-Chern classes}
\vskip 0.20cm
\noindent
(C.1) In this section, we show that there exists a construction
$\phi_{\text{BC}}$ which satisfies 
Axioms 1, 2 and 3 in (A.4). 
\vskip 0.20cm
To do so, as above, set
$$\phi_{\text{BC}}(E.,\rho.):=\int_{{\Bbb
P}^1}[\log|z|^2]\cdot \Big(\phi(DE_2,D\rho_2)\Big).$$ By the fact
that $dd^c[\log |z|^2]=\delta_0-\delta_\infty,$ and
characteristic forms are $d$-closed  and functorial,
(see e.g.,  Theorem A.3(1) and (2),) 
$\phi_{\text{BC}}(E.,\rho.)$ satisfies Axiom 1.
\vskip 0.20cm
\noindent
(C.2) By the fact that $\phi_{\text{BC}}(E.,\rho.)$ does not
depend on  $D\rho_2$, i.e.,
Claim B.1.1, and  characteristic
forms are functorial, i.e., Theorem A.1.3,
  $\phi_{\text{BC}}(E.,\rho.)$ satisfies Axiom 2 as well.
\vskip 0.20cm
\noindent
(C.3) If $(E.,\rho.)$ splits, we may take $(DE.,D\rho.)$ as the
pull-back of $(E.,\rho.)$ via the natural projection to $M\times
{\Bbb P}^1$. In this case, $\phi(DE_2,D\rho_2)$ does not depend
on $z$. Hence $$\phi_{\text{BC}}(E.,\rho.):=\int_{{\Bbb
P}^1}[\log|z|^2]\cdot \Big(\phi(DE_2,D\rho_2)\Big)=0.$$ This
shows that $\phi_{\text{BC}}(E.,\rho.)$ satisfies Axiom 3 too.
This then ends the proof of the uniqueness and the existence
of classical Bott-Chern secondary characteristic classes.
\vskip 0.20cm
\noindent
(C.4) We end this chapter with the following remarks.  At the very
beginning, it is not quite clear why the secondary characteristic
classes should be in
$\tilde A(M)$, a quotient of $A(M)$. One possible
explanation is that  with the limitation of the methods
employed, the uniqueness can only be established after modulo
exact forms. On the other hand, we may equally use other
spaces as well, e.g.,  $\hat
A(M):=A(M)/(\text{exact\ forms}+d-\text{closed\ forms})$. If we
wish,  we may justify this later choice of $\hat A(M)$ by arguing
that, after all,  
$\phi(E,\rho)$ is equal to $\phi(E,c\cdot\rho)$ for any positive
constant $c$, yet, in general,
$\phi_{\text{BC}}(E;\rho,c\cdot\rho)\not=0$ in $\tilde A(M)$, (so
it is too sensitive,) and
the $d$-closed forms have already been taken care of in the
theory of characteristic classes (so in the theory of secondary
characteristic classes, 
$d$-closed forms should be simply viewed  as zero). As a matter
of fact, later in Section 7.B, we  use yet  a third  space to
understand classical Bott-Chern classes.
\vfill\eject
\vskip 0.20cm
\noindent
{\bf 2. Axioms for relative Bott-Chern Secondary
Characteristic Classes}
\vskip 0.20cm
\noindent
A) Downstairs Rule
\vskip 0.20cm
\noindent
B) Projection Rule
\vskip 0.20cm
\noindent
C) Functorial Rule
\vskip 0.20cm
\noindent
D) Uniqueness Rule w.r.t. Hermitian Bundles
\vskip 0.20cm
\noindent
E) Uniqueness Rule w.r.t. Metrized Morphisms
\vskip 0.20cm
\noindent
In this chapter,  we naturally expose first five of six key axioms
for the so-called relative Bott-Chern secondary characteristic
classes associated to smooth proper morphisms of compact K\"ahler
manifolds.
\vskip 0.20cm
\noindent
{\bf A. Downstairs Rule}
\vskip 0.20cm
\noindent
(A.1) Let $X$ and $Y$ be two compact K\"ahler manifolds. Let
$f:X\to Y$ be a smooth, proper morphism between $X$ and $Y$. For any
$y\in Y$, denote by $X_y$ the fiber of $f$ at $y$. Fix a hermitian
metric
$\tau_f$ on $T_f$, the relative tangent vector bundle of $f$. Moreover, we assume that the
metric on $X_y$ induced by the restriction of $(T_f,\tau_f)$ is K\"ahler for all $y\in Y$.
\vskip 0.20cm
Let $(E,\rho)$ be a hermitian vector bundle on $X$. Assume that $E$ is $f$-acyclic, i.e., all
higher direct images $R^if_*E$ vanish, where $i>0$.
 As a direct
consequence of this assumption,  the direct image
$f_*E$ of $E$ is a vector bundle on $Y$. Moreover, $f_*E$ admits
a natural metric $L^2(\rho,\tau)$, the $L^2$-metric induced by
taking integration of $\rho$ along $f$ with respect to $\tau_f$.
In the sequel, we call $(f:X\to Y;E,\rho;T_f,\tau_f)$, or simply
$(E,\rho;f,\tau_f)$ a {\it properly metrized datum} if all the
above conditions are satisfied.
\vskip 0.20cm
 Associated to a properly metrized datum $(E,\rho;f,\tau_f)$ is
the so-called {\it relative Bott-Chern secondary characteristic
classes}
$\text{ch}_{\text{BC}}(E,\rho;f,\tau_f)\in\widetilde A(Y)$
which may be essentially characterized by the following 
axioms.
\vskip 0.20cm
\noindent
(A.2) {\it  Axiom 1.} ({\bf Downstairs Rule}) In $A(Y)$, 
$$dd^c\text{ch}_{\text{BC}}(E,\rho;f,\tau_f)=f_*\Big(\text{ch}(E,\rho)\cdot\text{td}(T_f,\tau_f)\Big)
-\text{ch}\Big(f_*E,L^2(\rho,\tau)\Big).$$
\vfill\eject
\vskip 0.20cm
\noindent
{\bf B. Projection Rule}
\vskip 0.20cm
\noindent
(B.1) Let $(F,\rho_F)$ be a hermitian vector bundle on $Y$. Then
$E\otimes f^*F$ is $f$-acyclic as well. (Indeed, we know that
$f_*(E\otimes f^*F)=f_*E\otimes F.$) So, $(E\otimes f^*F,\rho\otimes
f^*\rho_F;f,\tau_f)$ is again a properly metrized datum.
Thus, we also have 
the relative Bott-Chern secondary characteristic classes
$\text{ch}_{\text{BC}}(E\otimes f^*F,\rho\otimes
f^*\rho_F;f,\tau_f)\in \tilde A(Y)$. 
\vskip 0.20cm
\noindent
(B.2) {\it Axiom 2.} ({\bf Projection Rule}) In $\tilde A(Y)$,
$$\text{ch}_{\text{BC}}(E\otimes f^*F,\rho\otimes
f^*\rho_F;f,\tau_f)=\text{ch}_{\text{BC}}(E,\rho;f,\tau_f)\cdot
\text{ch}(F,\rho_F).$$
\vskip 0.20cm
\noindent
{\bf C. Functorial Rule}
\vskip 0.20cm
\noindent
(C.1) Let $g:Y'\to Y$ be a morphism between
two compact K\"ahler manifolds. Then we have the Cartesian product
$$\matrix X\times_YY'&\buildrel g_f\over\to&X\\
&&\\
f_g\downarrow&&\downarrow f\\
&&\\
Y'&\buildrel g\over\to&Y.\endmatrix$$
Define the pull-back of $(E,\rho;f,\tau_f)$ via $g_f$ by
$$g_f^*(E,\rho;f,\tau_f)=(g_f^*(E,\rho);f_g,g_f^*\tau_f).$$
Easily, we see that
\vskip 0.20cm
\noindent
(i) by the $f$-acyclic condition on $E$, $g_f^*E$ is
$f_g$-acyclic;
\vskip 0.20cm
\noindent 
(ii) $g_f^*\tau_f$ gives a hermitian metric on the relative
tangent bundle $T_{f_g}$ for the smooth, proper morphism $f_g$
which induces K\"ahler metrics on its fibers.
\vskip 0.20cm
Hence, $g_f^*(E,\rho;f,\tau_f)$ is again a properly metrized
datum, and we may  have the corresponding relative Bott-Chern
secondary characteristic classes
$\text{ch}_{\text{BC}}\Big(g_f^*(E,\rho;f,\tau_f)\Big)
\in\widetilde A(Y')$ as well.
\vskip 0.20cm
\noindent
(C.2) {\it Axiom 3.} ({\bf Functorial Rule}) In $\tilde A(Y')$, 
$$g_f^*\Big(\text{ch}_{\text{BC}}(E,\rho;f,\tau_f)\Big)=
\text{ch}_{\text{BC}}\Big(g_f^*(E,\rho;f,\tau_f)\Big).$$
\vfill\eject
\vskip 0.20cm
\noindent
{\bf D. Uniqueness Rule w.r.t. Hermitian Bundles}
\vskip 0.20cm
\noindent
(D.1) Let $E.: 0\to
E_1\to E_2\to E_3\to 0$ be an exact sequence of
$f$-acyclic vector bundles on $X$. Then the direct image of $E.$
gives an exact sequence of vector bundles on
$Y$ as follows;
$$f_*E.:\qquad 0\to f_*E_1\to f_*E_2\to f_*E_3\to 0.$$
If we put hermitian metrics $\rho_i$ on $E_i$, $i=1,2,3$,
we then get naturally
$L^2$-metrics $L^2(\rho_i,\tau_f)$ on $f_*E_i$. By Theorem 1.A.5, 
there exists the so-called (classical)
Bott-Chern secondary characteristic classes associated to both
$(E.,\rho.)$ and
$(f_*E.,L^2(\rho.,\tau_f))$. Denote these two classical
Bott-Chern secondary characteristic classes by
$\text{ch}_{\text{BC}}(E.,\rho.)\in \widetilde A(X)$ and
$\text{ch}_{\text{BC}}(f_*E.,L^2(\rho.,\tau_f))\in\widetilde
A(Y)$ respectively. Surely, for each $(E_i,\rho_i)$, $i=1,2,3$,
we have the corresponding properly metrized datum
$(E_i,\rho_i;f,\tau_f)$, and hence the relative Bott-Chern
secondary characteristic classes
$\text{ch}_{\text{BC}}(E_i,\rho_i;f,\tau_f)\in\widetilde A(Y)$.
\vskip 0.20cm
\noindent
(D.2) {\it Axiom 4.} ({\bf Uniqueness Rule w.r.t. Hermitian
Bundles}) In $\tilde A(Y)$, 
$$\eqalign{~&\text{ch}_{\text{BC}}(E_2,\rho_2;f,\tau_f)
-\text{ch}_{\text{BC}}(E_1,\rho_1;f,\tau_f)
-\text{ch}_{\text{BC}}(E_3,\rho_3;f,\tau_f)\cr
=&f_*\Big(\text{ch}_{\text{BC}}(E.,\rho.)\cdot\text{td}(T_f,\tau_f)\Big)
-\text{ch}_{\text{BC}}\Big(f_*E.,L^2(\rho.,\tau_f)\Big).\cr}$$
\vskip 0.20cm
\noindent
{\bf E. Uniqueness Rule w.r.t. Metrized Morphisms}
\vskip 0.20cm
\noindent
(E.1) Let $Z$ be a
compact K\"ahler manifold and $g:Y\to Z$ be a smooth,
proper morphism. So we
have a triangle of smooth fibrations of compact K\"ahler manifolds
among compact K\"ahler manifolds as follows;
$$\matrix X&\buildrel f\over\to&Y\\
&&\\
&g\circ f\searrow &\downarrow g\\
&&\\
&&Z.\endmatrix$$
Naturally, we have the short exact sequence of relative tangent bundles
$$T.: \hskip 1.5cm 0\to T_f\to T_{g\circ f}\to f^*T_g\to 0,$$
where $T_f$, $T_g$  and $T_{g\circ f}$ stand for the relative
tangent bundles associated to the smooth, proper morphisms $f$, $g$ and
$g\circ f$ respectively. Put hermitian metrics $\tau_f$, $\tau_g$
and $\tau_{g\circ f}$ on
$T_f$, $T_g$  and $T_{g\circ f}$ respectively such that the 
induced metrics on all 
corresponding fibers are K\"ahler.  By Theorem 1.A.5, we have
the corresponding (classical) Bott-Chern secondary characteristic
classes
$\text{td}_{\text{BC}}(T.,\tau.)\in\widetilde A(X)$ with respect
to the Todd characteristic class.
\vskip 0.20cm
\noindent
(E.2) Recall that $E$ is $f$-acyclic so that $f_*E$ is a vector
bundle on $Y$. In addition, assume that $E$ is $g\circ f$-acyclic
and
$f_*E$ is $g$-acyclic. As a direct consequence, 
$(g\circ f)_*E$ is a vector bundle on $Z$ and is equal to 
$g_*(f_*E)$. All this gives us the properly metrized data
$(E,\rho;f,\tau_f)$,
$(f_*E,L^2(\rho,\tau_f);g,\tau_g)$, and $(E,\rho;g\circ
f,\tau_{g\circ f})$, and hence their associated
 relative Bott-Chern secondary characteristic
classes
$\text{ch}_{\text{BC}}(E,\rho;f,\tau_f)\in \widetilde A(Y),$
$\text{ch}_{\text{BC}}(f_*E,L^2(\rho,\tau_f);g,\tau_g)\in\widetilde
A(Z)$ and
$\text{ch}_{\text{BC}}(E,\rho;g\circ f,\tau_{g\circ
f})\in\widetilde A(Z).$
\vskip 0.20cm
Moreover, on $Z$, we have a short exact sequence of vector 
bundles $$0\to 
g_*(f_*E)\to (g\circ f)_*E\to 0.$$ Put the $L^2$-metrics
$L^2(\rho,\tau_{g\circ f})$ and
$L^2(L^2(\rho,\tau_f),\tau_g)$ on $(g\circ f)_*E$ and on $g_*(f_*E)$ respectively, we then
have the associated (classical) Bott-Chern secondary characteristic classes,
which, as in Remark 1.1.3 in 1.A.5, we denote by
$$\text{ch}_{\text{BC}}\Big((g\circ f)_*E;L^2(\rho,\tau_{g\circ
f}),L^2(L^2(\rho,\tau_f),\tau_g)\Big)\in\widetilde A(Z).$$
\vskip 0.20cm
\noindent
(E.3) {\it Axiom 5.} ({\bf Uniqueness Rule w.r.t. Metrized
Morphisms}) In $\tilde A(Z)$, 
$$\eqalign{~
&\text{ch}_{\text{BC}}(E,\rho;g\circ f,\tau_{g\circ
f})\cr
&\qquad-g_*\Big(\text{ch}_{\text{BC}}(E,\rho;f,\tau_f)\cdot
\text{td}(T_g,\tau_g)\Big)
-\text{ch}_{\text{BC}}\Big(f_*E,L^2(\rho,\tau_f);g,\tau_g\Big)\cr
=&(g\circ f)_*\Big(\text{ch}(E,\rho)\cdot
\text{td}_{\text{BC}}(T.,\tau.)\Big)-
\text{ch}_{\text{BC}}\Big((g\circ f)_*E;L^2(\rho,\tau_{g\circ
f}),L^2(L^2(\rho,\tau_f),\tau_g)\Big).\cr}$$
\vskip 0.20cm
\noindent
(E.4) To state the last axiom, the deformation to the normal cone
rule, we need to make a lengthy discussion. So we delay it
untill the next chapter.
\vfill\eject
\vskip 0.20cm
\noindent
{\bf 3. Deformation to the Normal Cone Rule}
\vskip 0.20cm
\noindent
A) Deformation to the normal cone: an algebraic construction
\vskip 0.20cm
\noindent
B) Deformation to the normal cone: the associated metrics
\vskip 0.20cm
\noindent
C) Deformation to the normal cone rule
\vskip 0.20cm
\noindent
In this chapter, we give the final axiom for the relative
Bott-Chern secondary characteristic classes, the deformation to
the normal cone rule.
\vskip 0.20cm
\noindent
{\bf A. Deformation to the normal cone: an algebraic construction}
\vskip 0.20cm
\noindent
(A.1) Let $f:X\to Y$ and $g: Z\to Y$ be two smooth, proper
morphisms of compact K\"ahler manifolds and $i:X\to Z$ be a
codimension one closed immersion over $Y$, i.e., $i$ is a closed
immersion of codinemsion one such that $f=g\circ i$. Then we have
the following standard construction of the deformation to the
normal cone.
\vskip 0.20cm
Denote by $$\pi: W:=B_{X\times \{\infty\}}Z\times {\Bbb P}^1
\rightarrow Z\times {\Bbb P}^1,$$ the natural projection, where
$B_{X\times \{\infty\}}Z\times {\Bbb P}^1$ denotes the blowing-up 
of $Z\times {\Bbb P}^1$ along $X\times \{\infty\}$.  Denote the
exceptional divisor  of $\pi$ by ${\Bbb P}$. It is well-known that
the map $q_W:W\rightarrow {\Bbb P}^1$, obtained by composing $\pi$
with the projection
$q:Z\times {\Bbb P}^1\rightarrow {\Bbb P}^1$, is flat, and that 
for $t\in {\Bbb P}^1$:
$$q_W^{-1}(t)=\cases Z\times \{t\}, &$for $t\not=\infty$$,\\
{\Bbb P}\cup B_XZ,&$for $t=\infty$$.\endcases$$ Here $B_XZ$
denotes the blowing-up of
$Z$ along
$X$. Moreover, by the construction, ${\Bbb P}$ and $B_XZ$
intersect transversally, and 
${\Bbb P}\cap B_XZ$ is the exceptional divisor $X$ on
$B_XZ(=Z\times \{\infty\}$, due to the dimensional reason).
\vskip 0.20cm
Denote by $I:X\times {\Bbb P}^1\hookrightarrow W$  the induced
codimension one closed embedding. Easily we see that the image of
$I$ does not intersect with $B_XZ$, and the image $X\times
\{\infty\}$ in $W$ is a section of ${\Bbb P}$.
Denote  the induced fibration $Z\times \{t\}\to Y\times\{t\}$
by $g_t$ for $t\not=\infty$ and  set $g_\infty$ to be the
composition of the projection of ${\Bbb P}$ on $X$ with
$(X=)X\times\{\infty\}\to Y\times\{\infty\}(=Y)$. Denote by $f_t:
X\times \{t\}\to Y\times \{t\}$ the smooth morphisms induced
from $f$ for all
$t\in {\Bbb P}^1$, by $i_\infty'$ the inclusion of
$W_\infty={\Bbb P}+B_XZ$ into $W$, and by
$k$ (resp. $l$) the inclusion of ${\Bbb P}$ (resp. $B_XZ$) to
$W$.  So
we have the following commutative diagram;
$$\matrix X\times \{t\}&&\hookrightarrow&&
X\times {\Bbb P}^1&&\hookleftarrow&&X\times\{\infty\}\\
|&&&&&&&&|\\
|&i_t\searrow&&&I\downarrow&&&\swarrow i_\infty&|\\
|&&&&&k\swarrow&{\Bbb P}&&|\\
|&(t\not=\infty)&Z\times \{t\}&\hookrightarrow
&B_{X\times\{\infty\}}Z\times {\Bbb P}^1&\buildrel
i_\infty'\over\hookleftarrow&+&\searrow g_\infty&|\\
|&&&&&l\nwarrow&Z&&|\\
f_t&&&\searrow&\pi\downarrow&\swarrow&&&f_\infty\\ |&&&&&&&&|\\
|&g_t\swarrow&&& Z\times {\Bbb P}^1&&&\searrow&|\\
|&&&&&&&&|\\
|&&&&\downarrow&&&&|\\
\downarrow&&&&&&&&\downarrow\\
Y\times\{t\}&&\hookrightarrow&&Y\times{\Bbb
P}^1&&\hookleftarrow&& Y\times\{\infty\}\endmatrix$$
\vskip 0.20cm
\noindent
(A.2) Now let $E$ be a $g$-acyclic vector bundle on $Z$ such that
 $E(-X):=E\otimes {\Cal O}_Z(-X)$ is $g$-acyclic as well. In
particular, we have the short exact sequence
$$0\to E(-X)\to E\to i_*i^*E\to 0\eqno(3.1)$$ on $Z$. We want to
introduce a
${\Bbb P}^1$-deformation for (3.1) via $W\to {\Bbb
P}^1$.
\vskip 0.20cm
First, pulling it back to
$Z\times {\Bbb P}^1$, we have the short exact sequence
$0\to p_Z^*E(-X\times {\Bbb P}^1)\to p_Z^*E\to I'_*(I')^*E\to 0$
on $Z\times {\Bbb P}^1$. Here $p_Z:Z\times {\Bbb P}^1\to Z$
denotes the natural projection and $I':X\times {\Bbb
P}^1\hookrightarrow Z\times {\Bbb P}^1$ denotes the induced
injection from $i$. 
\vskip 0.20cm
To  get a ${\Bbb P}^1$-deformation for  (3.1) on  $W$, 
we  should take, instead of simply taking
$\pi^*(p_Z^*E(-X\times {\Bbb P}^1))$, the vector bundle
$\pi^*(p_Z^*E(Z\times\{\infty\}-X\times {\Bbb P}^1))$. 
(For a possible motivation, see e.g., the definition $DE_1$ in
1.B.1.) Then, we know that
$$\eqalign{~&\pi^*(p_Z^*E(Z\times\{\infty\}-X\times {\Bbb
P}^1))\cr
=&(\pi\circ p_Z)^*E\otimes \pi^*{\Cal O}_{Z\times {\Bbb
P}^1}(Z\times\{\infty\}-X\times {\Bbb P}^1)\cr
=&(\pi\circ
p_Z)^*E\otimes {\Cal O}_{W}({\Bbb P}+B_XZ-X\times
{\Bbb P}^1-{\Bbb P})\cr
=&(\pi\circ
p_Z)^*E\otimes {\Cal O}_{W}(B_XZ-X\times
{\Bbb P}^1)\cr
=&(\pi\circ
p_Z)^*E(B_XZ-X\times
{\Bbb P}^1).\cr}$$ So the correct choice of the ${\Bbb
P}^1$-deformation of (3.1) on $W$ we seek  is the following
exact sequence of coherent sheaves on $W$;
$$0\to (\pi\circ
p_Z)^*E(B_XZ-X\times
{\Bbb P}^1)\to (\pi\circ
p_Z)^*E(B_XZ)\to I_*I^* ((\pi\circ
p_Z)^*E(B_XZ))\to 0.\eqno(3.2)$$
\vskip 0.20cm  
Indeed, by the flatness of $q_W:W\to {\Bbb P}^1$, 
the restrictions of (3.2) to the fibers $W_t$ of $q_W$ for all
$t\in {\Bbb P}^1$ are  exact.
 Thus, in particular, for each
$t\not=\infty$ in
${\Bbb P}^1$, from (3.2), we have the induced
exact sequence
$$0\to E_t(-X)\to E_t\to (i_t)_*i_t^*E_t\to 0
\qquad\text{over}\quad Z\times\{t\}\eqno(3.3)$$ as
$(\pi\circ
p_Z)^*E(B_XZ)\big|_{Z\times\{t\}}=E_t$,
$(\pi\circ
p_Z)^*E(B_XZ-X\times {\Bbb
P}^1)\big|_{Z\times\{t\}}=E_t(-X)$
 with $E_t$ the pull-back of
$E$ under the canonical identity $X\times \{t\}\simeq X$.
Similarly, for the fiber at $\infty$, if we set
$E(B_XZ)\big|_{{\Bbb P}}:=E_\infty$, 
$E(B_XZ)\big|_{B_XZ}=:E_\infty'$ and 
$E(B_XZ-X\times {\Bbb
P}^1)\big|_{B_XZ}=:E_\infty''$. Then
$E(B_XZ-X\times {\Bbb P}^1)\big|_{{\Bbb P}}=E_\infty(-X)$, and
(3.2) splits into two exact sequences 
$$0\to E_\infty(-X)\to E_\infty\to
(i_\infty)_*i_\infty^*E_\infty\to 0\qquad\text{over}\quad
{\Bbb P}\eqno(3.4)$$ and $$0\to E_\infty'
\to E_\infty''\to 0\to 0 \qquad\text{over}\quad B_XZ.\eqno(3.5)$$
Here in the last statement, we use that fact that $I(X\times {\Bbb
P})$ is away from $B_XZ$ in $W$. Thus, in particular, 
on $B_XZ$, $E_\infty'=E_\infty''$. Now we further assume that $E$
on $Z$ is chosen so that
$E_\infty$ and $E_\infty(-X)$ are all $g_\infty$-acyclic.
\vskip 0.20cm
Note that $X\times {\Bbb P}^1$ is away from
$B_XZ$, the restriction of $E_\infty-E_\infty(-X)$ to $B_XZ$ gives
a zero element in the
$K$-group of ${\Bbb P}\cap B_XZ$. As a direct consequence, 
in $K(Y)$, we get, for all $t\in {\Bbb P}^1$,
$$(g_t)_*(E_t)-(g_t)_*(E_t(-X))=(g_\infty)_*(E_\infty)-
(g_\infty)_*(E_\infty(-X)).\eqno(3.6)$$
\vskip 0.20cm
\noindent
{\bf B. Deformation to the normal cone: the associated metrics}
\vskip 0.20cm
\noindent
(B.1)  With the same notation as in the previous section, choose a
K\"ahler metric $\tau_W$ on $W$. Then $\tau_W$ naturally induces
K\"ahler metrics $\tau_{g_t}$ and $\tau_{g_\infty}$ on
$Z_t:=Z\times\{t\}$ and ${\Bbb P}$ respectively for all
$t\not=\infty$. Note that for a smooth morphism, the relative
tangent bundle is a subbundle of the tangent bundle of the total
space,  by taking restriction again, we then get hermitian metrics
$\tau_t$ on $T_{g_t}$ for all $t\in {\Bbb P}^1$. Easily, we see
that the induced metrics for  fibers of $g_t$ from $\tau_t$ are
all K\"ahler. 
\vskip 0.20cm
Now let $E$ be a $g$-acyclic vector bundles on $Z$ such that
$E_t$ and $E_t(-X)$ are  $g_t$-acyclic as
well for all $t\in {\Bbb P}^1$. Fix hermitian metrics $\rho$
and $\rho'$ on $E$ and on $E(-X)$ respectively. Use the same
notation to denote the pull-back of $(E,\rho)$ onto $W$. 
Choose the Fubini-Study metric on
${\Cal O}_{{\Bbb P}^1}(\infty)$ and a metric on ${\Cal
O}_{W}(-X\times {\Bbb P}^1)$ such that in a neighborhood $U$ of
$B_XZ$, which is away from
$X\times {\Bbb P}^1$, the natural isomorphism ${\Cal
O}_{W}(-X\times {\Bbb P}^1)\simeq {\Cal O}_W$ induces an
isometry, once we put the trivial metric on ${\Cal O}_W$. Denote
these final induced metrics on
$E(B_XZ)$ and $E(B_XZ-X\times {\Bbb P}^1)$ by $D\rho$ and
$D\rho'$ respectively. (Question: Figure out how the metric on
${\Cal O}_W(B_XZ)$ is defined?)
\vskip 0.20cm
Denote the induced metrics via restrction to $E_t$ and $E_t(-X)$
(resp. to $E_\infty'$ and $E_\infty''$ on $B_XE$) by $\rho_t$ and
$\rho_t'$ respectively for all
$t\in {\Bbb P}^1$, (resp. $\rho_\infty''$ and $\rho_\infty'''$).
Easily, we see that $\rho_0=\rho$ and $\rho_0'=\rho'$, and 
$(E_\infty',\rho_\infty'')$ is isomorphic to
$(E_\infty'',\rho_\infty''')$ by the
construction.  In this way, $$(E(-X),\rho')\hookrightarrow
(E,\rho)\qquad\text{on}\quad Z$$ is deformed to
$$(E_\infty(-X),\rho_\infty')\hookrightarrow
(E_\infty,\rho_\infty)\quad \text{on}\quad {\Bbb P}$$ (and 
$$(E_\infty',\rho_\infty'')\simeq(E_\infty'',\rho_\infty''')\qquad
\text{on}\quad
B_XZ).$$ 
\vskip 0.20cm
In particular, we have, for all $t\in {\Bbb P}^1$,
\vskip 0.20cm
\noindent
(i) smooth, proper morphisms $g_t$ for compact K\"ahler
manifolds, together with hermitian metrics $\tau_t$ on relative
tangent bundles $T_t:=T_{g_t}$, whose induced metrics on all
fibers of $g_t$ are K\"ahler as well;
\vskip 0.20cm
\noindent
(ii) $g_t$-acyclic hermitian vector bundles $(E_t,\rho_t)$ and 
$(E_t(-X),\rho_t')$.
\vskip 0.20cm
Thus, as a direct consequence, we have the properly metrized
data $(E_t,\rho_t;g_t,\tau_{g_t})$ and
$(E_t(-X),\rho_t';g_t,\tau_{g_t})$ for all $t\in {\Bbb P}^1$, and
their associated relative Bott-Chern secondary characteristic
classes 
$\text{ch}_{\text{BC}}(E_t,\rho_t;g_t,\tau_{g_t})\in \tilde A(Y)$
and 
$\text{ch}_{\text{BC}}(E_t(-X),\rho_t';g_t,\tau_{g_t})\in
\tilde A(Y)$ for all $t\in {\Bbb P}^1$. At this point, we may say
that the final axiom for the relative Bott-Chern secondary
characteristic classes is the one to understand the relation
between
$$\text{ch}_{\text{BC}}(E_0,\rho_0;g_0,\tau_{g_0})-\text{ch}_{\text{BC}}(E_0(-X),\rho_0';g_0,\tau_{g_0})$$
and
$$\text{ch}_{\text{BC}}(E_\infty,\rho_\infty;
g_\infty,\tau_{g_\infty})
-\text{ch}_{\text{BC}}(E_\infty(-X),\rho_\infty';
g_\infty,\tau_{g_\infty})$$ in $\tilde A(Y).$ 
\vskip 0.20cm
\noindent
(B.2) To facilitate the ensuring discussion, we now find out how
$\text{ch}_{\text{BC}}(E_t,\rho_t;g_t,\tau_{g_t})\in\tilde A(Y)$
and 
$\text{ch}_{\text{BC}}(E_t(-X),\rho_t';g_t,\tau_{g_t})\in
\tilde A(Y)$ change with respect to  $t\in {\Bbb
A}^1$.
\vskip 0.20cm
By the construction, for $t\not=\infty$, we may simply view
$g_t:Z\times\{t\}\to Y\times \{t\}$ as the original $g:Z\to Y$.
In this way,
$\text{ch}_{\text{BC}}(E_t,\rho_t;g_t,\tau_{g_t})\in \tilde A(Y)$
may be interpretated as the relative Bott-Chern secondary
characteristic classes for the original morphism $g$ but with
metrics $\rho_t$ on $E$ and $\tau_t$ on $T_g$. Moreover, we may
get this type of changing metrics in two steps: 
First, from
$(\rho_0,\tau_0)$ to $(\rho_0,\tau_t)$; then from
$(\rho_0,\tau_t)$ to $(\rho_t,\tau_t)$.
In this way,  by applying Axiom 5 in step 1 and Axiom 4 in step 2,
we have the following relations:
$$\eqalign{~&\text{ch}_{\text{BC}}(E_0,\rho_0;g_t,\tau_{g_t})
-\text{ch}_{\text{BC}}(E_0,\rho_0;g_0,\tau_{g_0})\cr
=&g_*\Big(\text{ch}(E;\rho_0)\cdot
\text{td}_{\text{BC}}(T_g;\tau_t,\tau_0)\Big)-\text{ch}_{\text{BC}}
\Big(g_*E;L^2(\rho_0,\tau_t),L^2(\rho_0,\tau_0)\Big),\cr}$$
and $$\eqalign{~&\text{ch}_{\text{BC}}(E_t,\rho_t;g_t,\tau_{g_t})
-\text{ch}_{\text{BC}}(E_0,\rho_0;g_t,\tau_{g_t})\cr
=&g_*\Big(\text{ch}_{\text{BC}}(E;\rho_t,\rho_0)\cdot
\text{td}(T_g,\tau_t)\Big)-\text{ch}_{\text{BC}}
\Big(g_*E;L^2(\rho_t,\tau_t),L^2(\rho_0,\tau_t)\Big).\cr}$$
That is to say,
$$\eqalign{~&\text{ch}_{\text{BC}}(E_t,\rho_t;g_t,\tau_{g_t})
-\text{ch}_{\text{BC}}(E_0,\rho_0;g_0,\tau_{g_0})\cr
=&g_*\Big(\text{ch}_{\text{BC}}(E;\rho_t,\rho_0)\cdot
\text{td}(T_g,\tau_t)+\text{ch}(E;\rho_0)\cdot
\text{td}_{\text{BC}}(T_g;\tau_t,\tau_0)\Big)\cr
&-\text{ch}_{\text{BC}}
\Big(g_*E;L^2(\rho_t,\tau_t),L^2(\rho_0,\tau_0)\Big).\cr}\eqno(3.7)$$
Similarly,
$$\eqalign{~&\text{ch}_{\text{BC}}(E_t(-X),\rho_t';g_t,\tau_{g_t})
-\text{ch}_{\text{BC}}(E_0(-X),\rho_0';g_0,\tau_{g_0})\cr
=&g_*\Big(\text{ch}_{\text{BC}}(E(-X);\rho_t',\rho_0')\cdot
\text{td}(T_g,\tau_t)+\text{ch}(E(-X);\rho_0')\cdot
\text{td}_{\text{BC}}(T_g;\tau_t,\tau_0)\Big)\cr
&-\text{ch}_{\text{BC}}
\Big(g_*(E(-X));L^2(\rho_t',\tau_t),L^2(\rho_0',\tau_0)\Big).\cr}
\eqno(3.8)$$
Therefore, 
$$\eqalign{~&\text{ch}_{\text{BC}}(E_t,\rho_t;g_t,\tau_{g_t})
-\text{ch}_{\text{BC}}(E_t(-X),\rho_t';g_t,\tau_{g_t})\cr
&\qquad-\Big(\text{ch}_{\text{BC}}(E_0,\rho_0;g_0,\tau_{g_0})
-\text{ch}_{\text{BC}}(E_0(-X),\rho_0';g_0,\tau_{g_0})\Big)\cr
=&g_*\Big((\text{ch}_{\text{BC}}(E;\rho_t,\rho_0)-\text{ch}_{\text{BC}}(E(-X);\rho_t',\rho_0'))
\cdot
\text{td}(T_g,\tau_t)\cr
&\qquad+(\text{ch}(E;\rho_0)-\text{ch}(E(-X);\rho_0'))\cdot
\text{td}_{\text{BC}}(T_g;\tau_t,\tau_0)\Big)\cr
&-\Big(\text{ch}_{\text{BC}}
\Big(g_*E;L^2(\rho_t,\tau_t),L^2(\rho_0,\tau_0)\Big)-\text{ch}_{\text{BC}}
\Big(g_*(E(-X));L^2(\rho_t',\tau_t),L^2(\rho_0',\tau_0)\Big)\Big).\cr}\eqno(3.9)$$
So if $q_W:W\to {\Bbb P}^1$ was smooth, we would then have found
that the relation between
$$\text{ch}_{\text{BC}}(E_\infty,\rho_\infty;g_\infty,\tau_{g_\infty})
-\text{ch}_{\text{BC}}(E_\infty(-X),\rho_\infty';g_\infty,\tau_{g_\infty})$$
and $$\text{ch}_{\text{BC}}(E_0,\rho_0;g_0,\tau_{g_0})
-\text{ch}_{\text{BC}}(E_0(-X),\rho_0';g_0,\tau_{g_0})$$ can be
deduced from Axioms 4 and 5. 
\vskip 0.20cm
\noindent
{\bf C. Deformation to the normal rule}
\vskip 0.20cm
\noindent
(C.1) Unfortunately, as indicated before, $q_W:W\to {\Bbb P}^1$ is
far from being smooth. So we cannot simply deduce a relation
between the relative Bott-Chern secondary characteristic classes
for the fibers at 0 and $\infty$ directly by only using previous
5 axioms for  relative Bott-Chern secondary characteristic
classes with respect to smooth morphisms. Nevertheless, it seems
most likely that with a suitable generalization of relative
Bott-Chern classes, say for proper and flat morphisms, one then
can understand how the relative Bott-Chern classes change even
across $\infty$. All this is indeed the motivation for
introducing the so-called ternary characteristic classes in
[We2], which themselves need an additional work to justify. So
here, we take an alternative approach by sticking on smooth,
proper morphisms. 
\vskip 0.20cm
\noindent
(C.2) The idea is also quite simple. That is,
instead of trying to find out what is exactly the difference
between
$$\text{ch}_{\text{BC}}(E_\infty,\rho_\infty;g_\infty,\tau_{g_\infty})
-\text{ch}_{\text{BC}}(E_\infty(-X),\rho_\infty';g_\infty,\tau_{g_\infty})$$
and $$\text{ch}_{\text{BC}}(E_0,\rho_0;g_0,\tau_{g_0})
-\text{ch}_{\text{BC}}(E_0(-X),\rho_0';g_0,\tau_{g_0}),$$
we view $\{\text{ch}_{\text{BC}}(E_t,\rho_t;g_t,\tau_{g_t})
-\text{ch}_{\text{BC}}(E_t(-X),\rho_t';g_t,\tau_{g_t})\}_{t\in
{\Bbb P}^1}$ as a family on ${\Bbb P}^1$ and try to establish the
continuity of such a family of forms, or better of classes.
\vskip 0.20cm
It looks that we are now in a position to state the final
axiom.  No, it is not so. There is yet another
difficulty: It is well-known that for the family $W\to
Y\times {\Bbb P}^1$, usually the induced
$L^2$-metrics on the direct images for $g_t$ when $t\to\infty$
have singularities. So we need to take care of such 
singularities as well. (Recall that by Axiom 1 in 2.A.2 for
relative Bott-Chern secondary characteristic classes,
$$dd^c\text{ch}_{\text{BC}}(E_t,\rho_t;g_t,\tau_{g_t})
=(g_t)_*\Big(\text{ch}(E_t,\rho_t)\cdot\text{td}(T_t,\tau_{g_t})\Big)
-\text{ch}\Big((g_t)_*E,L^2(\rho_t,\tau_{g_t})\Big),$$ and
$$\eqalign{~&dd^c\text{ch}_{\text{BC}}(E_t(-X),\rho_t';g_t,\tau_{g_t})\cr
&\qquad=(g_t)_*\Big(\text{ch}(E_t(-X),\rho_t)\cdot\text{td}(T_t,\tau_{g_t})\Big)
-\text{ch}\Big((g_t)_*E(-X),L^2(\rho_t',\tau_{g_t})\Big).)\cr}$$
 For this
purpose,  consider the virtual vector bundle $(g_t)_*(E_t)-
(g_t)_*(E_t(-X))$ on $Y$ for
$t\in {\Bbb P}^1$. By the discussion in (B.1), in particular,
(3.6), we know that for
$t\not=\infty$,
$(g_t)_*(E_t)=(g_0)_*(E_0)$, $(g_t)_*(E_t(-X))=(g_0)_*(E_0(-X))$
once we identify $Y\times \{0\}$ with $Y\times\{t\}$, and
as virtual vector bundles on $Y$,
$$(g_\infty)_*(E_\infty)- (g_\infty)_*(E_\infty(-X))=
(g_t)_*(E_t)-(g_t)_*(E_t(-X)).$$ So in particular, there exists a
natural number $r$ such that, locally over $Y$, as holomorphic
vector bundles, we have  the local isomorphism
$${\Cal O}_Y^{\oplus r}\oplus(g_\infty)_*(E_\infty)\oplus
(g_t)_*(E_t(-X))=
{\Cal O}_Y^{\oplus r}\oplus(g_t)_*(E_t)\oplus
(g_\infty)_*(E_\infty(-X)).$$ Hence, by the fact that Chern form
depends only on the $C^\infty$-structure of a vector bundle,
we may choose hermitian metrics $\gamma_t$ and $\gamma_t'$ on
$(g_t)_*E$ and $(g_t)_*(E_t(-X))$ respectively for all $t\in {\Bbb
P}^1$ such that
$$\text{ch}(E_t,\gamma_t)-
\text{ch}(E_t(-X),\gamma_t')=\text{ch}(E_\infty,\gamma_\infty)-
\text{ch}(E_\infty(-X),\gamma_\infty').\eqno(3.10)$$
With this, via $Y\times\{0\}=Y\times \{\infty\}$, we have the
 classical Bott-Chern secondary characteristic classes
$\text{ch}_{\text{BC}}((g_t)_*(E_t;
L^2(\rho_t,\tau_t),\gamma_t)$ and
$\text{ch}_{\text{BC}}((g_t)_*(E_t(-X));
L^2(\rho_t',\tau_t),\gamma_t')$.
Now we may consider the family
$$\eqalign{~&\{\Big(\text{ch}_{\text{BC}}(E_t,\rho_t;g_t,\tau_{g_t})
-\text{ch}_{\text{BC}}(E_t(-X),\rho_t';g_t,\tau_{g_t})\Big)\cr
&+\Big(\text{ch}_{\text{BC}}((g_t)_*(E_t);
L^2(\rho_t,\tau_t),\gamma_t)-\text{ch}_{\text{BC}}((g_t)_*(E_t(-X));
L^2(\rho_t',\tau_t),\gamma_t')\Big)\}_{t\in {\Bbb
P}^1}.\cr}$$ As a matter of fact,  the final axiom for the
relative Bott-Chern secondary characteristic classes requires
that this family is  continuous. 
\vskip 0.20cm
\noindent
(C.3) {\it Axiom 6.} ({\bf Deformation to the normal cone rule})
In $\tilde A(Y)$,
$$\eqalign{~&\lim_{t\to
\infty}\Big(\big(\text{ch}_{\text{BC}}(E_t,\rho_t;g_t,\tau_{g_t})
-\text{ch}_{\text{BC}}(E_t(-X),\rho_t';g_t,\tau_{g_t})\big)\cr
&+\big(\text{ch}_{\text{BC}}((g_t)_*(E_t);
L^2(\rho_t,\tau_t),\gamma_t)-\text{ch}_{\text{BC}}((g_t)_*(E_t(-X));
L^2(\rho_t',\tau_t),\gamma_t')\big)\Big)\cr
=&\big(\text{ch}_{\text{BC}}(E_\infty,\rho_\infty;g_\infty,\tau_{g_\infty})-
\text{ch}_{\text{BC}}(E_\infty(-X),\rho_\infty';g_\infty,\tau_{g_\infty})\big)+\cr
&\big(\text{ch}_{\text{BC}}((g_\infty)_*(E_\infty);
L^2(\rho_\infty,\tau_\infty),\gamma_\infty)
-\text{ch}_{\text{BC}}((g_\infty)_*(E_\infty(-X));
L^2(\rho_\infty',\tau_\infty),\gamma_\infty')\big).\cr}$$
\vskip 0.20cm
\noindent
(C.4) All this then gives us the six key axioms for the
relative Bott-Chern secondary characteristic classes. The reader
may see that Axioms 2 and 3 (resp. Axioms 4 and 5) may be grouped
together naturally. So there are essentially four groups of
axioms. The first group, the downstairs rule, relates the
secondary theory with its origin; the second group tells us how
to deal with the changes from the base; the third group tells us
how to deal with the changes from the total space; while the last
tells us how to relate the theory for different relative
dimensions naturally.
\vfill\eject
\vskip 0.20cm
\noindent
{\bf 4. Uniqueness of Relative Bott-Chern Secondary
Characteristic  Classes}
\vskip 0.20cm
\noindent
A)  Weak Uniqueness Theorem 
\vskip 0.20cm
\noindent
B)  Strong Uniqueness Theorem 
\vskip 0.20cm
\noindent
In this chapter, we will state the  uniqueness
theorem for  relative Bott-Chern secondary characteristic
classes. There are two versions of it.
 Roughly speaking, the first uniqueness
theorem claims that after modulo
$d$-closed forms, the classes of the relative Bott-Chern
secondary characteristic classes are unique, while the second
uniqueness theorem claims that after modulo exact forms, the relative Bott-Chern
secondary characteristic classes are essentially unique, i.e.,
the difference of any two possible relative Bott-Chern
secondary characteristic classes can be understood precisely.
\vskip 0.20cm
\noindent
{\bf A.  Weak Uniqueness Theorem} 
\vskip 0.20cm
\noindent
(A.1) We have already stated the axioms for relative Bott-Chern
secondary characteristic classes associated to smooth, proper
metrized morphisms and relative acyclic hermitian vector
bundles in Chapter 2. Now we state the first uniqueness result
for them. Recall that
$(f:X\to Y; E,\rho;T_f,\tau_f)$ is called a properly metrized
datum  if $f:X\to Y$ is a smooth morphism of compact
K\"ahler manifolds,
$(E,\rho)$ is an $f$-acyclic hermitian vector bundle, and  
$\tau_f$ is a hermitian metric on the relative tangent bundle
$T_f$ of $f$ such that the induced metrics on all fibers of $f$
are K\"ahler. 
\vskip 0.20cm
\noindent
(A.2) {\bf Theorem.} ({\bf The Weak Uniqueness for
Relative Bott-Chern Secondary Characteristic Classes}) {\it
Suppose that  there are two constructions
$\text{ch}_{\text{BC}}$ and
$\text{ch}_{\text{BC}}'$ which satisfy six axioms for 
relative Bott-Chern secondary characteristic classes, then, for
all properly metrized data $(f:X\to Y;E,\rho;T_f,\tau_f)$, the
images of $\text{ch}_{\text{BC}}(E,\rho;f,\tau_f)$ and
$\text{ch}_{\text{BC}}'(E,\rho;f,\tau_f)$ in $\hat
A(Y):=A(Y)/(\text{exact\ forms}+ d-\text{closed\ forms})$ are the
same.} 
\vskip 0.20cm
\noindent
{\bf B.  Strong Uniqueness Theorem}
\vskip 0.20cm
\noindent
(B.1) To state the second uniqueness theorem, we
need  a preparation.
\vskip 0.20cm
Let $B$ be a subring of ${\Bbb R}$, and let $P(x)\in B[[x]]$ be
any symmetric power series. Then for any  vector bundle
$E$ on $Y$, by the splitting principle for vector bundles, there
exists a unique additive characteristic class $P({E})\in
H^{*}(X,{\Bbb R})$, the de Rham cohomology of $X$. 
\vskip 0.20cm
Now assume  that there exists a
construction $\text{ch}_{\text{BC}}$ such that it satisfies the
six axioms for relative Bott-Chern secondary characteristic
classes. Then, for any properly metrized datum $(f:X\to
Y;E,\rho;T_f,\tau_f)$, we have an element 
$\text{ch}_{\text{BC}}(E,\rho;f,\tau_f)\in
\tilde A(Y)$.
\vskip 0.20cm
With such an existence,  then we claim that
$\text{ch}_{\text{BC}}(E,\rho;f,\tau_f)+f_*\Big(\text{ch}(E)\cdot
\text{td}(T_f)\cdot P(T_f)\Big)\in \tilde A(Y)$ for any fixed
additive characteristic class $P$ as above also
satisfies the listed six axioms  for relative Bott-Chern
secondary characteristic classes. Indeed, easily, by the fact that
$f_*\Big(\text{ch}(E)\cdot
\text{td}(T_f)\cdot P(T_f)\Big)$ is in $H^*(Y,{\Bbb R})$, we see
that
$\text{ch}_{\text{BC}}(E,\rho;f,\tau_f)+f_*\Big(\text{ch}(E)\cdot
\text{td}(T_f)\cdot P(T_f)\Big)\in \tilde A(Y)$ satisfies the
first five axioms. Moreover note that the deformation to the
normal cone does not change the cohomology classes, which is the
key property used in the proof of Grothendieck-Riemann-Roch
theorem in algebraic geometry, we equally can check the axiom 6
for $\text{ch}_{\text{BC}}(E,\rho;f,\tau_f)+f_*\Big(\text{ch}(E)\cdot
\text{td}(T_f)\cdot P(T_f)\Big)$. (See e.g., Lemma 5.B.7.4 below.)
In this sense, we may say that there is no uniqueness for relative
Bott-Chern secondary characteristic classes.  
\vskip 0.20cm
Even the latest statement is absolutely correct, we still can
understand precisely the structure of relative Bott-Chern
secondary characteristic classes. Roughly speaking, we may say
that the twisting by an additive characteristic class stated
above is the only flaw for establishing the uniqueness for
relative Bott-Chern secondary characteristic classes. 
\vskip 0.20cm
\noindent
(B.2) {\bf Theorem.} ({\bf The Srong Uniqueness for
Relative Bott-Chern Secondary Characteristic Classes}) {\it 
Suppose that there are two constructions
$\text{ch}_{\text{BC}}$ and
$\text{ch}_{\text{BC}}'$ which satisfy six axioms for 
relative Bott-Chern secondary characteristic classes, then, 
there exists an additive characteristic class $R$
such that,  for all properly metrized data $(f:X\to
Y;E,\rho;T_f,\tau_f)$, in $\tilde A(Y)$,
$$\text{ch}_{\text{BC}}'(E,\rho;f,\tau_f)=
\text{ch}_{\text{BC}}(E,\rho;f,\tau_f)+f_*\Big(\text{ch}(E)\cdot
\text{td}(T_f)\cdot R(T_f)\Big).$$}
\vskip 0.20cm
\noindent
{\it Remark 4.1.} The weak uniqueness theorem says that the
relative Bott-Chern  classes are unique
modulo $d$-closed forms, while the strong uniqueness theorem says
that the relative Bott-Chern  classes are
unique up to a certain well-structured topological term.
\vskip 0.20cm
Obviously, the weak  uniqueness  is a direct
consequence of the strong uniqueness, so it is
sufficient to prove the later one.
We will do it in the following two chapters.
\vfill\eject
\vskip 0.20cm
\noindent
{\bf 5. Some intermediate results}
\vskip 0.20cm
\noindent
A) Statements of intermediate results
\vskip 0.20cm
\noindent
B) The proofs 
\vskip 0.20cm
\noindent
In this chapter, we prove some intermediate results for relative
Bott-Chern secondary characteristic classes, which will be used
in the proof of the uniqueness theorems. So in Chapter 5 and
Chapter 6, we will assume that there is a construction
$\text{ch}_{\text{BC}}$ which satisfies all six axioms for the
relative Bott-Chern secondary characteristic classes.
\vskip 0.20cm
\noindent
{\bf A. Statements of intermediate results}
\vskip 0.20cm
\noindent
(A.1) Suppose that for any properly metrized datum $(f:X\to
Y;E,\rho;T_f,\tau_f)$, there exist  elements
$\text{ch}_{\text{BC}}(E,\rho;f,\tau_f)\in \tilde A(Y)$ and
$\text{ch}_{\text{BC}}'(E,\rho;f,\tau_f)\in\tilde A(Y)$
which satisfy the  six
axioms stated in Chapters 2 and 3. Let $P$ be an additive
characteristic class and set
$$\text{Err}(E,\rho;f,\tau_f;P):=\text{ch}_{\text{BC}}(E,\rho;f,\tau_f)
-\text{ch}_{\text{BC}}'(E,\rho;f,\tau_f)+f_*\Big(\text{ch}(E)\cdot
\text{td}(T_f)\cdot P(T_f)\Big).$$ We want
 to show that there exists a unique universal additive
characteristic class $R$ such that
$\text{Err}(E,\rho;f,\tau_f;R)=0$ for all properly metrized data
$(E,\rho;f,\tau_f)$. For this purpose, let
us state some intermediate results.
\vskip 0.20cm
First we study how Err depends on hermitian metrics on
vector bundles.
\vskip 0.20cm
\noindent {\bf Proposition.} {\it Let $f:X\rightarrow Y$ be a
smooth, proper morphism of smooth complex K\"ahler manifolds. Fix a
hermitian metric $\tau_f$ on the relative  tangent vector bundle
$T_f$ of $f$ such that the induced metrics on all fibers of $f$
are K\"ahler. Then for any short exact sequence of 
$f$-acyclic hermitian vector bundles
$$E.:0\rightarrow E_1\rightarrow E_2\rightarrow E_3\rightarrow
0,$$ with  hermitian metrics $\rho_i$ on $E_i$
for $i=1,2,3,$ 
$$\text{Err}(E_1,\rho_1;f,\rho_f;P)+\text{Err}(E_3,\rho_3;f,\rho_f;P)
=\text{Err}(E_2,\rho_2;f,\rho_f;P).$$ In particular,
$\text{Err}(E,\rho;f,\rho_f;P)$ does not depend on the choice of
the metric
$\rho$. Furthermore, 
$\text{Err}(E,\rho;f,\rho_f;P)$
lies in the classes of $dd^c$-closed  forms.}
\vskip 0.20cm
\noindent
(A.2) Now we study how Err depends on hermitian metrics on
manifolds.
\vskip 0.20cm
\noindent
{\bf Proposition.} {\it Let $f:X\rightarrow Y$ and $g:Y\rightarrow
Z$ be two smooth, proper morphisms of compact K\"ahler manifolds which
admit hermitian metrics $\tau_f$, $\tau_g$ and
$\tau_{g\circ f}$  on the relative tangent vector  bundles of
$f$, $g$ and $g\circ f$ respectively, such that the induced
metrics on all fibers are K\"ahler. Let
$(E,\rho)$ be an
$f$- and $g\circ f$-acyclic hermitian  vector bundle on $X$ such
that $f_*E$ is $g$-acyclic as well. Then
$$\text{Err}(E,\rho;g\circ f,\tau_{g\circ f};P)=
\text{Err}(f_*E,f_*\rho;g,\tau_g;P)
+g_*(\text{Err}(E,\rho;f,\tau_f;P)\cdot
\text{td}(T_g,\tau_g)).$$ In particular,
$\text{Err}(E,\rho;f,\tau_f;P)$ does not depend on the metric
$\tau_f$.} 
\vskip 0.20cm
\noindent
{\it Remark 5.1.} Because of these two propositions, we from now
on in this section denote
$\text{Err}(E,\rho;f,\tau_f;P)$ simply by $\text{Err}(E;f;P)$. 
\vskip 0.20cm
\noindent
(A.3) We may go slight further. Indeed, we have the following
\vskip 0.20cm
\noindent
{\bf Proposition.} {\it  Let $f:X\rightarrow Y$ be a
smooth, proper morphism of compact K\"ahler manifolds. There is a
natural morphism
$$\text{Err}(\cdot;f,P): K(X)\rightarrow H^*(Y,{\Bbb R}),$$ 
such that
$\text{Err}(E;f,P)=\text{Err}(E;f;P)$ for all $f$-acyclic vector
bundles $E$ on $X$.}
\vskip 0.20cm
\noindent
{\it Remark 5.2.} Because of this prosition,  we may equally write
$\text{Err}(\cdot;f,P)$ as $\text{Err}(\cdot;f;P)$.
\vskip 0.20cm
\noindent
(A.4) We have a functorial property as well.
\vskip 0.20cm
\noindent{\bf Proposition.} {\it Let $f:X\rightarrow Y$ be a
smooth, proper morphism of compact K\"ahler manifolds. Then, for 
the Cartesian diagram
  $$\matrix Y'\times_Y X&\buildrel {g_f}\over\rightarrow &X\\
f_g\downarrow &&\downarrow f\\
Y'&\buildrel g\over\rightarrow &Y\endmatrix$$ induced from a
base change $g:Y'\to Y$ of compact K\"ahler manifolds,
 $$g^*\text{Err}(E;f;P)=\text{Err}(g_f^*E;f_g;P).$$} 
\vskip 0.20cm
\noindent
(A.5) We now give a projection formula for Err.
\vskip 0.20cm
\noindent
{\bf Proposition.} {\it Let $f:X\rightarrow Y$ be a
smooth, proper morphism of compact K\"ahler manifolds. Let $E$
and $F$ be vector bundles on $X$ and $Y$ respectively. Then 
$$\text{Err}(E\otimes
f^*F;f;P)=\text{Err}(E;f;P)\cdot\text{ch}(F).$$}
\vskip 0.20cm
\noindent
(A.6) We now study what happens when $f$ is a ${\Bbb P}^1$-bundle
on $Y$.
\vskip 0.20cm
\noindent
{\bf Proposition.} {\it There is a unique characteristic
class $R$ for vector bundles of rank $\leq 2$ such that for any 
${\Bbb P}^1$-bundle $f:X={\Bbb
P}_Y(F)\rightarrow Y,$
$\text{Err}(\cdot;f;R)\equiv 0$.}
\vskip 0.20cm
\noindent
(A.7) Finally, we consider Err for closed immersions.
For doing so, we introduce  a new Err term: Let
$i:X\hookrightarrow Z$ be a closed immersion with the smooth
structure morphisms
$f:X\rightarrow Y$ and $g:Z\rightarrow Y$ of compact K\"ahler
manifolds, set
$$\text{Err}(E;i;P):=\text{Err}(E;f;P)-\text{Err}(i_*E;g;P).$$
By Proposition A.3, this definition makes sense, even though
$i_*E$ is usually  a coherent sheaf only.
\vskip 0.20cm
\noindent
{\bf Proposition.} {\it Let $i:X\hookrightarrow Z$
be a codimension-one regular closed immersion of compact K\"ahler
manifolds over a compact K\"ahler
manifold $Y$ with smooth structure morphisms $f:X\rightarrow Y$
and $g:Z\rightarrow Y.$  Then for any vector bundle $E$ on $X$,
$$\text{Err}(E;i;P)=0,$$ for any
additive characteristic class $P$, which, for rank two
vector bundles, coincides with
$R$ in Proposition A.6.}
\vskip 0.20cm
\noindent
{\bf B. The Proofs}
\vskip 0.20cm
\noindent
(B.1) {\it Proof of Proposition (A.1).}
In fact, in $\tilde A(Y)$, 
$$\eqalign
{~&\text{Err}(E_2,\rho_2;f,\tau_f;P)-
\Big(\text{Err}(E_1,\rho_1;f,\tau_f;P)
+\text{Err}(E_3,\rho_3;f,\tau_f;P)\Big)\cr
=&\text{ch}_{\text{BC}}(E_2,\rho_2;f,\tau_f)
-\text{ch}_{\text{BC}}'(E_2,\rho_2;f,\tau_f)
+f_*\big(\text{ch}(E_2)\cdot
\text{td}(T_f)\cdot P(T_f)\big)\cr
&-\Big(\text{ch}_{\text{BC}}(E_1,\rho_1;f,\tau_f)
-\text{ch}_{\text{BC}}'(E_1,\rho_1;f,\tau_f)
+f_*\big(\text{ch}(E_1)\cdot
\text{td}(T_f)\cdot P(T_f)\big)\cr
&+\text{ch}_{\text{BC}}(E_3,\rho_3;f,\tau_f)
-\text{ch}_{\text{BC}}'(E_3,\rho_3;f,\tau_f)
+f_*\big(\text{ch}(E_3)\cdot
\text{td}(T_f)\cdot P(T_f)\big)\Big)\cr
&\qquad(\text{by\ definition})\cr}$$
$$\eqalign{=&\text{ch}_{\text{BC}}(E_2,\rho_2;f,\tau_f)
-\Big(\text{ch}_{\text{BC}}(E_1,\rho_1;f,\tau_f)
+\text{ch}_{\text{BC}}(E_3,\rho_3;f,\tau_f)\Big)\cr
&-\Big(\text{ch}_{\text{BC}}'(E_2,\rho_2;f,\tau_f)
-\big(\text{ch}_{\text{BC}}'(E_1,\rho_1;f,\tau_f)
+\text{ch}_{\text{BC}}'(E_3,\rho_3;f,\tau_f)\big)\Big)\cr
&+f_*\Big(\big(\text{ch}(E_2)-\text{ch}(E_1)
-\text{ch}(E_3)\big)\cdot
\text{td}(T_f)\cdot P(T_f)\Big)\cr
=&f_*\Big(\text{ch}_{\text{BC}}(E.,\rho.)
\cdot\text{td}(T_f,\tau_f)\Big)
-\text{ch}_{\text{BC}}(f_*(E.),L^2(\rho.;\tau_f))\cr
&-\Big(f_*\Big(\text{ch}_{\text{BC}}(E.,\rho.)
\cdot\text{td}(T_f,\tau_f)\Big)
-\text{ch}_{\text{BC}}(f_*(E.),L^2(\rho.;\tau_f))\Big)\cr
&+f_*\Big(0\cdot
\text{td}(T_f)\cdot P(T_f)\Big)\qquad(\text{by\ Axiom\
4\ in\ 2.D.2})\cr =&0.\cr}$$
This is simply the first statement of the proposition.
With this, by taking $E_3$ to be zero, we get  the second
statement. Finally,
$$\eqalign{~&dd^c\text{Err}(E,\rho;f,\tau_f;P)\cr
=&dd^c\Big(\text{ch}_{\text{BC}}(E,\rho;f,\tau_f)
-\text{ch}_{\text{BC}}'(E,\rho;f,\tau_f)
+f_*\big(\text{ch}(E_2)\cdot
\text{td}(T_f)\cdot P(T_f)\big)\Big)\cr
&\qquad(\text{by\ definition})\cr
=&dd^c\Big(\text{ch}_{\text{BC}}(E,\rho;f,\tau_f)\Big)
-dd^c\Big(\text{ch}_{\text{BC}}'(E,\rho;f,\tau_f)\Big)
+dd^c\Big(f_*\big(\text{ch}(E_2)\cdot
\text{td}(T_f)\cdot P(T_f)\big)\Big)\cr
=&f_*\big(\text{ch}(E,\rho)\cdot\text{td}(T_f,\tau_f)\big)
-\text{ch}(f_*E,L^2(\rho;\tau_f))\cr
&-\Big(f_*\big(\text{ch}(E,\rho)\cdot\text{td}(T_f,\tau_f)\big)
-\text{ch}(f_*E,L^2(\rho;\tau_f))\Big)+0\quad(\text{by\ Axiom\
1\ in\ 2.A.2})\cr
=&0.\cr}$$ 
This completes the proof of Proposition A.1.
\vskip 0.20cm
\noindent
(B.2) {\it Proof of Proposition A.2.} Setting $g=\text{Id}_Y$,
the identity map of $Y$, we find that the second statement of this
proposition is a consequence of the first one.
\vskip 0.20cm
Now we prove the first statement.
$$\eqalign{~&\text{Err}(E,\rho;g\circ f,\tau_{g\circ
f};P)\cr
&\qquad-\Big(
\text{Err}(f_*E,f_*\rho;g,\tau_g;P)
+g_*(\text{Err}(E,\rho;f,\tau_f;P)
\text{td}(T_g,\tau_g))\Big)\cr
=&\text{ch}_{\text{BC}}(E,\rho;g\circ f,\tau_{g\circ f})
-\text{ch}_{\text{BC}}'(E,\rho;g\circ
f,\tau_{g\circ f})\cr
&\qquad+(g\circ f)_*\Big(\text{ch}(E)\cdot
\text{td}(T_{g\circ f})\cdot P(T_{g\circ f})\Big)\cr
&-\Big(\text{ch}_{\text{BC}}(f_*E,L^2(\rho,\tau_f);g,\tau_g)
-\text{ch}_{\text{BC}}'(f_*E,L^2(\rho,\tau_f);g,\tau_{g})
\cr
&\qquad+f_*\Big(\text{ch}(f_*E)\cdot
\text{td}(T_{g})\cdot P(T_{g})\Big)\Big)\cr
&-g_*\Big(\big(\text{ch}_{\text{BC}}(E,\rho;f,\tau_{f})
-\text{ch}_{\text{BC}}'(E,\rho;f,\tau_{f})\cr
&\qquad
+f_*\Big(\text{ch}(E)\cdot
\text{td}(T_{f})\cdot
P(T_{f})\Big)\big)\cdot\text{td}(T_g,\tau_g)\Big)\quad(\text{by\
definition})\cr}$$ $$\eqalign{
=&\Big(\text{ch}_{\text{BC}}(E,\rho;g\circ f,\tau_{g\circ f})
-\text{ch}_{\text{BC}}(f_*E,L^2(\rho,\tau_f);g,\tau_g)\cr &\qquad
-g_*\Big(\text{ch}_{\text{BC}}(E,\rho;f,\tau_{f})
\cdot\text{td}(T_g,\tau_g)\Big)\cr
&-\Big(\text{ch}_{\text{BC}}'(E,\rho;g\circ
f,\tau_{g\circ f})
-\text{ch}_{\text{BC}}'(f_*E,L^2(\rho,\tau_f);g,\tau_{g})
\cr
&\qquad-g_*\Big(\text{ch}_{\text{BC}}'(E,\rho;f,\tau_{f})\cdot
\text{td}(T_g,\tau_g)\Big)\cr
&+\Big((g\circ f)_*\Big(\text{ch}(E)\cdot
\text{td}(T_{g\circ f})\cdot P(T_{g\circ f})\Big)
-g_*\Big(\text{ch}(f_*E)\cdot
\text{td}(T_{g})\cdot P(T_{g})\Big)\cr
&\qquad
-g_*\Big(f_*\Big(\text{ch}(E)\cdot
\text{td}(T_{f})\cdot
P(T_{f})\Big)\cdot
\text{td}(T_g,\tau_g)\Big)\cr
=&(g\circ f)_*\Big(\text{ch}(E,\rho)\cdot 
\text{td}_{\text{BC}}(T.,\tau.)\Big)\cr
&\qquad-
\text{ch}_{\text{BC}}\Big((g\circ f)_*E;L^2(\rho,\tau_{g\circ
f}),L^2(L^2(\rho,\tau_f),\tau_g)\Big)\cr
&-\Big((g\circ f)_*\Big(\text{ch}(E,\rho)\cdot 
\text{td}_{\text{BC}}(T.,\tau.)\Big)\cr
&\qquad-
\text{ch}_{\text{BC}}\Big((g\circ f)_*E;L^2(\rho,\tau_{g\circ
f}),L^2(L^2(\rho,\tau_f),\tau_g)\Big)\Big)\cr
&+\Big((g\circ f)_*\Big(\text{ch}(E)\cdot
\text{td}(T_{g\circ f})\cdot P(T_{g\circ f})\Big)
-g_*\Big(\text{ch}(f_*E)\cdot
\text{td}(T_{g})\cdot P(T_{g})\Big)\cr
&\qquad
-g_*\Big(f_*\Big(\text{ch}(E)\cdot
\text{td}(T_{f})\cdot
P(T_{f})\Big)\cdot
\text{td}(T_g,\tau_g)\Big)\cr
&\qquad(\text{by\ Axiom\ 4\ in\ 2.D.2})\cr
=&0+(g\circ f)_*\Big(\text{ch}(E)\cdot
\text{td}(T_{g\circ f})\cdot P(T_{g\circ f})\Big)
-g_*\Big(\text{ch}(f_*E)\cdot
\text{td}(T_{g})\cdot P(T_{g})\Big)\cr
&\qquad
-g_*\Big(f_*\Big(\text{ch}(E)\cdot
\text{td}(T_{f})\cdot
P(T_{f})\Big)\cdot
\text{td}(T_g)\Big).\cr}$$ 
From here it is sufficient to prove the following;
\vskip 0.20cm
\noindent
{\bf Lemma.} {\it With the same notation as above, 
$$\eqalign{(g\circ f)_*&\Big(\text{ch}(E)\cdot
\text{td}(T_{g\circ f})\cdot P(T_{g\circ f})\Big)\cr
=&g_*\Big(\text{ch}(f_*E)\cdot
\text{td}(T_{g})\cdot P(T_{g})\Big)
+g_*\Big(f_*\Big(\text{ch}(E)\cdot
\text{td}(T_{f})\cdot
P(T_{f})\Big)\cdot
\text{td}(T_g)\Big).\cr}$$}
\vskip 0.20cm
\noindent
{\it Proof.} By definition, $P$ is additive and
$\text{td}$ is multiplicative. Hence applying them to the exact
sequence of relative tangent bundles
$$0\to T_f\to T_{g\circ f}\to f^*T_g\to 0,$$ we have
$P(T_{g\circ f})=P(T_f)+f^*P(T_g),\quad \text{td}(T_{g\circ
f})=\text{td}(T_f)\cdot f^*\text{td}(T_g).$ Here we use the fact
that as characteristic classes, $P$ and $\text{td}$ have the
functorial property.
Hence 
$$\eqalign{~&(g\circ f)_*\Big(\text{ch}(E)\cdot
\text{td}(T_{g\circ f})\cdot P(T_{g\circ f})\Big)\cr
=&(g\circ f)_*\Big(\text{ch}(E)\cdot
\text{td}(T_f)\cdot f^*\text{td}(T_g)
\cdot \big(P(T_f)+f^*P(T_g)\big)\Big)\cr
=&(g\circ f)_*\Big(\text{ch}(E)\cdot
\text{td}(T_f)\cdot f^*\text{td}(T_g)
\cdot P(T_f)\Big)\cr
&\qquad+(g\circ f)_*\Big(\text{ch}(E)\cdot
\text{td}(T_f)\cdot f^*\text{td}(T_g)
\cdot f^*\big(P(T_g)\big)\Big)\cr
=&g_*\Big(f_*\Big(\text{ch}(E)\cdot
\text{td}(T_f)
\cdot P(T_f)\cdot f^*\text{td}(T_g)\Big)\Big)\cr
&\qquad
+g_*\Big(f_*\Big(\text{ch}(E)\cdot
\text{td}(T_f)\cdot f^*\big(\text{td}(T_g)
\cdot P(T_g)\big)\Big)\Big)\cr
=&g_*\Big(f_*\Big(\text{ch}(E)\cdot
\text{td}(T_f)
\cdot P(T_f)\Big)\cdot \text{td}(T_g)\Big)\cr
&\qquad+
g_*\Big(f_*\Big(\text{ch}(E)\cdot
\text{td}(T_f)\Big)\cdot \text{td}(T_g)
\cdot P(T_g)\Big)\cr
&\qquad(\text{by\ the\ projection\
formula})\cr
=&g_*\Big(f_*\Big(\text{ch}(E)\cdot
\text{td}(T_f)
\cdot P(T_f)\Big)\cdot \text{td}(T_g)\Big)+
g_*\Big(\text{ch}(f_*(E))\cdot \text{td}(T_g)
\cdot P(T_g)\Big)\cr
&\qquad(\text{by\ the\ 
Grothendieck-Riemann-Roch\ Theorem\ in\ Algebraic\
Geometry}).\cr}$$ This completes the proof of the lemma ane hence
Proposition A.2.
\vskip 0.20cm
\noindent
(B.3) {\it Proof of Proposition A.3.} The first point we need to
address is that by Proposition A.1, Err is only
$dd^c$-closed, but in this proposition, we instead choose the de
Rham cohomology groups as the range. This may be solved as
follows.
\vskip 0.20cm
For a compact K\"ahler manifold $Y$, denote the space of
$dd^c$-closed forms by $A^*_{dd^c}(Y)$ and set
$H_{dd^c}^*(Y,{\Bbb R})$ to be the quotient space of
$A^*_{dd^c}(Y)$ modulo the $\partial-$ and $\bar\partial$-forms.
It is well-known that for a compact K\"ahler manifold, a  form
$\omega$ is
$d$-closed and is either $d-$, or $\partial-$, or
$\bar\partial$-exact if and only if there exists a form
$\gamma$ such that $dd^c\gamma=\omega$. (See e.g. 7.A.1 later.)
Hence if $\omega$ is
$dd^c$-closed, then we may choose $\gamma$ to be
$\phi+\partial\alpha+\bar\partial\beta$ with $\phi$ harmonic.
Therefore, we see that $H_{dd^c}^*(Y,{\Bbb R})$ is isomorphic to
the cohomology of $Y$.
\vskip 0.20cm
 The second point here is that we should construct a
map Err from the Grothendieck
$K$-group of coherent sheaves of $X$ to $H_{dd^c}(Y, {\Bbb R})$
which coincides with the previous Err  for $f$-acyclic vector
bundles. But this is rather formal.
\vskip 0.20cm
\noindent
(i) For complex K\"ahler
manifolds, the Grothendieck $K$ group is isomorphic to the
$K$ group generated by vector bundles. So we will write both
of these $K$-groups as $K(X)$; 
\vskip 0.20cm
\noindent
(ii) for $X$, a compact K\"ahler manifold, smooth and
proper over
$Y$,  $K(X)$ is generated by $f$-acyclic vector bundles.
Indeed, for any coherent sheaf $F$ on $X$, there exists an
 vector bundle resolution $0\to F\to E_0\to
E_2\to\dots\to E_m\to 0$ such that all $E_i$'s are $f$-acyclic;
\vskip 0.20cm
\noindent
(iii) any two $f$-acyclic vector
bundle resolutions of a fixed coherent sheaf in
(ii) are dominanted by a common third one.
\vskip 0.20cm
Now, we are ready to define a morphism
$\text{Err}(\cdot;f,P):K(X)\to H^*(Y,{\Bbb R})$ as follows:
For any coherent sheaf $F$, let $0\to F\to E_0\to E_1\to\dots\to
E_m\to 0$ be a resolution of $F$ by $f$-acyclic vecor
bundles. Set
$$\text{Err}(F;f;P):=\sum_{i=0}^m(-1)^{i}
\text{Err}(E_i;f;P).$$
\vskip 0.20cm
Err is well-defined, i.e., it does not depend on the $f$-acyclic
resolutions of vector bundles we choose. Indeed, if $F$ is a
vector bundle, then by Proposition A.2 and (ii) and (iii), we see
Err is well-defined. All this, together with  (i), implies that
the same is true for coherent sheaves, which then completes the
proof of Proposition A.3.
\vskip 0.20cm
\noindent
(B.4) {\it Proof of Proposition A.4.} 
It is clear that we only need to show this for $f$-acyclic
vector bundles. But then 
$$\text{Err}(E,\rho;f,\tau_f;P)=\text{ch}_{\text{BC}}(E,\rho;f,\tau_f)
-\text{ch}_{\text{BC}}'(E,\rho;f,\tau_f)+f_*\Big(\text{ch}(E)\cdot
\text{td}(T_f)\cdot P(T_f)\Big).$$ Now, by Axiom 3
for relative Bott-Chern secondary characteristic classes in 2.C.2,
we see that $\text{ch}_{\text{BC}}(E,\rho;f,\tau_f)$ and
$\text{ch}_{\text{BC}}'(E,\rho;f,\tau_f)$ are compactible with
base change. Obviously, so does the term $f_*\Big(\text{ch}(E)\cdot
\text{td}(T_f)\cdot P(T_f)\Big)$. This then proves Proposition
A.4.
\vskip 0.20cm
\noindent
{\it Remark 5.3.} In the past, the functorial rule for 
relative Bott-Chern secondary characteristic classes are
systematically denied by a group of very ill-motivated
mathematicians, who simply claim that the functorial
rule I used is completely wrong. Unfortunately, they stand on the
wrong side. Later we will give more consequences for the
functorial rule to indicate its absolute importance.
\vskip 0.20cm
\noindent
(B.5) {\it Proof of Proposition A.5.} Obviously we only need to
show it for $f$-acyclic vector bundles $E$ on $X$. Choose
hermitian metrics $\rho_E$ and $\rho_F$ on $E$ and $F$
respectively, and a hermitian metric $\tau_f$ on $T_f$ such that
all the induced metrics on the fibers of $f$ are K\"ahler. Then,
by Propostion A.1, 
$$\eqalign{~&\text{Err}(E\otimes f^*F;f;P)\cr
=&\text{Err}(E\otimes f^*F;\rho\otimes
f^*F;f,\tau_f;P)\cr
=&\text{ch}_{\text{BC}}(E\otimes f^*F;\rho\otimes
f^*F;f,\tau_f)-\text{ch}_{\text{BC}}'(E\otimes f^*F;\rho\otimes
f^*F;f,\tau_f)\cr
&\qquad+f_*\Big(\text{ch}(E\otimes f^*F)\cdot
\text{td}(T_f)\cdot P(T_f)\Big)\quad(\text{by\ definition})\cr
=&\Big(\text{ch}_{\text{BC}}(E;\rho;f,\tau_f)-
\text{ch}_{\text{BC}}'(E;\rho;f,\tau_f)\Big)\cdot\text{ch}(F,\rho_F)
\cr
&\qquad+f_*\Big(\text{ch}(E)\cdot f^*\text{ch}(F)\cdot
\text{td}(T_f)\cdot P(T_f)\Big)\quad(\text{by\ Axiom\ 2\ in\
2.B.2})\cr =&\text{Err}(E;f;0)\cdot\text{ch}(F,\rho_F)
+f_*\Big(\text{ch}(E)\cdot
\text{td}(T_f)\cdot P(T_f)\Big)\cdot
\text{ch}(F)\cr
&\qquad(\text{by\ definition\ and\ the\ projection\
formula})\cr
=&\text{Err}(E;f;P)\cdot\text{ch}(F).\cr}$$ This
completes the proof of the proposition.
\vskip 0.20cm
\noindent
(B.6) {\it Proof of Proposition A.6.} In the proof of this
proposition, all axioms but the last one, i.e., the rule of the 
deformation to the normal cone, will be used.
\vskip 0.20cm
First, it is well-known that, as a $K(Y)$-module, $K(X)$ is
generated by line bundles ${\Cal O}_X$ and
${\Cal O}_X(-1)$. Obviously, ${\Cal O}_X,$ and
${\Cal O}_X(-1)$ are both $f$-acyclic. So by Propositions A.3 and
A.5, it is sufficient to show that there exists a unique additive
characteristic class
$R$ such that
$$\text{Err}({\Cal O}_X;f;R)=0,\qquad\text{and}\qquad
\text{Err}({\Cal O}_X(-1);f;R)=0.$$ We prove this by a direct
calculation and  Proposition A.4, i.e.,  the
functorial rule.
\vskip 0.20cm
For the time being, we will assume that $F$ is generated by
global sections. Such a condition can be automatically granted if
we tensor $F$ by a very ample line bundle on $Y$. Moreover, it is
well-known that such a change will not change the structure
$f$, nor the sheaves ${\Cal O}$ and ${\Cal O}(-1)$.
\vskip 0.20cm
However, with this latest condition, then, we may find a certain
natural classifying map $\phi:Y\to G(2,2+n)$, where $G(2,2+n)$
denotes the Grassmannian of dimension two subspaces in a complex
$(2+n)$-dimensional space. Moreover, $F$ can be reconstructed as
the pull-back of the canonical rank two subbundle $S_n$ of the
rank $(2+n)$ trivial bundle on $G(2,2+n)$. By the functorial
rule, or better, Proposition A.4, we can then without loss of
generality work on the canonical ${\Bbb P}^1$ bundle $p_n:
{\Bbb P}(S_n)\to G(2,2+n).$ 
\vskip 0.20cm
Note that in a natural way, $\{p_{n}\}_{n\geq 0}$ forms an
inverse limit system: We indeed has the Cartesian product
$$\matrix {\Bbb P}(S_n)&\buildrel
i_{n,m}\over\hookrightarrow&{\Bbb P}(S_{n+m})\\ &&\\
p_n\downarrow&&\downarrow p_{n+m}\\
&&\\
G(2,2+n)&\buildrel
j_n\over\hookrightarrow&G(2,2+(n+m))\endmatrix$$
such that
$i_{n,m}^*{\Cal O}_{{\Bbb
P}(S_{n+m})}(k)={\Cal O}_{{\Bbb
P}(S_{n})}(k).$  Hence, by the
functorial rule, i.e., Proposition A.5, we see that Err for each
$p_n$ for ${\Cal O}_{{\Bbb
P}(S_{n})}$ or ${\Cal O}_{{\Bbb
P}(S_{n})}(-1)$ forms an inverse limit system in
$\{H^*(G(2,2+n),{\Bbb R})\}_{n\geq 0}$.
\vskip 0.20cm
But it is well-known that $$H^*(G(2,2+n),{\Bbb
R})={\Bbb R}[c_1(S_n),c_2(S_n),c_1(Q_n),\dots,
c_n(Q_n)]/(c(S_n)\cdot c(Q_n)=1)$$ where $Q_n$
is the canonical quotient bundle on $G(2,2+n)$ which may be
defined via the canonical  exact sequence
$$0\to S_n\to G(2,2+n)\times {\Bbb C}^{2+n}\to Q_n\to 0,$$ 
$c(S_n)=1+c_1(S_n)+c_2(S_n)$ and $c(Q_n)=1+c_1(Q_n)+
c_2(Q_n)+\dots+c_n(Q_n)$. (See e.g., [BT, \S23].)
Note that the relation above says that $c_j(Q_n)$ can be written
as a polynomial of $c_1(S_n)$ and $c_2(S_n)$ for $j=1,\dots, n$.
More precisely, for $k=2,\dots,n$, we have the recurrence formula
$c_k(Q_n)=-c_1(S_n)\cdot c_{k-1}(Q_n)-
c_2(S_n)\cdot c_{k-2}(Q_n)$ with $c_0(Q_n)=1,$
$c_1(Q_n)=-c_1(S_n)$. Hence, in particular,
 $\lim_{\leftarrow}H^*(G,{\Bbb R})={\Bbb R}[[c_1,c_2]]$.
So,  Err for ${\Cal O}$ and ${\Cal O}(-1)$ can be written as some
universal classes in a ring of power series of $c_1$ and $c_2$.
Set $Err({\Cal O};\pi;P):=P_1(c_1,c_2)$ and
$Err({\Cal O}(-1);\pi;P):=P_2(c_1,c_2)$ be the two power series.
\vskip 0.20cm
Moreover, it is clear that if we change $S_n$ by
tensoring an ample line bundle $L$ on $S_n$, ${\Cal O}, \pi$ and
$P$ will not be changed. But $c_1$ and $c_2$ are changed to
$c_1+2c_1(L)$ and $c_2+c_1(L)(c_1+c_1(L))$ respectively. This
implies that 
$$c_2(S_n)-{1\over
4}c_1^2(S_n)=c_2(S_n\otimes L)-{1\over
4}c_1^2(S_n\otimes L)$$ remains the same. From this, we see that  
$P_1(c_1,c_2)$ is acturally a power series in $c_2-{1\over
4}c_1^2$, which we denote by $P_1$ as well by an abuse of
notation.
\vskip 0.20cm
On the other hand, we have the exact sequence $$0\rightarrow
{\Cal O}_X\rightarrow p^*S^\vee
\otimes {\Cal O}_X(1)\rightarrow {T}_{p}\rightarrow 0.$$
Hence, ${T}_{p}\simeq p^*(\text{det}\,S)^\vee\otimes
{\Cal O}_X(2)$. So ${\Cal O}(-1)\otimes
{1\over 2}p^*(\text{det}S)={1\over 2}T_p^\vee$. But if we tensor
$S$ by a line bundle, $T_p$,  associated to $p$ which remains
the same, will not be changed. Thus, by Proposition A.5, we see
that if we multiply the Err by
$\exp({1\over 2}c_1)$, which takes care of the part of ${1\over
2}p^*(\text{det}S)$
in ${1\over 2}T_p^\vee={\Cal O}(-1)\otimes {1\over
2}p^*(\text{det}S)$, 
$P_2(c_1,c_2)=\text{Err}({\Cal O}(-1);p;R)$ becomes  a
power series in $c_2-{1\over 4}c_1^2$ as well, which we still
denote by
$P_2$  by an abuse of notation.
\vskip 0.20cm
With all this, now we are ready to determine Err
precisely. For this purpose, first, let
$A=c_1({\Cal O}(1))-{1\over 2}p^*c_1(S)$. Then, we have
$$c_1({T}_{p},\tau_{p})=2A,\quad\text{and}\quad 
A^2=c_1({\Cal O}(1))^2-c_1({\Cal O}(1))\cdot
p^*c_1({S})=p^*(-c_2({F})+{1\over 4} c_1({S})^2).$$ 
Therefore,  by the fact
that the direct image of the structure sheaf is the structure
sheaf on the base, $p_*A=1$,
$$p_*A^{2m}=0,\quad\text{and}\quad
p_*(A^{2m+1})=(-c_2({F})+{1\over 4}
c_1({F})^2)^mp_*A=(-c_2({F})+{1\over 4} c_1({F})^2)^m\eqno(5.1)$$
for all positive integer $m$. 
\vskip 0.20cm
Now, we want to see how Err depends on $P$. Obviously, if we
change $P$ by adding $P'$, then the corresponding Err for
${\Cal O}_X$ changes by $$p_*({{2A}\over
{1-e^{-2A}}}P'(2A))$$ and similarly,  for ${\Cal O}_X(-1)$,
(up to the factor $\text{exp}({1\over 2}c_1)$ as discussed above,)
the error changes by
$$p_*({{2Ae^{-A}}\over {1-e^{-2A}}}P'(2A)).$$
In particular, for ${\Cal O}(-1)$,  the
factor before $P'$ is an even function in $A$, and hence is a
series in even powers.  So, by (5.1), we may choose a unique odd
power series $R^{\text{odd}}$, such that for any even power
series $P'$,  $$\text{Err}({\Cal O}_X(-1);p;
R^{\text{odd}}+P')=Err({\Cal O}_X(-1);p;
R^{\text{odd}})=0.$$
\vskip 0.20cm
Similarly, as for ${\Cal O}$,  once the odd part
of $P$ is fixed, by (5.1) again, there is one and
only one even power series $R^{\text{even}}$ such that
$$\text{Err}({\Cal O}_X;p; R^{\text{odd}}+R^{\text{even}})=0,$$
as now only the odd terms of ${{2A}\over
{1-e^{-2A}}}(R^{\text{odd}}(2A)+P'(2A))$ matters.
\vskip 0.20cm
Therefore, finally, we see that there exists a unique power
series $R=R^{\text{even}}+R^{\text{odd}}$ such that
$\text{Err}(\cdot;p;R)=0$ for all ${\Bbb P}^1$-bundles.
This completes the proof of Proposition A.6.
\vskip 0.20cm
\noindent
(B.7) {\it Proof of Proposition A.7.} We divide the proof of this
proposition into two steps;
\vskip 0.20cm
\noindent
(i) Consider the special situation where  codimension-one
closed immersions are  sections of  ${\Bbb P}^1$-bundles.
In this case, we have the following
\vskip 0.20cm
\noindent
{\bf Lemma 1.} {\it Let $f:X\to Y$ be a smooth, proper morphism of
compact K\"ahler manifolds, and $\pi:{\Bbb P}_X^1\to X$ be a
${\Bbb P}^1$-bundle over $X$. Assume that $i:X\to {\Bbb P}^1_X$
is a section, so that $\pi\circ i=\text{Id}_X$. Then for the
closed immersion $i$ over $Y$, and any vector bundle $E$ on $X$,
we have
$$\text{Err}(E;i;P)=0.$$}
\vskip 0.20cm
\noindent
(ii) Deduce a general codimension-one closed immersion to
a section of a ${\Bbb P}^1$-bundle by the deformation to
the normal cone technique. More precise, we have the following
\vskip 0.20cm
\noindent
{\bf Lemma 2.} {\it With the same notation as in
3.A.1, i.e., suppose we have the diagram
$$\matrix {\Bbb
P}^1_X&\buildrel i_\infty\over\hookleftarrow&X&\buildrel
i_0=i\over \hookrightarrow&Z\\
&&&&\\
&g_\infty\searrow&f_\infty\downarrow f_0&\swarrow g_0\\
&&&&\\
&&Y,&&\endmatrix$$ then for any vector bundle $E$ on $Z$, we
have
$$\text{Err}\Big((i_0)_*(i^*E);g_0;P\Big)=
\text{Err}\Big((i_\infty)_*(i^*E);g_\infty;P\Big).$$}
\vskip 0.20cm
Before proving these two lemmas, let us show how Proposition A.7
follows from them.
\vskip 0.20cm
First of all, by Lemma 2, we see that for all coherent sheaves
$F$ on $Z$, $$\text{Err}((i_0)_*i^*F;g_0;P)=
\text{Err}((i_\infty)_*i^*F;g_\infty;P).$$
Obviously, all coherent sheaves on $X$ can be realized as the
pull-back of coherent sheaves from $Z$, so, for all coherent
sheaves $F'$ on $X$, 
$$\text{Err}((i_0)_*F';g_0;P)=
\text{Err}((i_\infty)_*F';g_\infty;P).$$ 
 This then certainly implies
that for all vector bundles $E'$ on $X$, 
$$\text{Err}((i_0)_*E';g_0;P)=
\text{Err}((i_\infty)_*E';g_\infty;P).\eqno(5.2)$$
\vskip 0.20cm
On the other hand, by definition,  $f_0=f_\infty$, we have 
$$\text{Err}(E';f_0;P)=\text{Err}(E';f_\infty;P).\eqno(5.3)$$
Thererfore, 
$$\eqalign{~&\text{Err}(E';i_0;P)\cr
=&\text{Err}(E';f_0;P)-
\text{Err}((i_0)_*E';g_0;P)\qquad(\text{by\ definition})\cr
=&\text{Err}(E';f_\infty;P)-
\text{Err}((i_0)_*E';g_0;P)\qquad(\text{by\ (5.3)})\cr
=&\text{Err}(E';f_\infty;P)-
\text{Err}((i_\infty)_*E';g_\infty;P)\qquad(\text{by\ (5.2)})\cr
=&\text{Err}(E';i_\infty;P)\qquad(\text{by\
definition})\cr
=&0\quad(\text{by\ Lemma\ 1}).\cr}$$ This then is exactly what we
want to show in Proposition A.7. 
\vskip 0.20cm
Now we go back to show Lemma 1 and Lemma 2.
\vskip 0.20cm
\noindent
(i) {\it Proof of Lemma 1.} Here that
 $i$ is a section of a ${\Bbb P}^1$-bundle
over
$X$ plays a very important rule. Indeed, if  $g$
denotes the natural composition ${\Bbb P}^1_X\buildrel \pi\over\to
X\buildrel f\over\to Y$, then by definition,
$$\eqalign{~&\text{Err}(E;i;P)\cr
=&\text{Err}(E;f;P)-\text{Err}(i_*E;g;P)\cr
=&\text{Err}(E;f;P)-
\text{Err}(i_*E;f\circ
\pi;P)\cr
=&\text{Err}(E;f;P)-\Big(\text{Err}(\pi_*\circ
i_*E;f;P)+f_*\big(\text{Err}(i_*E;\pi;P)\cdot
\text{td}(T_{f})\big)\Big)\cr
&\qquad(\text{by\ Proposition
A.2})\cr
=&\text{Err}(E;f;P)-\Big(\text{Err}((\text{Id}_X)_*E;
f;P)+f_*\big(\text{Err}(i_*E;\pi;R)\cdot
\text{td}(T_{f})\big)\Big)\cr
&\qquad(\text{as}\ \pi\circ i=\text{id}_X\ \text{and\ by\ our\
assumption\ on}\  P)\cr
=&-f_*\big(\text{Err}(i_*E;\pi;R)\cdot
\text{td}(T_{f})\big)\cr
=&-f_*\big(0\cdot
\text{td}(T_{f})\big)\quad(\text{by\ Proposition\
A.6})\cr =&0.\cr}$$ 
\vskip 0.20cm
\noindent
{\it Remark 5.4.} Note that in this proof, we only use
the fact that $i$ is a section of a ${\Bbb P}^1$-bundle and for
${\Bbb P}^1$-bundles, Err is simply zero. Thus if
$i$ is a section of a ${\Bbb P}^n$-bundle and, for all ${\Bbb
P}^n$-bundles, Err is zero as well,  we may have a similar
discussion. We will  use this remark in the next chapter.
\vskip 0.20cm
\noindent
(ii) {\it Proof of Lemma 2.} This is where we should use the
deformation to the normal cone rule.
So we recall the following commutative diagram.
$$\matrix X\times \{t\}&&\hookrightarrow&&
X\times {\Bbb P}^1&&\hookleftarrow&&X\times\{\infty\}\\
|&&&&&&&&|\\
|&i_t\searrow&&&I\downarrow&&&\swarrow i_\infty&|\\
|&&&&&&{\Bbb P}&&|\\
|&(t\not=\infty)&Z\times \{t\}&\hookrightarrow
&B_{X\times\{\infty\}}Z\times {\Bbb
P}^1&\hookleftarrow&+&\searrow g_\infty&|\\ |&&&&&&Z&&|\\
f_t&&&\searrow&\pi\downarrow&\swarrow&&&f_\infty\\
|&&&&&&&&|\\
|&g_t\swarrow&&& Z\times {\Bbb P}^1&&&\searrow&|\\
|&&&&&&&&|\\
|&&&&\downarrow&&&&|\\
\downarrow&&&&&&&&\downarrow\\
Y\times\{t\}&&\hookrightarrow&&Y\times{\Bbb
P}^1&&\hookleftarrow&& Y\times\{\infty\}\endmatrix$$
By pulling back $E$ from $Z$ onto $W=B_{X\times\{\infty\}}Z\times
{\Bbb P}^1$ via the composition of natural maps
$W\buildrel\pi\over\to Z\times {\Bbb P}^1\to Z$, we obtain a
vector bundle on
$W$, denoted by $E$ as well by an abuse of notation.  Twisted by
$B_XZ$, we then get the
exact sequence of coherent sheaves on $W$;
$$0\to E(B_XZ-X\times {\Bbb P}^1)\to E(B_XZ)\to
I_*I^*\Big(E(B_XZ)\Big)\to 0.\eqno(5.4)$$ Note that here in
particular
$E(B_XZ-X\times {\Bbb P}^1)$ and $E(B_XZ)$ are vector bundles on
$W$. Easily, we have
$$\eqalign{I_*I^*\Big(E(B_XZ)\Big)\Big|_{X\times
\{t\}}=&(i_t)_*(i^*E),\quad\text{if}\ t\not=\infty;\cr
I_*I^*\Big(E(B_XZ)\Big)\Big|_{{\Bbb P}}=&(i_\infty)_*(i^*E),\cr}$$
and 
$$E(B_XZ-X\times {\Bbb P}^1)=\pi^*(E(-X\times {\Bbb P}^1)\otimes
{\Cal O}_{Z\times {\Bbb P}^1}(Z\times \{\infty\}))
=\pi^*(p_Z^*(E(-X))\otimes p_{{\Bbb P}^1}^*({\Cal O}_{{\Bbb
P}^1}(\infty)).$$ Here $p_Z$ and $p_{{\Bbb P}^1}$ denote the
projection of $Z\times {\Bbb P}^1$ to $Z$ and ${{\Bbb P}^1}$
respectively.
\vskip 0.20cm
On the other hand, by the fact that $W\to {\Bbb P}^1$ is
flat, the restrictions of (5.4) to $Z\times\{t\}$ and ${\Bbb
P}\cup B_XZ$ are again exact. In particular, with the same
notation as in 3.A.2, we have the short exact sequences
$$0\to E_t(-X)\to E_t\to (i_t)_*i_t^*E_t\to 0\qquad\text{on}\quad
Z\times\{t\};$$
$$0\to E_\infty(-X)\to E_\infty\to
(i_\infty)_*i_\infty^*E_\infty\to 0\qquad\text{on}\quad
{\Bbb P};$$ and
$$0\to E_\infty'\to E_\infty''\to 0\to
0\qquad\text{on}\quad B_XZ.$$
Thus, by Proposition A.3, it is
sufficient to show that
$$\eqalign{~&\text{Err}\Big(E(B_XZ-X\times {\Bbb
P}^1)\Big|_{Z\times\{0\}};g_0;P\Big)
-\text{Err}\Big(E(B_XZ)\Big|_{Z\times\{0\}};g_0;P\Big)\cr
=&\text{Err}\Big(E(B_XZ-X\times {\Bbb
P}^1)\Big|_{{\Bbb
P}};g_\infty;P\Big)-\text{Err}\Big(E(B_XZ)\Big|_{{\Bbb
P}};g_\infty;P\Big),\cr}$$ or better, to
show that $$\text{Err}\Big(E(-X);g_0;P\Big)
-\text{Err}\Big(E;g_0;P\Big)
=\text{Err}\Big(E_\infty(-X);g_\infty;P\Big)-\text{Err}\Big(E_\infty
;g_\infty;P\Big).$$ Here $E_\infty$ denotes $E(B_XZ)\Big|_{{\Bbb
P}}$.  From here, by Proposition A.3 again, without loss of
generality, we may assume that $E_t$ and $E_t(-X)$ are all
$g_t$-acyclic.
\vskip 0.20cm
\noindent
{\bf Lemma 3.} {\it With the same notation as above, for all
$t\not=\infty$,
$$\text{Err}\Big(E(-X);g_0;P\Big)
-\text{Err}\Big(E;g_0;P\Big)
=\text{Err}\Big(E_t(-X);g_t;P\Big)
-\text{Err}\Big(E_t;g_t;P\Big).$$}
\vskip 0.20cm
\noindent
{\it Proof.} We first show that
$$\text{Err}\Big(E;g_0;P\Big)=\text{Err}\Big(E_t;g_t;P\Big).
\eqno(5.5)$$
Choose a metric $\rho$ on $E$ as a vector bundle on $Z$. By
Proposition A.1,
$$\eqalign{~&\text{Err}\Big(E;g_0;P\Big)\cr
=&\text{ch}_{\text{BC}}
(E_0,\rho;g_0,\tau_{g_0})-
\text{ch}_{\text{BC}}'
(E_0,\rho;g_0,\tau_{g_0})
+(g_0)_*\Big(\text{ch}(E_0)\cdot \text{td}(T_{g_0})\cdot
P(T_{g_0})\Big),\cr}$$ and 
$$\eqalign{~&\text{Err}\Big(E_t;g_t;P\Big)\cr
=&\text{ch}_{\text{BC}}
(E_t,\rho_t;g_t,\tau_{g_t})-
\text{ch}_{\text{BC}}'
(E_t,\rho_t;g_t,\tau_{g_t})
+(g_t)_*\Big(\text{ch}(E_t)\cdot \text{td}(T_{g_t})\cdot
P(T_{g_t})\Big).\cr}$$ 
\vskip 0.20cm
Note that in the construction of the deformation to the normal
cone, the base curve is a rational curve ${\Bbb P}^1$. As a direct consequence,
the de Rham cohomology classes will not be changed for this
deformation family. This then shows that in $H^*(Y,{\Bbb R})$,
for all $t\not=\infty$,
$$(g_t)_*\Big(\text{ch}(E_t)\cdot \text{td}(T_{g_t})\cdot
P(T_{g_t})\Big)=(g_0)_*\Big(\text{ch}(E_0)\cdot \text{td}(T_{g_0})\cdot
P(T_{g_0})\Big).\eqno(5.6)$$
\vskip 0.20cm
Note that $E_0\simeq E_t$ and $T_{g_0}\simeq T_{g_t}$, so by
changing from 0 to $t$, we may understand that we are working on
the same fibration $g:Z\to Y$ for the same $g$-acyclic vector
bundle $E$. Hence, from $(E_0,\rho_0;g_0,\tau_0)$ to
$(E_t,\rho_t;g_t,\tau_t)$, the only
changes come from the metrics: for
$E$, the metric $\rho$ is changed to $\rho_t$ while for $T_g$, the
metric $\tau_{g_0}$ is changed to $\tau_{g_t}$. Thus by Axiom
4 and Axiom 5 for relative Bott-Chern secondary characteristic
classes in 2.D.2 and 2.E.3 respectively, we see that 
the change of relative Bott-Chern classes
from $0$ to $t$ for $\text{ch}_{\text{BC}}$ is the same as the change of relative Bott-Chern classes
from $0$ to $t$ for $\text{ch}_{\text{BC}}'$, as described in
3.B.2. More precisely, by (3.3.7),
$$\eqalign{~&\text{ch}_{\text{BC}} (E_0,\rho;g_0,\tau_{g_0})-
\text{ch}_{\text{BC}}
(E_t,\rho_t;g_t,\tau_{g_t})\cr
=&-g_*\Big(\text{ch}_{\text{BC}}(E;\rho_t,\rho_0)\cdot
\text{td}(T_g,\tau_t)+\text{ch}(E;\rho_0)\cdot
\text{td}_{\text{BC}}(T_g;\tau_t,\tau_0)\Big)\cr
&\qquad -\text{ch}_{\text{BC}}
\Big(g_*E;L^2(\rho_t,\tau_t),L^2(\rho_0,\tau_0)\Big)\cr
=&\text{ch}_{\text{BC}}'
(E_0,\rho;g_0,\tau_{g_0})-\text{ch}_{\text{BC}}'
(E_t,\rho;g_t,\tau_{g_t}).\cr}$$ This, together with (5.6), then
certainly gives (5.5). Similarly, we have
$$\text{Err}\Big(E(-X);g_0;P\Big)=\text{Err}\Big(E_t(-X);g_t;P\Big).
\eqno(5.7)$$
So finally, by (5.5) and (5.7), we complete  the proof
of Lemma 3.
\vskip 0.20cm
Now we are finally ready to complete the proof of Proposition A.7.
\vskip 0.20cm
For $t\not=\infty$ in the projective line ${\Bbb P}^1$, we have
$$\eqalign{~&\Big(\text{Err}(E_0(-X);g_0;P)-\text{Err}(E_0;g_0;P)\Big)\cr
&\qquad-\Big(
\text{Err}(E_\infty(-X);g_\infty;P)-\text{Err}(E_\infty;g_\infty;P)\Big)\cr
=&\Big(\text{Err}(E_t(-X);g_t;P)-\text{Err}(E_t;g_t;P)\Big)\cr
&\qquad-\Big(
\text{Err}(E_\infty(-X);g_\infty;P)-\text{Err}(E_\infty;g_\infty;P)\Big)\cr
&\qquad(\text{by\
Lemma\ 3})\cr
=&\Big(\text{ch}_{\text{BC}}(E_t(-X),\rho_t';g_t,\tau_{g_t})
-\text{ch}_{\text{BC}}'(E_t(-X),\rho_t';g_t,\tau_{g_t})\cr
&\qquad+(g_t)_*
\big(\text{ch}(E_t(-X))\cdot \text{td}(T_{g_t})\cdot
P(T_{g_t})\big)\Big)\cr
&-\Big(\text{ch}_{\text{BC}}(E_t,\rho_t;g_t,\tau_{g_t})
-\text{ch}_{\text{BC}}'(E_t,\rho_t;g_t,\tau_{g_t})\cr
&\qquad+(g_t)_*
\big(\text{ch}(E_t)\cdot \text{td}(T_{g_t})\cdot
P(T_{g_t})\big)\Big)\cr
&-\Big(\text{ch}_{\text{BC}}(E_\infty(-X),\rho_\infty';g_\infty,
\tau_{g_\infty})
-\text{ch}_{\text{BC}}'(E_\infty(-X),\rho_\infty';g_\infty,\tau_{g_\infty})\cr
&\qquad+(g_\infty)_*
\big(\text{ch}(E_\infty(-X))\cdot \text{td}(T_{g_\infty})\cdot
P(T_{g_\infty})\big)\Big)\cr
&+\Big(\text{ch}_{\text{BC}}(E_\infty,\rho_\infty;g_\infty,\tau_{g_\infty})
-\text{ch}_{\text{BC}}'(E_\infty,\rho_\infty;g_\infty,
\tau_{g_\infty})\cr
&\qquad+(g_\infty)_*
\big(\text{ch}(E_\infty)\cdot \text{td}(T_{g_\infty})\cdot
P(T_{g_\infty})\big)\Big)\cr
&\qquad(\text{by\ definition\ of\
Err})\cr
=&\Big(\Big(\text{ch}_{\text{BC}}(E_t(-X),\rho_t';g_t,\tau_{g_t})
-\text{ch}_{\text{BC}}(E_t,\rho_t;g_t,\tau_{g_t})\Big)\cr
&\qquad-\Big(\text{ch}_{\text{BC}}'(E_t(-X),\rho_t';g_t,\tau_{g_t})
-\text{ch}_{\text{BC}}'(E_t,\rho_t;g_t,\tau_{g_t})\Big)\Big)\cr
&-\Big(\Big(\text{ch}_{\text{BC}}(E_\infty(-X),\rho_\infty';g_\infty,\tau_{g_\infty})
-\text{ch}_{\text{BC}}(E_\infty,\rho_\infty;g_\infty,\tau_{g_\infty})\Big)\cr
&\qquad-\Big(\text{ch}_{\text{BC}}'(E_\infty(-X),\rho_\infty';g_\infty,\tau_{g_\infty})
-\text{ch}_{\text{BC}}'(E_\infty,\rho_\infty;g_\infty,\tau_{g_\infty})\Big)\Big)\cr
&+\Big((g_t)_*
\big(\text{ch}(E_t(-X)-E_t)\cdot \text{td}(T_{g_t})\cdot
P(T_{g_t})\big)\cr
&\qquad-(g_\infty)_*
\big(\text{ch}(E_\infty(-X)-E_\infty)\cdot
\text{td}(T_{g_\infty})\cdot P(T_{g_\infty})\big)\Big)\cr
=&\Big(\Big(\text{ch}_{\text{BC}}(E_t(-X),\rho_t';g_t,\tau_{g_t})
-\text{ch}_{\text{BC}}(E_t,\rho_t;g_t,\tau_{g_t})\Big)\cr
&\qquad-\Big(\text{ch}_{\text{BC}}'(E_t(-X),\rho_t';g_t,\tau_{g_t})
-\text{ch}_{\text{BC}}'(E_t,\rho_t;g_t,\tau_{g_t})\Big)\Big)\cr
&-\Big(\Big(\text{ch}_{\text{BC}}(E_\infty(-X),\rho_\infty';g_\infty,\tau_{g_\infty})
-\text{ch}_{\text{BC}}(E_\infty,\rho_\infty;g_\infty,\tau_{g_\infty})\Big)\cr
&\qquad-\Big(\text{ch}_{\text{BC}}'(E_\infty(-X),\rho_\infty';g_\infty,\tau_{g_\infty})
-\text{ch}_{\text{BC}}'(E_\infty,\rho_\infty;g_\infty,\tau_{g_\infty})\Big)\Big)\cr
&+\Big((g_0)_*
\big(\text{ch}(E_0(-X)-E_0)\cdot \text{td}(T_{g_0})\cdot
R(T_{g_0})\big)\cr
&\qquad-(g_\infty)_*
\big(\text{ch}(E_\infty(-X)-E_\infty)\cdot
\text{td}(T_{g_\infty})\cdot
P(T_{g_\infty})\big)\Big)\cr
&\qquad(\text{by\ (5.6)}).\cr}
$$ From here, if we use Axiom 6 of the deformation to the normal
cone rule for the relative Bott-Chern secondary characteristic
classes in 3.C.3, by taking the limit  $t\to\infty$, we see that
in the above expression, the contribution of
$\text{ch}_{\text{BC}}$ exactly  cancels out the contribution from
$\text{ch}_{\text{BC}}'$. This then implies that
$$\eqalign{\Big(\text{Err}&(E_0(-X);g_0;P)-\text{Err}(E_0;g_0;P)\Big)\cr
&\qquad-\Big(
\text{Err}(E_\infty(-X);g_\infty;P)-\text{Err}(E_\infty;g_\infty;P)\Big)\cr
=&(g_0)_*
\big(\text{ch}(E_0(-X)-E_0)\cdot \text{td}(T_{g_0})\cdot
R(T_{g_0})\big)\cr
&\qquad-(g_\infty)_*
\big(\text{ch}(E_\infty(-X)-E_\infty)\cdot
\text{td}(T_{g_\infty})\cdot P(T_{g_\infty})\big).\cr}\eqno(5.8)$$
\vskip 0.20cm
\noindent
{\bf Lemma 4.} {\it With the same notation as above, 
$$(g_0)_*
\Big(\text{ch}(E_0(-X)-E_0)\cdot \text{td}(T_{g_0})\cdot
R(T_{g_0})\Big)=(g_\infty)_*
\Big(\text{ch}(E_\infty(-X)-E_\infty)\cdot
\text{td}(T_{g_\infty})\cdot P(T_{g_\infty})\Big).$$}
\vskip 0.20cm
\noindent
{\it Proof.} Let $N_{i_t}$ denotes the normal bundle of
$i_t$ for all $t\in {\Bbb P}^1$. Then we have the following
exact sequence
$$0\to T_f\to i_t^*T_{g_t}\to N_{i_t}\to 0.$$ Therefore,
$i_t^*\text{td}(T_{g_t})=\text{td}(T_f)\cdot\text{td}(N_{i_t})$
and $i_t^*P(T_{g_t})=P(T_f)+P(N_{i_t}).$
Hence, for all points $t$ in the projective line ${\Bbb
P}^1$,
$$\eqalign{~&(g_t)_*
\Big(\text{ch}(E_t(-X)-E_t)\cdot \text{td}(T_{g_t})\cdot
P(T_{g_t})\Big)\cr
=&(g_t)_*
\Big(\text{ch}((i_t)_*i_t^*E)\cdot
\text{td}(T_{g_t})\cdot P(T_{g_t})\Big)\cr
=&(g_t)_*
\Big((i_t)_*\big(\text{ch}(i_t^*E)\cdot
\text{td}(N_{i_t})^{-1}\big)\cdot
\text{td}(T_{g_t})\cdot P(T_{g_t})\Big)\cr
&\qquad(\text{by\ Grothendieck-Riemann-Roch\ Theorem\ in\
Algebraic\ Geometry\ for}\ i_t)\cr
=&(g_t)_*
\Big((i_t)_*\big(\text{ch}(i_t^*E)\cdot
\text{td}(N_{i_t})^{-1}\cdot
i_t^*\text{td}(T_{g_t})\cdot i_t^*P(T_{g_t})\big)\Big)\cr
&\qquad(\text{by\ the\ projection\ formula})\cr
=&(g_t)_*
\Big((i_t)_*\big(\text{ch}(i_t^*E)\cdot
\text{td}(N_{i_t})^{-1}\cdot
\text{td}(T_f)\cdot\text{td}(N_{i_t})
\cdot(P(T_f)+P(N_{i_t}))
\big)\Big)\cr
=&f_*\Big(\text{ch}(i_t^*E)\cdot
\text{td}(N_{i_t})^{-1}\cdot
\text{td}(T_f)\cdot\text{td}(N_{i_t})
\cdot(P(T_f)+P(N_{i_t}))\Big)\cr
=&f_*\Big(\text{ch}(i^*E)\cdot
\text{td}(T_f)\cdot (P(T_f)+P(N_{i_t}))\Big).\cr}$$ On the other
hand,
$P(N_{i_\infty})=P(N_{i_t})$. Indeed, since
$X\times {\Bbb P}^1$ does not meet $B_XZ$, and $W\to 
{\Bbb P}^1$ is flat,  the restriction of the normal
bundle $N_I$ of the closed immersion $I:X\times {\Bbb P}^1$ to
$X\times\{t\}$ is simply $N_t$ for all $t\in {\Bbb P}^1$.
Therefore, the cycles corresponding to $\text{ch}(N_{i_t})$ and
$\text{ch}(N_\infty)$ are the same for all $t$. So do
the corresponding de Rham cohomology classes. This then implies
$P(N_{i_\infty})=P(N_{i_t})$ in $H^*(Y,{\Bbb R})$, and hence
completes the proof of Lemma 4.
\vskip 0.20cm
Thus, we finally have;
$$\eqalign{~&\text{Err}(i_*i^*E;g;P)-
\text{Err}((i_\infty)_*i_\infty^*E;g_\infty;P)\cr
=&\Big(\text{Err}(E_0(-X);g_0;P)-\text{Err}(E_0;g_0;P)\Big)\cr
&\qquad-\Big(
\text{Err}(E_\infty(-X);g_\infty;P)-\text{Err}(E_\infty;g_\infty;P)\Big)
\quad(\text{by\
definition})\cr =&(g_0)_*
\big(\text{ch}(E_0(-X)-E_0)\cdot \text{td}(T_{g_0})\cdot
R(T_{g_0})\big)\cr
&\qquad-(g_\infty)_*
\big(\text{ch}(E_\infty(-X)-E_\infty)\cdot
\text{td}(T_{g_\infty})\cdot
P(T_{g_\infty})\big)\quad(\text{by\ (5.8)})\cr
=&0\qquad(\text{by\ Lemma\ 4}),\cr}$$ which completes the proof
of  Proposition A.7. 
\vfill\eject
\vskip 0.20cm
\noindent
{\bf 6. Proof of the Uniqueness for Relative Bott-Chern Secondary
Characteristic Classes}
\vskip 0.20cm
\noindent
A) Rank $n$-Projective Bundles
\vskip 0.20cm
\noindent
B) Closed Immersions of Higher Codimension
\vskip 0.20cm
\noindent
C) Proof of Uniqeness Theorems
\vskip 0.20cm
\noindent
In this chapter, we will complete the proof of the
uniqueness theorems  for relative Bott-Chern secondary
characteristic classes. For doing so, we need to study Err in
two cases:  ${\Bbb P}^n$-bundles and general closed
immersions of codimension $m$. We will use a trick of Bott to
deduce Err for a ${\Bbb P}^n$-bundle and use a trick of
Faltings to deduce Err for a
closed immersion of higher codimension to the cases of  ${\Bbb
P}^1$-bundles and  closed immersions of codimension $1$. Once
this is achieved, we then can use Propositions 5.A.6 and 5.A.7 to
finish the proof.
\vskip 0.20cm
\noindent
{\bf A. Rank $n$-Projective Bundles}
\vskip 0.20cm
\noindent
(A.1) For simplicity, in this section,  we will use  $p_n$ to
denote the projection from any ${\Bbb P}^n$-bundle to its base,
and use
$i_1$ to denote any codimension one closed immersion.
\vskip 0.45cm 
In order to prove Err's are zero for 
projections $p_n$ from ${\Bbb P}^n$-bundles, we use an induction
on $n$. In a more practical term,  we   deduce the 
problems for Err with respect to $p_n$ to these for just $p_1$ and
$i_1$. 
\vskip 0.45cm 
If $n=1$, by Proposition 5.A.6, we
know that there exists a characteristic class $R$ for rank two
bundles such that $$\text{Err}(\cdot;p_1;R)\equiv 0.\eqno(6.1)$$ 
Assume now that for any $m<n$, we have a characteristic class
$R_m$ for rank $m$ vector bundles such that
\vskip 0.20cm
\noindent
($1_m$) $R_m$ is additive and $R_m(E)=R_{m-1}(E)$ for all
rank $m-1$ vector bundles; 
\vskip 0.20cm
\noindent
($2_m$) $\displaystyle{\text{Err}(\cdot;p_m;R_m)\equiv 0.}$ 
\vskip 0.20cm
We then want to prove that there exists a characteristic class
$R_n$ for rank $n$ vector bundles such that
\vskip 0.20cm
\noindent
($1_n$) $R_n$  is additive and $R_n(E)=R_{n-1}(E)$ for all
rank $n-1$ vector bundles; 
\vskip 0.20cm
\noindent
($2_n$) $\displaystyle{\text{Err}(\cdot;p_n;R_n)\equiv 0.}$ 
\vskip 0.20cm
\noindent
(A.2) In order to prove this,  consider the generator of $K(X)$
for
$X={\Bbb P}_Y({F})$, where ${F}$ is a rank $n+1$ vector bundle on
$Y$.  Note that by using
a classifying map to Grassmannian, we may assume that ${F}$ has
a rank 1 sub-line bundle ${L}$ such that $F/L$ is
again a vector bundle. In particular, we have the following
closed embedding of codimension one:
$$\matrix {\Bbb P}_Y(F/L)&&\buildrel i=i_1\over\hookrightarrow&&
{\Bbb P}_Y(F)\\
&&&&\\
&p_{n-1}\searrow&&\swarrow p_n&\\
&&&&\\
&&Y.&&\endmatrix$$
\vskip 0.20cm 
\noindent 
{\bf Lemma 1.} {\it As a $K(Y)$-module, $K(X)$ is generated by
${\Cal O}_X(-1)$ and the direct image of $i_*\Big(K({\Bbb
P}_Y(F/L))\Big)$.}
\vskip 0.45cm
\noindent
{\it Proof.} This is a direct consequence of the facts
that
\vskip 0.20cm
\noindent
(i) we have a short exact sequence $$0\to L\to F\to F/L\to 0$$
and hence $\text{ch}(F)=\text{ch}(F/L)+\text{ch}(L);$
\vskip 0.20cm
\noindent
(ii) $K({\Bbb P}_Y^n)$ is generated by ${\Cal O}(i)$,
$i=1,\dots,n$ and $i_*({\Cal O}(a))={\Cal O}(a)-{\Cal
O}(a-1)$ from the structure exact sequence  
$$0\to {\Cal O}(-1)\to {\Cal O}\to i_*{\Cal O}\to 0.$$
\vskip 0.20cm
\noindent
(A.3) With this lemma, in order to prove the uniqueness
theorems, we only  need to show that there exists an $R_n$
satisfies $(1_n)$ and
$(2_n)$ for ${\Cal O}_X(-1)$ and all the elements in
$i_*\Big(K({\Bbb P}_Y(F/ L))\Big)$.
\vskip 0.20cm
We first deal   with the elements in $i_*\Big(K({\Bbb P}_Y(F/
L))\Big)$. For this purpose, we  use Proposition 5.A.7. In
fact, since  $$i_1: {\Bbb P}_Y(F/ L)\hookrightarrow{\Bbb
P}_Y(F)$$ is a codimension-one closed imbedding, we have
Err$(\alpha; i_1;R_{n})=0$ for all $\alpha\in K({\Bbb
P}_Y(F/L))$.  But by  definition,
$$\text{Err}(\alpha;i_1;R_n)= \text{Err}(\alpha;p_{n-1};R_n)
-\text{Err}((i_1)_*\alpha;p_n;R_n).$$ Here $p_{n-1}$ (resp. $p_n$)
denotes the natural projection from ${\Bbb P}_Y(F)$ (resp. ${\Bbb
P}_Y(F/ L))$ to $Y$. This implies that
$\text{Err}(\alpha;p_{n-1};R_n)=\text{Err}((i_1)_*\alpha;p_n;R_n).$
\vskip 0.20cm
Now, by the induction hypothesis $(1_m)$ and ($2_m)$, $m\leq n-1$,
$$\text{Err}(\alpha;p_{n-1};R_n=R_{n-1})\equiv 0,$$ hence we
have
$$\text{Err}((i_1)_*\alpha;p_n;R_n)\equiv 0,$$ which exactly
means that  $(1_n)$ and $(2_n)$ are valid for the
elements in the direct image of $K({\Bbb P}_Y(F/ L))$.
\vskip 0.20cm 
\noindent
(A.4) Now let us study $\text{Err}({\Cal O}_X(-1),p_n;R_n)$. 
For this special purpose, we use a trick of Bott following
Faltings [F]. 
\vskip 0.20cm 
Let $\text{Flag}_Y(F)$ be the Flag variety of $F$ on $Y$. That 
is, the variety which classifies complete filtrations of $F$: 
$$0=F_0\subset 
F_1\subset ...\subset F_{n+1}=F,$$ where  the successive 
vector bundle quotients are of rank 1. 
There is a natural morphism from $\text{Flag}_Y(F)$ to $X$ which is just the composition
of the forgetting maps. Hence the morphism from
$\text{Flag}_Y(F)$ to $X$ is a composition of ${\Bbb
P}^m$-bundles with $m<n$. Therefore, by Proposition 5.A.2, and the
induction hypothesis,  Err becomes zero
for the morphism $\text{Flag}_Y(F)\rightarrow X$ and any additive
$R$ such that $R=R_m$ for rank $m<n$ vector bundles.
\vskip 0.20cm
On the other hand,   consider the pull-back of the line bundle
${\Cal O}_X(-1)$ over $\text{Flag}_Y(F)$. By a simple
calculation using the projection formula and the fact that the
direct image of the structure sheaf on the total space
$\text{Flag}_Y(F)$ is simply the structure sheaf on the base $X$,
we see  that the push-forward to
$X$ of this pull-back line bundle, denoted by ${\Cal O}'(-1)$,  on
$\text{Flag}_Y(F)$ coincides with ${\Cal O}_X(-1)$ itself. 
Note that ${\Cal O}'(-1)$ is $(\text{Flag}_Y(F)\to X)$-acyclic and
its direct image ${\Cal O}_X(-1)$ is $(X\to Y)$-acyclic. Thus by
Proposition 5.A.2, we have
$$\eqalign{~&\text{Err}({\Cal O}'(-1);\text{Flag}_Y(F)\to
X\buildrel p_n\over
\to Y;R)\cr =&(p_n)_*\Big(\text{Err}({\Cal
O}'(-1);\text{Flag}_Y(F)\to X;R)\cdot
\text{td}(T_{p_n})\Big)+\text{Err}({\Cal O}_X(-1);p_n;R)\cr
=&(p_n)_*\Big(0\cdot
\text{td}(T_{p_n})\Big)+\text{Err}({\Cal
O}_X(-1);p_n;R)\quad(\text{by\ induction\ hypothesis})\cr
=&\text{Err}({\Cal O}_X(-1);p_n;R).\cr}\eqno(6.2)$$
Therefore, it is sufficient to show  that, for the natural
morphism $\text{Flag}_Y(F)\rightarrow Y$ that Err for ${\Cal
O}'(-1)$, the pull-back of ${\Cal O}_X(-1)$, is zero.
\vskip 0.20cm 
In order to deal with the  morphism $\text{Flag}_Y(F)\rightarrow
Y$, we introduce another decomposition: Let  $\text{Flag}_Y'(F)$ be
the flag variety which classifies the following partial
filtrations of $F$: $$0=F_0\subset F_2\subset ...\subset
F_{n+1}=F,$$ where the rank of  $F_k$ is $k$. Then the
composition of the natural morphism from $\text{Flag}_Y(F)$ to
$\text{Flag}_Y'(F)$ with the natural morphism from
$\text{Flag}_Y'(F)$ to
$Y$ is just $\text{Flag}_Y(F)\rightarrow Y$. 
But, the morphism from $\text{Flag}_Y(F)$ to $\text{Flag}_Y'(F)$ 
is a ${\Bbb P}^1$-bundle. Therefore, Err for ${\Cal O}'(-1)$ on
$\text{Flag}_Y(F)$ (with respect to the morphism 
$\text{Flag}_Y(F)\to\text{Flag}_X'F$) vanishes, by our Proposition
5.A.7 for ${\Bbb P}^1$-bundles. On the other hand, the
push-forward of ${\Cal O}'(-1)$ via the ${\Bbb P}^1$-bundle
$\text{Flag}_Y(F)\to\text{Flag}_X'F$ to $\text{Flag}_X'F$
is the zero bundle, which certainly is $(\pi:\text{Flag}_X'F\to
Y)$-acyclic and satisfies the trivial condition that
$\text{Err}(0;\text{Flag}_X'F\to
Y;R)=0$ by Proposition 5.A.3.
 Therefore, by
Proposition 5.A.2 again, we have
$$\eqalign{~&\text{Err}({\Cal O}'(-1);\text{Flag}_XF\to
Y;R)\cr
=&\pi_*\Big(\text{Err}({\Cal O}'(-1);\text{Flag}_XF\to
\text{Flag}_X'F;R)\cdot \text{td}(T_{\text{Flag}_X'F\to
Y})\Big)\cr
&\qquad+\text{Err}(0;{\text{Flag}_XF\to
Y};R)\cr
=&\pi_*\Big(0\cdot \text{td}(T_{\text{Flag}_X'F\to
Y})\Big)+0\cr 
&\qquad(\text{by\ the\ fact\ that}\
\text{Flag}_XF\to
\text{Flag}_X'F\ \text{is\ a}\ {\Bbb P}^1-\text{bundle})\cr
=&0.\cr}$$ This, by (6.2), then gives $\text{Err}({\Cal
O}_X(-1);X\to Y;R)=0,$ and hence completes the proof of
$(2_n)$. In this way, we see that there exists a unique additive 
chacteristic class $R$ such that $\text{Err}(\cdot;p_n;R)=0$
for all ${\Bbb P}^n$-bundles $p_n$.
\vskip 0.20cm
\noindent
{\bf B. Closed Immersions of higher codimension}
\vskip 0.20cm 
\noindent
(B.1) In this section,  we deal with a regular closed
immersion of higher codimension. Similarly, we use the
deformation to the normal cone technique so as  to
deduce a general immersion  to the following situation:
\vskip 0.45cm 
\noindent
(i) a section of a projective
bundle;
\vskip 0.20cm
\noindent
(ii) codimension-one closed immersions.
\vskip 0.45cm 
\noindent
(B.2) Let $$\matrix X&&\buildrel i_n\over\hookrightarrow&&Z\\
&&&&\\
&f\searrow&&\swarrow g&\\
&&&&\\
&&Y&&\endmatrix$$ be a closed immersion $i_n:X\hookrightarrow Z$
of codimension $n$, smooth over $Y$ via $f$ and $g$.
Similarly, as for $i_1$ before, for any vector bundle $E$ on $X$
with $R$ as in Section A, set
$$\text{Err}(E;i_n;R):=\text{Err}(E;f;R)-\text{Err}((i_n)_*E;
g;R).$$ We here want to show that
$\text{Err}(\cdot;p_n;R)=0$ and $\text{Err}(\cdot;i_1;R)=0$
implies $\text{Err}(\cdot;i_n;R)=0.$
\vskip 0.20cm
\noindent
(B.3) Let us now suppose that we can reduce the problem for a
general closed immersion  $i_n$ to the case 1(i) and 1(ii)
above,  then by definition, in case (i),
$$\text{Err}(E;i_n;R)=\text{Err}(E;X\buildrel f\over\to
Y;R)-\text{Err}((i_n)_*E; {\Bbb P}_X^{n}\to X\to Y;R).$$
Note that now $i_n$ is simply a section of ${\Bbb P}^n_X\to X$,
hence with the same proof as for  Lemma 4.B.7.1,  (see e.g.
 Remark 5.4,) we have
$\text{Err}(\cdot;i_n;R)=0$ for $i_n$ a
section of a projective bundle. Note also that by Proposition
5.A.7, Err is zero for any closed immersion of codimension one, 
so we should find a way to deduce a general closed immersion to a
section of a ${\Bbb P}^n$-bundle and codimension one
closed immersions. 
\vskip 0.20cm
\noindent
(B.4) To deduce the case of an arbitrary closed immersion to 1(i)
and 1(ii) above, we use deformation to the normal cone theory as
usual. For this,  recall the following basic fact concerning
the  deformation to the normal cone. 
\vskip 0.20cm
Let $f:X\to Y$ and $g: Z\to Y$ be two smooth, proper morphisms of compact
K\"ahler manifolds and $i:X\to Z$ be a codimension $n$ closed
immersion over $Y$, i.e., $i$ is a closed immersion of
codinemsion $n$ such that $f=g\circ i$. Then we have the
following standard construction of the deformation to the normal
cone.
\vskip 0.20cm
Denote by $\pi: W:=B_{X\times \{\infty\}}Z\times {\Bbb
P}^1\rightarrow Z\times {\Bbb P}^1,$  where
$B_{X\times \{\infty\}}Z\times {\Bbb P}^1$ denotes the blowing-up of
$Z\times {\Bbb P}^1$ along $X\times \{\infty\}$.  
Denote the exceptional divisor by ${\Bbb
P}$. It is well-known that the map $q_W:W\rightarrow {\Bbb P}^1$, obtained by
composing $\pi$ with the projection
$q:Z\times {\Bbb P}^1\rightarrow {\Bbb P}^1$, is flat, and that for $z\in
{\Bbb P}^1$:
$$q^{-1}(z)=\cases Z, &$for $z\not=\infty$$,\\
{\Bbb P}\cup B_XZ,&$for $z=\infty$$.\endcases$$ 
Here $B_XZ$ denotes the blowing-up of $Z$ along
$X$. By the construction,
${\Bbb P}\cap B_XZ$ is the exceptional diviosr of $B_XZ$. In
particular, ${\Bbb
P}^n={\Bbb P}_X({N_{i_n}}\oplus {\Cal
O}_X)$ with $N_{i_n}$ the normal bundle corresponding to  $i_n$.
\vskip 0.20cm
Denote by $I:X\times {\Bbb P}^1\hookrightarrow W$  the induced 
codimension $n$ closed embedding. Easily we see that the image of
$I$ does not intersect with $B_XZ$, and the image $X\times
\{\infty\}$ in $W$ is a section of ${\Bbb P}$.
\vskip 0.45cm 
\noindent
{\bf Lemma.} {\it With the same notation as above,
the following two morphisms $$X \buildrel
 {i_\infty} \over \hookrightarrow {\Bbb P}_X({N_{i_n}}\oplus {\Cal
O}_X)\buildrel {j_\infty}\over 
\hookrightarrow W\qquad \text{and}\qquad X \buildrel {i_0} \over
\hookrightarrow W_0 \buildrel {j_0} 
\over\hookrightarrow W$$ induce the same morphism for 
$K$-groups.} 
\vskip 0.45cm 
\noindent
{\it Proof.} Clearly, as $W$ is flat over the rational curve
${\Bbb P}^1$, and $K$-group is compactible with rational
equivalence on cycles by the fact that there is a ${\Bbb
Q}$-isomorphism between $K$-groups  and Chow groups,
 we see that
$X\hookrightarrow W_t$ for all $t\in {\Bbb P}^1$ induce the same
map on $K$-groups.
\vskip 0.20cm
Now let $0\to E.\to I_*F\to 0$ be any resolution of a vector
bundle on $X$, (viewed as a vector bundle on $X\times {\Bbb
P}^1$.) Then by the flatness of $W\to {\Bbb P}^1$, 
 the restrictions of resolutions $0\to E.\to I_*F\to 0$ to $W_t$
for all $t\in {\Bbb P}^1$ are all exact. 
On the other hand, by the construction, $X\times {\Bbb P}^1$ is
disjoint from the  component $B_XZ$ in $W_\infty$, so we
see that on $B_XZ$, such a restriction is acyclic. That is to
say, it gives zero element in $K(B_XZ)$ and as well as in
$K({\Bbb P}\cap B_XZ)$. Thus by definition,  we see that the
$K$-push-forward for $X\to W_\infty\to W$ may be simply
calculated by the one induced via
$X\to {\Bbb P}\to W$. This completes the proof of the lemma.
\vskip 0.45cm 
In this way, by definition, we see that 
$$\text{Err}(E;{j_\infty\circ i_\infty};R)=\text{Err}(E;{j_0 
\circ i_0};R).$$ By definition again, we know that
$$\text{Err}(E;{j_\infty\circ i_\infty};R)=
\text{Err}(E; i_\infty;R)+\text{Err}(i_{\infty *}E; j_\infty;R)$$
and $$\text{Err}(E; {j_0\circ i_0};R)=\text{Err}(E;
i_0;R)+\text{Err}(i_{0 *}E; j_0;R).$$ Thus to  complete the proof,
it is sufficient to prove that $$\text{Err}(\cdot;
i_\infty;R)\equiv 0,$$ 
$$\text{Err}(\cdot; j_\infty;R)\equiv 0,$$  and
$$\text{Err}(\cdot; 
j_0;R)\equiv 0.$$ Note that each of the  three closed immersions,
$i_\infty$, $j_0$ and $j_\infty$,   is either a codimension-one
closed immersion or a section of a ${\Bbb P}^n$-bundle, i.e.,
they belong exactly to 1(i) and 1(ii) above.  So we
see that $\text{Err}(\cdot;i_n;R)\equiv 0$ for any closed
immersion $i_n$ of codimension $n$.
\vskip 0.20cm
\noindent
{\bf C. The Proof of Uniqueness Theorems} 
\vskip 0.20cm
\noindent
(C.1) Clearly, the weak uniqueness theorem  is a direct
consequence of the strong uniqueness theorem, so it is
sufficient to prove the later one.
\vskip 0.20cm
\noindent
(C.2) For doing so, factor $f:X\to Y$ as a regular closed
immersion
$i_m:X\to Z$ followed by a projection $p_n:{\Bbb P}^n_Y\to Y$.
Then by definition, for any vector bundle $E$ on $X$,
$$\eqalign{~&\text{Err}(E;f;R)\cr
=&\text{Err}(E;i_m;R)+\text{Err}(i_*E;p_n;R)\cr
=&0+\text{Err}(i_*E;p_n;R)\qquad(\text{by\ (B.4)})\cr
=&0+0\qquad(\text{by\ (A.3)\ and\ (A.4)})\cr
=&0.\cr}$$ This together with
Proposition 5.A.3 implies that $$\text{Err}(\cdot;f;R)\equiv 0.$$
Thus by definition, we have
$$\text{ch}_{\text{BC}}'(E,\rho;f,\tau_f)=
\text{ch}_{\text{BC}}(E,\rho;f,\tau_f)+f_*\Big(\text{ch}(E)\cdot
\text{td}(T_f)\cdot R(T_f)\Big)$$ for any smooth proper morphism
$f$ and any $f$-acyclic vector bundle $E$. This then surely
completes the proof of the uniqueness theorems.
\vfill\eject
\vskip 0.20cm
\noindent
{\bf 7. Existence of Relative Bott-Chern
Secondary Characteristic Classes}
\vskip 0.20cm
\noindent
A) Arithmetic intersection and arithmetic characteristic
classes
\vskip 0.20cm
\noindent
B) An Effective Construction of Relative Bott-Chern
Secondary Characteristic Classes
\vskip 0.20cm
\noindent
In this chapter, we prove a weak version of the  existence of
relative Bott-Chern secondary characteristic classes by
effectively constructing some classes of differential forms,
which is sufficient and necessary for our application to the
arithmetic Grothendieck-Riemann-Roch theorem.  
\vskip 0.20cm
\noindent
{\bf A. Arithmetic intersection and arithmetic characteristic
classes}
\vskip 0.20cm
\noindent
(A.1) In this section, we first recall the theory of arithmetic
intersection and arithmetic characteristic classes developed
by Arakelov [Ar1,2], Deligne [De2] and Gillet-Soul\'e [GS1,2]. All
results in this section are mainly due to [GS1,2]. The notable
difference is that here we only work over ${\Bbb C}$. So some
modifications are needed. (Contary to the general principle used
in this paper, in this section, we will only tell the reader how
an arithmetic intersection and arithmetic characteristic
classes are introduced, instead of indicating a more sound why.) 
\vskip 0.20cm
Let $X$ be a complex compact manifold of dimension $d$. Denote
 the space of differential forms of degree $n$ on $X$ by 
${A}^n(X):=\oplus_{p+q=n}{A}^{p,q}(X)$. There  are natural
boundary morphisms
$\partial:{A}^{p,q}(X)\rightarrow {A}^{p+1,q}(X),
\bar\partial:{A}^{p,q}(X)\rightarrow {A}^{p,q+1}(X)$, 
and the usual differential $d:{A}^n(X)\rightarrow {A}^{n+1}(X)$.
We say that a  linear function $T$ on ${A}^n(X)$ is a {\it
current}, if $T$ is continuous in the sense of  Schwartz: for any
sequence
$\{\omega_r\}\subset {A}^n(X)$ with the supports contained in
certain fixed compact subset $K$, 
$T(\omega_r)\rightarrow 0$ if all the coefficients of
$\omega_r$ together with their derivatives tend uniformly to
zero when $r\rightarrow \infty.$ The set of currents forms a
topological dual space ${A}(X)^*$  of ${A}(X)$. Denote by
${D}_n(X):={A}^n(X)^*.$  There is a natural decomposition
${D}_n(X)=\oplus_{p+q=n}{D}_{p,q}(X),$ where 
${D}_{p,q}(X)$ is the dual of ${A}^{p,q}(X).$ It is convenient
to set ${D}^{p,q}(X):= {D}_{d-p,d-q}(X).$ Then
 $\partial,\bar\partial$, and $d$ induce morphisms
$\partial',\bar\partial',d'$ from  $D^{p,q}$ to  
$D^{p+1,q}(X),$ $D^{p,q+1}(X),$ and $D^{p+1,q+1}(X)$
respectively, e.g.,  
$(\partial'T)(\alpha):=T(\partial \alpha).$
\vskip 0.20cm
\noindent
{\bf Examples.} (i) There is a natural inclusion
$$\matrix A^{p,q}(X)&\hookrightarrow & D^{p,q}(X)\\
\omega&\mapsto&[\omega]\endmatrix,$$ where
$[\omega](\alpha):=\int_X
\omega\wedge \alpha$ for any $\alpha\in A^{d-p,d-q}(X).$ We
say that a current $T$ is {\it smooth} if there exists a
smooth form $\omega$ such that $T=[\omega].$ In particular, if
$p+q=n$, it follows by Stokes' theorem that 
$[d\omega](\alpha)=\int_Xd\omega\wedge \alpha
=\int_Xd(\omega\wedge \alpha)-\int_X(-1)^n\omega\wedge d\alpha
=(-1)^{n+1}\int_X\omega\wedge
d\,\alpha=(-1)^{n+1}(d'[\omega])(\alpha).$ Therefore,  we
let $\partial,\bar\partial,d$ on the currents be
$(-1)^{n+1}\partial',\  (-1)^{n+1}\bar\partial',\ (-1)^{n+1}d'$
respectively, and let $d^c:= {1\over {4\pi
i}}(\partial-\bar\partial)$.  Then $dd^c=-{1\over{2\pi
i}}\partial\bar\partial$ is a real operator, and we have the
following commutative diagram:
$$\matrix A^{p,q}(X)&\hookrightarrow &D^{p,q}(X)\\
\partial \downarrow &&\downarrow \partial\\
A^{p+1,q}(X)&\hookrightarrow &D^{p+1,q}(X).\endmatrix$$ 
\vskip 0.20cm
\noindent
(ii) Let $i: Y\hookrightarrow X$ be an irreducible 
subvariety of codimension 
$p$. We get a current $\delta_Y\in D^{p,p}(X)$ by letting 
$\delta_Y(\alpha):=\int_{Y^{\text{ns}}}
i^*\alpha$ for any $\alpha\in A^{d-p,d-p}(X).$ Here
$Y^{\text{ns}}$ denotes the non-singular locus of $Y$. We call
this current the {\it Dirac symbol} of $Y$. 
\vskip 0.20cm
Concerning the relations between smooth differential forms and
currents, we have the following well-known facts: 
\vskip 0.20cm
\noindent
(i)  With the boundary morphisms
$\partial,\ \bar\partial,\ d$, the cohomologg goups of $X$ for
differential forms are isomorphic to these for  currents.
\vskip 0.20cm
\noindent
(ii)  Let $\gamma$ be a current on $X$ such that
$dd^c\gamma$ is smooth. Then there exist currents $\omega,\
\alpha,\ \beta$ such that $\gamma=\omega+\partial \alpha+
\bar\partial \beta$, with $\omega$ smooth.
\vskip 0.20cm
\noindent
(iii) As a current, if $\omega$ is smooth and
$\omega=\partial u+\bar\partial v$,  then there
exist smooth currents $\alpha,\ \beta$ such that $\omega=\partial \alpha+\bar\partial
\beta.$
\vskip 0.20cm
\noindent
(iv) If $X$ is a K\"ahler manifold, and $\eta\in
D^{p,q}(X),\ p,\,q\geq 1$, is $d$-closed and is either $d,
\ \partial,\ \bar\partial$ exact. Then there exists $\gamma\in D^{p-1,q-1}(X)$ such that
$dd^c\gamma=\eta.$ In particular, if $\eta=0$, we may choose
$\gamma=\omega+\partial
\alpha+\bar\partial \beta$ with $\omega$  a harmonic form.
\vskip 0.20cm
\noindent
(A.2) Let $Y$ be a codimension $p$ analytic subvariety of $X$.
Following [GS1], we say that a current
$g\in D^{p-1,p-1}(X)$ is a {\it Green's current} of $Y$ if
$dd^cg=[\omega]-\delta_Y$ for some
$\omega\in A^{p,p}(X).$ It is well-known that if $X$ is a
K\"ahler manifold, then Green's currents for
analytic subvarieties of $X$ exist. 
And if $g_1$ and $g_2$ are two Green's currents for
$Y$, then, by A.2(ii),  $g_1-g_2=[\eta]+\partial S_1 +\bar\partial
S_2,$ where $\eta\in A^{p-1,p-1}(X).$
\vskip 0.20cm
\noindent
{\bf Example.} ({\bf The Poincar\'e-Lelong equation}) Let
$(L,\rho)$ be  a hermitian line bundle on $X$ and $s$ a
non-zero meromorphic  section of ${L}$. Then
$-\log |s|_\rho^2\in L^1(X)$, and hence induces a distribution
$[-{\log}|s|_\rho^2]\in D^{0,0}(X)$. Moreover,
$$dd^c[-{\log}|s|_\rho^2]
=[c_1({L},\rho)]-\delta_{\text{div}(s)}.$$ So,
$[-{\log}|s|_\rho^2]$ is a  Green's current of
$\text{div}(s)$, the divisor of $s$.
\vskip 0.20cm
\noindent
(A.3)  From now on, assume that $X$ is a projective
complex manifold. For any irreducible subvariety $Y$,
following [GS1], we say a smooth form $\alpha$ on $X-Y$ has
{\it logarithmic  growth} along $Y$, if there exists a proper
morphism
$\pi:\tilde X\rightarrow X$ such that
$E:=\pi^{-1}(Y)$ is a divisor with normal crossings, $\pi: \tilde X-E\simeq X-Y$ and 
$\alpha$ is  the direct image of a form $\beta$ on $\tilde X-E$
by $\pi$ with the following property: Near each $x\in \tilde
X$, let
$z_1\dots z_k=0$ be a local defining equation of $E$. Then,
there exists $d$-closed smooth forms $\alpha_i$ and a smooth
form $\gamma$ such that
$\beta=\sum_{i=1}^k\alpha_i{\log}|z_i|^2+\gamma.$ 
\vskip 0.20cm
Obviously, such an $\alpha$ is always locally integrable on $X$, 
and hence defines a current $[\alpha]$, which is the direct image by $\pi$ of the current $[\beta]$.
\vskip 0.20cm
By [GS1], we know that for every irreducible subvariety 
$Y\subset X$, there exists a smooth form $g_Y$ on $X-Y$ with
logarithm growth along $Y$ such that $[g_Y]$ is a Green's
current for $Y$. (This may be viewed as a generalization of
the Poincar\'e-Lelong equation.) Moreover, if  $\alpha$
is a form on $X-Y$ with logathmic growth along $Y$, then 
\vskip 0.20cm
\noindent
(i)  if $f:X'\rightarrow X$ is
a morphism of smooth  projective varieties such that
$f^{-1}(Y)$ does not contain any  component of $X'$, then the
form $f^*(\alpha)$ is of logarithmic  growth along $f^{-1}(Y)$;
\vskip 0.20cm
\noindent
(ii)   $d[\alpha]=[d\alpha].$
\vskip 0.20cm
\noindent
(A.4) We introduce now the arithmetic Chow groups and their
cohomological properties following [GS1].
\vskip 0.45cm
Let $X$ be a regular projective variety over ${\Bbb C}$. Set
$$\tilde A^{p,p}(X):=A^{p,p}(X)/(\text{Im}
\partial+\text{Im}\bar\partial);\quad\tilde
A(X):=\oplus_p\tilde A^{p,p}(X);$$ $$\tilde
D^{p,p}(X):=D^{p,p}(X)/(\text{Im}
\partial+\text{Im}\bar\partial);\quad\tilde
D(X):=\oplus_p\tilde D^{p,p}(X).$$ 
\vskip 0.45cm
Denote by $Z^p(X)$ the free group generated by subvarieties of
$X$. We say that an element
$(Z,g_Z)\in Z^p(X)\oplus
\tilde D^{p-1,p-1}(X)$ is an {\it arithmetic
$p$-cycle} if
$g_Z$ is a Green's current of $Z$, i.e.
$dd^cg_z=\omega(Z,g_Z)-\delta_Z$ for some
$\omega_Z:=\omega(Z,g_Z)\in A^{p,p}(X).$ Denote by
$Z_{\text{Ar}}^p(X)$ the abelian group generated by arithmetic
$p$-cycles. 
\vskip 0.20cm
Next, define arithmetic rational equivalence among  
arithmetic cycles. Let $i:Y\hookrightarrow X$ be a subvariety
of codimension $p-1$. Then,  there is a resolution  of
singularities of $Y$, $\pi:\tilde Y\rightarrow Y$ with
$\pi$ proper. Thus any rational function $f\in k(Y)^*$ defines
a rational function $\tilde f$ on $\tilde Y$ such that
$\log |\tilde f|^2$ is $L^1$ on $\tilde Y$. Hence $\tilde
f$ is contained in $D^{0,0}(\tilde Y)$. Let $\tilde i:\tilde
Y\rightarrow X$ be the natural induced morphism, then $\tilde
i_*[\log |\tilde f|^2]\in D^{p-1,p-1} (X),$ and
is independent of the choice of $\tilde Y$. Hence we may
simply denote it by $i_*[{\log} |f|^2].$ Moreover, by the
Poincar\'e-Lelong equation in Example A.2, 
$$\text{div}_{\text{Ar}}(f):=(\text{div}(f),
-i_*[{\log}|f|^2])\in Z^p_{{\text{Ar}}}(X),$$  and is defined to
be {\it arithmetically rationally equivalent to zero}. Let
$R_{\text{Ar}}^p(X)$ be  the subgroup of $Z_{\text{Ar}}^p(X)$
generated by  
$\text{div}_{\text{Ar}}(f)$
for $f\in k(W)^*$, with  $W$ any codimension-$(p-1)$ 
subvariety. Define the {\it p-th arithmetic Chow
group}, denoted by $\text{CH}_{\text{Ar}}^p(X)$, to be
 the quotient group
$Z_{\text{Ar}}^p(X)/R_{\text{Ar}}^p(X).$
Set
$\text{CH}_{\text{Ar}}(X):=\oplus_p\text{CH}^p_{\text{Ar}}(X).$
\vskip 0.20cm
\noindent
(A.5) Let $X$ be a smooth projective variety over
${\Bbb C}$, and 
$Y\subset X$ a closed irreducible subvariety and $f:Z\rightarrow
X$  a proper morphism of irreducible projective varieties  
such that $f(Z)\not\subset Y$. Let $g_Y$ be a Green's current
of $Y$ with logarithmic growth along Y associated to $(\tilde
X,E)$ as in A.3. Denote the associated current
by $[g_Y]$. Then by resolving the singularities of
$Z$, we may construct a commutative diagram
$$\matrix \tilde Z&
\buildrel j\over \rightarrow &\tilde X\\
p\downarrow&\searrow q&\downarrow \pi\\
Z&\buildrel f\over\rightarrow&X,\endmatrix$$ such that
$D=j^{-1}(E)$ is a divisor with normal crossings, $\tilde Z$
is projective and smooth, and $p$ is birational. 
By A.3(i), $q^*g_Y$ is of logarithmic growth along $q^{-1}(Y)$,
so it is integrable and 
$[g_Y]\wedge\delta_Z:=q_*[q^*g_Y]$ defines a current in $X$.
Furthermore, if
$g_Z$ is an arbitrary Green's current of $Z$, we define the
*-product of $[g_Y]$ and $g_Z$ by 
$$[g_Y]*g_Z:=[g_Y]\wedge \delta_Z+[dd^cg_Y+
\delta_Y]\wedge g_Z.$$
\vskip 0.45cm
One checks that  the following holds;
\vskip 0.20cm
\noindent  
(i) If $Y$ and $Z$ intersect properly, i.e.,
$Y\cap F=\cup_iS_i$ with $\text{codim}_XS_i=\text{codim}_XY
+\text{codim}_XZ$, and as algebraic cycles,
$[Y][Z]=\sum_i\mu_iS_i$, then
 $dd^c([g_Y]*g_Z)=[\omega_Y\wedge\omega_Z]
-\sum_k\mu_k\delta_{S_k}.$ 
\vskip 0.20cm
\noindent  {(ii)} For any two Green's currents $g_Y,g_Y'$ of
$Y$ with logarithmic growth, as an element of $\tilde D(X)$,
$[g_Y]*g_Z=[g_Y']*g_Z$ for any Green's current
$g_Z$. Hence, we may also define the *-product  among
general Green's currents. In particular, 
$g_Y*g_Z=g_Z*g_Y,$ and $(g_Y*g_Z)*g_W=g_Y*(g_Z*g_W)$ whenever the
products make sense.
\vskip 0.20cm
\noindent
(A.6) Easily, we have the following  morphisms involving
$\text{CH}_{\text{Ar}}(X)$;
\vskip 0.20cm
\noindent {(i)}
$\zeta:\text{CH}_{\text{Ar}}^p(X)\rightarrow \text{CH}(X),\ \
\ \  (Z,g_Z)\mapsto Z.$
\vskip 0.20cm
\noindent {(ii)} $a:\tilde A^{p-1,p-1}(X)\rightarrow
\text{CH}_{\text{Ar}}^p(X),\ \ \ \ \alpha\mapsto (0,\alpha).$
\vskip 0.20cm
\noindent {(iii)} $\omega:\text{CH}_{\text{Ar}}^p(X)\rightarrow
A^{p,p}(X),\ \ \ \ (Z,g_Z)\mapsto dd^cg_Z+\delta_Z.$
\vskip 0.20cm
One checks that we have the following exact
exact sequence
$$ \tilde A^{p-1,p-1}(X)
\buildrel a\over\rightarrow \text{CH}^p_{\text{Ar}}(X)\buildrel\zeta\over\rightarrow \text{CH}^p(X)
\rightarrow 0.$$
\vskip 0.45cm
\noindent
(A.7) Now we are ready  following [GS1] to introduce an
arithmetic intersection theory for regular  projective
varieties over ${\Bbb C}$. 
\vskip 0.20cm
\noindent 
{\bf Theorem.} ([GS1]) {\it Let $X$ be a regular projective
variety over ${\Bbb C}$. Then 
\vskip 0.20cm
\noindent  {(i)} for each pair of natural numbers
$(p,q)$, there is a pairing
$$\matrix \text{CH}_{\text{Ar}}^p(X)\otimes
\text{CH}_{\text{Ar}}^q(X)&\rightarrow
&\text{CH}_{\text{Ar}}^{p+q}(X)_{\Bbb Q}\\
\alpha\otimes \beta&\mapsto&\alpha\beta.\endmatrix$$ The
pairing is uniquely determined by the following  property: 
If $Y$ and $Z$ are subvarieties of $X$ which intersect 
properly, and $g_Y$
and $g_Z$ are Green's currents for $Y$ and $Z$, then 
$$([Y],g_Y)([Z],g_Z):=([Y][Z],g_Y*g_Z);$$ 
\vskip 0.20cm
\noindent  
{(ii)} the product above makes
$\text{CH}_{\text{Ar}}(X)_{\Bbb
Q}:=\oplus_p\text{CH}_{\text{Ar}}^p(X)_{\Bbb Q}$
 a commutative, associative ${\Bbb Q}$-algebra;
\vskip 0.20cm
\noindent  
{(iii)} the natural morphism
$$(\zeta,\omega):\oplus_p\text{CH}_{\text{Ar}}^p(X)_{\Bbb
Q}\rightarrow
\oplus_p(\text{CH}^p(X)\oplus Z^{p,p}(X))_{\Bbb Q}$$ is a
${\Bbb Q}$-algebra homomorphism,  where $Z^{p,p}(X):=\text{the \
closed\ forms\ in}\ A^{p,p}(X)$. Moreover, $a(\phi)\cdot
(Z,g)=a(\phi\wedge
\omega(Z,g))$ for all $\phi\in\tilde A(X)$ and $(Z,g)\in
\text{CH}_{\text{Ar}}(X)$, and in particular,
 $a(\oplus_p H^{p,p}(X))$ is a square zero ideal in
$\text{CH}_{\text{Ar}}(X)_{\Bbb Q}$.}
\vskip 0.20cm
\noindent
(A.8) Let $f:X\to Y$ be a smooth morphism of regular
projective varieties over ${\Bbb C}$. Then for any subvariety
$Z$ in $Y$, if $g_Z$ is a Green's curent of $Z$,
 $f^*g_Z$ is  a Green's current for $f^{-1}(Z)$. One
checks that the arithmetic rational equivalence is compactible
with such a pull-back. Hence we have a well-defined
pull-back morphism
$f^*:\text{CH}_{\text{Ar}}^p(Y)\rightarrow
\text{CH}_{\text{Ar}}^p(X)$, which 
 induces a ${\Bbb Q}$-algebra morphism
$\text{CH}_{\text{Ar}}(Y)_{\Bbb Q}\rightarrow
\text{CH}_{\text{Ar}}(X)_{\Bbb Q}$
\vskip 0.20cm
We may also construct a puh-forward morphism. This is done as
follows. First  construct a map from $Z_{\text{Ar}}^p(X)$ to
$Z_{\text{Ar}}^{p-r}(Y)$ as follows, where $r$ denotes the
relative dimension of $f$: Let 
$(Z,g_Z)\in Z_{\text{{Ar}}}^p(X),$ with $Z$ irreducible, i.e.
$Z={\overline {\{z\}}}$ with $z$  the generic point of $Z$. Set
$$f_*(Z):=\cases [k(z):k(f(z))]\,{\overline {\{f(z)\}}}, &$if
$\text{dim}\,f(z)=\text{dim}\,z$;$\\ 0,&$
otherwise$.\endcases$$ Then, for  Green's currents, we know
that for any 
$\eta\in \tilde
A^{\text{dim}Y-p,\text{dim}Y-p}(X)$, 
$$\eqalign {(f_*\delta_Z)(\eta)=&\delta_Z(f^*\eta)
=\int_{Z}f^*\eta=\int_{Z}
f^*(\eta|_{f(Z)})\cr
=&\cases \text{deg}(Z/f(Z))\,\int_{f(Z)}\eta,&$if $Z\rightarrow
f(Z)$ is finite$;\\
0,&$otherwise.$\endcases\cr}$$
Hence $f_*\delta_Z=\delta_{f_*(Z)},$  and 
 $$dd^c(f_*g_Z)=[f_*\omega_Z]-\delta_{f_*(Z)}.$$ That is, $f_*g_Z$ defines a Green's current
of $f_*(Z)$. Therefore we may set $$f_*(Z,g_Z):=(f_*Z,f_*g_Z)\in
Z_{\text{Ar}}^{p-r}(Y).$$ Furthermore, it is not difficult to
check  that this  definition is compatible with the arithmetic
rational equivalence, and hence we get a push-out morphism $f_*$
for arithmetic Chow groups
$f_*:\text{CH}_{\text{Ar}}^p(X)\rightarrow 
\text{CH}_{\text{Ar}}^{p-r}(Y),$ which induces a ${\Bbb
Q}$-algebra morphism
$\text{CH}_{\text{Ar}}(X)_{\Bbb Q}\rightarrow
\text{CH}_{\text{Ar}}(Y)_{\Bbb Q}$. Moreover, 
 we have the {\it projection formula}
$$f_*(f^*(\alpha)\cdot\beta)=\alpha\cdot f_*(\beta).$$
\vskip 0.20cm
One may also introduce a pull-back morphism for arithmetic Chow
rings with respect to closed immersions $i:X\hookrightarrow Z$
of regular projective varieties over ${\Bbb C}$. In fact if $W$
is a subvariety of $Z$ which properly intersects with $X$, then
it is easily to check that for any Green's current $g_W$ of $W$,
$i^*g_Z$ is a Green's current of $W\cap Z$. Hence, 
$i^*(W,g_W)\in \text{CH}_{\text{Ar}}(Z)$. For general arithmetic
cycles, one may use the arithmetic intersection to define their
pull-back onto $Z$. So a suitable  Chow type moving lemma is
needed. All in all, as what does in [GS1], for $i$, we finally can
define a ${\Bbb Q}$-algebra morphism
$i^*:\text{CH}_{\text{Ar}}(Z)_{\Bbb Q}\to
\text{CH}_{\text{Ar}}(X)_{\Bbb Q}$ which coincides $i^*$ defined
above for arithmetic cycles whose corresponding algebraic cycles
are properly intersect with $Z$. 
\vskip 0.20cm
Thus, in particular, if $f:X\to Y$ is a  morphism between regular
projective varieties over ${\Bbb C}$, then  write $f$ as a
composition of a closed immersions $i$ and a smooth moprphism
$p$. Set
$f^*:=p^*\circ i^*: \text{CH}_{\text{Ar}}(Y)\to
\text{CH}_{\text{Ar}}(X)$. One checks that such an $f^*$ does not
depend on the decomposition $f=p\circ i$. 
\vskip 0.20cm
Moreover, if $f:X\to Y$ is smooth and $g:Y'\to Y$ is proper
morphism among regular projective varieties over ${\Bbb C}$.
Denote by
$g_f:X\times_YY'\to X$ and $f_g:X\times_YY'\to Y'$ the
projections as before. Then for any element $(Z,g)\in
\text{CH}_{\text{Ar}}(X)$, we have 
$g^*f_*(Z,g)=(f_g)_*(g_f)^*(Z,g)\in  \text{CH}_{\text{Ar}}(Y').$
\vskip 0.20cm
\noindent
(A.9) Let $X$ be regular projective variety over ${\Bbb C}$.
Then, following [GS2], we define the {\it arithmetic 
$K$-group} $K_{\text{Ar}}(X)$  as the abelian group generated by 
hermitian vector bundles $(E,\rho)$ on $X$ and $\eta\in \tilde
A(X)$ modulo  the subgroup generated by the following relations: 
For any short exact sequence of vector bundles on $X$,
$${E}.:\qquad 0\rightarrow E_1\rightarrow E_2\rightarrow
E_3\rightarrow 0,$$ let $\rho_i$ be 
hermitian metrics on  $E_i$, then
$$((E_1,\rho_1);\eta_1)+((E_3,\rho_3);\eta_3)=((E_2,\rho_2);
-\text{ch}_{\text{BC}} (E.,\rho_1,\rho_2,\rho_3)+\eta_1+\eta_3).$$ 
Here $\text{ch}_{\text{BC}}(E.;\rho.)$
denotes the classical Bott-Chern secondary characteristic 
classes associated with $(E.,\rho.)$ on $X$ with respect to ch.
(See e.g., Remark 1.1.1.)
\vskip 0.45cm
To motivate the definition of arithmetic characteristic
classes, we start with an example.  Let $(L,\rho)$ be a
hermitian line bundle on  $X$. Then for
any non zero section $s$ of $L$, by the
Poincar\'e-Lelong equation in Example A.2, 
$(\text{div}(s),-[\log |s|^2_\rho])$ is an element of 
$\text{CH}^1_{\text{Ar}}(X)$. Define
$c_{\text{Ar},1}(L,\rho)$ as the class of this element in the
arithmetic Chow ring.
From here, as one may imagine, in general, we may use the
splitting principle together with the classical Bott-Chern
secondary characteristic classes to construct  arithmetic
characteristic classes for  hermitian vector bundles.
\vskip 0.20cm
\noindent
(A.10) Let $B$ be a subring  of real number field ${\Bbb R}$,
and let $\phi\in B[[T_1,\dots,T_n]]$ be a symmetric power 
series. Then we have the following;
\vskip 0.20cm
\noindent
{\bf Theorem.} ([GS2]) {\it Associated to every hermitian
vector bundle $(E,\rho)$ of rank $n$ on $X$ is an
arithmetic characteristic class $\phi_{\text{Ar}}(E,\rho)\in
\text{CH}_{\text{Ar}}(X)_{\Bbb Q}$ which satisfies the
following properties:
\vskip 0.20cm
\noindent
{(i)} If $f:Y\to X$ is a morphism,
$f^*(\phi_{\text{Ar}}(E,\rho))=\phi_{\text{Ar}}(f^*E,f^*\rho).$
\vskip 0.20cm
\noindent
{(ii)}  If
$(E,\rho)=(L_1,\rho_1)\oplus\dots
\oplus (L_n,\rho_n)$ is an orthogonal direct sum of 
hermitian line bundles,
$\phi_{\text{Ar}}(E,\rho)=\phi(c_{{\text{Ar}},1}(L_1,\rho_1),\dots,c_{{\text{Ar}},1}(L_n,\rho_n)).$
\vskip 0.20cm
\noindent
{(iii)}  In $\tilde A(X)$, 
$\omega(\phi_{\text{Ar}}(E,\rho))=\phi(E,\rho),$ and in
$\text{CH}(X)_{\Bbb Q}$,
 $\zeta(\phi_{\text{Ar}}(E,\rho))=\phi(E).$
\vskip 0.20cm
\noindent
{(iv)} Let $\quad E.:\quad 0\rightarrow E_1\rightarrow E_2
\rightarrow E_3\rightarrow 0\quad$ be an exact sequence of vector 
bundles on $X$ together with
 hermitian metrics $\rho_i$ on $E_i$ for $i=1,2,3.$ Then
$$\phi_{\text{Ar}}(E_2,\rho_2)=\phi_{\text{Ar}}(E_1\oplus
E_3,\rho_1\oplus
\rho_3)+a(\phi_{\text{BC}}(E.,\rho_1,\rho_2,\rho_3)).$$}
\vskip 0.20cm
In particular, we have a well-defined arithmetic Chern
characteristic class
$$\text{ch}_{\text{Ar}}:K_{\text{Ar}}(X)_{\Bbb
Q}\to\text{CH}_{\text{Ar}}(X)_{\Bbb Q}.$$ It is a result of
Gillet and Soul\'e that $K_{\text{Ar}}(X)_{\Bbb
Q}$ admits a so-called $\lambda$-ring structure, and
$\text{ch}_{\text{Ar}}$ indeed gives a ring isomorphim. We
will not recall the details, instead, we state the following
equivalent:
\vskip 0.20cm
\noindent
{\bf Theorem}$'$. ([GS2]) {\it Let $X$ be a regular projective
variety over ${\Bbb C}$. Then
$$\text{ch}_{\text{Ar}}:K_{\text{Ar}}(X)_{\Bbb
Q}\to\text{CH}_{\text{Ar}}(X)_{\Bbb Q}$$ is a group
isomorphism. Moreover, 
$$\text{ch}_{\text{Ar}}
\Bigl((E,\rho)\otimes 
(E',\rho')\Bigr)
=\bigl(\text{ch}_{\text{Ar}}({E},\rho)\bigr)\cdot
\bigl(\text{ch}_{\text{Ar}}({E}',\rho')\bigr).$$}
\vskip 0.20cm
\noindent
(A.11) We end this section by recalling a result due to
Faltings.
\vskip 0.20cm
Let $f:X\to Y$ and $g: Z\to Y$ be two smooth, proper
morphisms of compact K\"ahler manifolds and $i:X\to Z$ be a
codimension one closed immersion over $Y$, i.e., $i$ is a closed
immersion of codinemsion one such that $f=g\circ i$. 
Denote by $$\pi: W:=B_{X\times \{\infty\}}Z\times {\Bbb P}^1
\rightarrow Z\times {\Bbb P}^1,$$ the natural projection, where
$B_{X\times \{\infty\}}Z\times {\Bbb P}^1$ denotes the blowing-up 
of $Z\times {\Bbb P}^1$ along $X\times \{\infty\}$.  
Denote the exceptional divisor by ${\Bbb P}$ of $\pi$. We know
that the map $q_W:W\rightarrow {\Bbb P}^1$, obtained by
composing $\pi$ with the projection
$q:Z\times {\Bbb P}^1\rightarrow {\Bbb P}^1$, is flat, and that 
for $t\in {\Bbb P}^1$:
$$q_W^{-1}(t)=\cases Z\times \{t\}, &$for $t\not=\infty$$,\\
{\Bbb P}\cup B_XZ,&$for $t=\infty$$.\endcases$$ Here $B_XZ$
denotes the blowing-up of
$Z$ along
$X$. Moreover, by the construction, ${\Bbb P}$ and $B_XZ$
intersect transversally, and 
${\Bbb P}\cap B_XZ$ is the exceptional divisor $X$ on
$B_XZ$.
\vskip 0.20cm
Denote by $I:X\times {\Bbb P}^1\hookrightarrow W$  the induced
codimension one closed embedding. Then the image of
$I$ does not intersect with $B_XZ$, and the image $X\times
\{\infty\}$ in $W$ is a section of ${\Bbb P}$.
Denote  the induced fibration $Z\times \{t\}\to Y\times\{t\}$
by $g_t$ for $t\not=\infty$ and  set $g_\infty$ to be the
composition of the projection of ${\Bbb P}$ on $X$ with
$(X=)X\times\{\infty\}\to Y\times\{\infty\}(=Y)$. Denote by $f_t:
X\times \{t\}\to Y\times \{t\}$ the smooth morphisms induced
from $f$ for all $t\in {\Bbb P}^1$. 
\vskip 0.20cm
Let $E$ be a  vector bundles on $Z$.
Then we have  the following
exact sequence of coherent sequences on $W$;
$$0\to (\pi\circ
p_Z)^*E(B_XZ-X\times
{\Bbb P}^1)\to (\pi\circ
p_Z)^*E(B_XZ)\to I_*I^* ((\pi\circ
p_Z)^*E(B_XZ))\to 0.\eqno(7.1)$$
\vskip 0.20cm  
Thus, by the flatness of $q_W:W\to {\Bbb P}^1$, we know that
the restriction of (7.1) to the fibers $W_t$ of $q_W$ for all
$t\in {\Bbb P}^1$ are  exact.
 Thus, in particular, for each
$t\not=\infty$ in
${\Bbb P}^1$, from (7.1), we have the induced
exact sequences
$$0\to E_t(-X)\to E_t\to (i_t)_*i_t^*E_t\to 0
\qquad\text{over}\quad Z\times\{t\}\eqno(7.2)$$
 with $E_t$ the pull-back of
$E$ under the canonical identity $X\times \{t\}\simeq X$.
Similarly, for the fiber at $\infty$, if we set
$E(B_XZ)\big|_{{\Bbb P}}:=E_\infty$, 
$E(B_XZ)\big|_{B_XZ}=:E_\infty'$ and 
$E(B_XZ-X\times {\Bbb
P}^1)\big|_{B_XZ}=:E_\infty''$. Then
$E(B_XZ-X\times {\Bbb P}^1)\big|_{{\Bbb P}}=E_\infty(-X)$, and
(7.1) splits into two exact sequences 
$$0\to E_\infty(-X)\to E_\infty\to
(i_\infty)_*i_\infty^*E_\infty\to 0\qquad\text{over}\quad
{\Bbb P}\eqno(7.3)$$ and $$0\to E_\infty'
\to E_\infty''\to 0\to 0 \qquad\text{over}\quad
B_XZ.\eqno(7.4)$$  In
particular, we see that on
$B_XZ$, $E_\infty'=E_\infty''$. 
\vskip 0.20cm
Choose a hermitian metric $\tau_W$ on $T_W$, the tangent
bundle of $W$. Then
$\tau_W$ naturally induces hermitian metrics
$\tau_t$ on $T_{g_t}$ for all $t\in {\Bbb P}^1$. Similarly,
$\tau_W$ induces a hermitian metric $\tau_G$ on 
$T_G(-\log\infty)$,
the logarithmic relative tangent bundle associated to the
morphism
$G:W\to Y\times {\Bbb P}^1$, which may be naturally embedded
in $T_W$ (see e.g. [De1]).
\vskip 0.20cm
Fix hermitian metrics $\rho$
and $\rho'$ on $E$ and on $E(-X)$ respectively. Use the same
notation to denote the pull-back of $(E,\rho)$ onto $W$. 
Choose the Fubini-Study metric on
${\Cal O}_{{\Bbb P}^1}(\infty)$ and a metric on ${\Cal
O}_{W}(-X\times {\Bbb P}^1)$ such that in a neighborhood $U$ of
$B_XZ$, which is away from
$X\times {\Bbb P}^1$, the natural isomorphism ${\Cal
O}_{W}(-X\times {\Bbb P}^1)\simeq {\Cal O}_W$ induces an
isometry, once we put the trivial metric on ${\Cal O}_W$. Denote
these final induced metrics on
$E(B_XZ)$ and $E(B_XZ-X\times {\Bbb P}^1)$ by $D\rho$ and
$D\rho'$ respectively. 
\vskip 0.20cm
Denote the induced metrics via restriction to $E_t$ and $E_t(-X)$
(resp. to $E_\infty'$ and $E_\infty''$ on $B_XE$) by $\rho_t$ and
$\rho_t'$ respectively for all
$t\in {\Bbb P}^1$, (resp. $\rho_\infty''$ and $\rho_\infty'''$).
Easily, we see that $\rho_0=\rho$ and $\rho_0'=\rho'$, and 
$(E_\infty',\rho_\infty'')$ is isomorphic to
$(E_\infty'',\rho_\infty''')$ by the
construction.  In this way, $$(E(-X),\rho')\hookrightarrow
(E,\rho)\qquad\text{on}\quad Z$$ is deformed to
$$(E_\infty(-X),\rho_\infty')\hookrightarrow
(E_\infty,\rho_\infty)\quad \text{on}\quad {\Bbb P}$$ (and 
$$(E_\infty',\rho_\infty'')\simeq(E_\infty'',\rho_\infty''')\qquad
\text{on}\quad
B_XZ).$$ 
\vskip 0.20cm
\noindent
{\bf Proposition.} ([F]) {\it With the same notation as above,
for all $t\in {\Bbb P}^1$,
$$\eqalign{~&\Big(\text{ch}_{\text{Ar}}(E_t,\rho_t)-
\text{ch}_{\text{Ar}}(E_t(-X),\rho_t')\Big)
\cdot
\text{td}_{\text{Ar}}(T_t,\tau_{g_t})\cr
=&i_t^*\Big(\Big(\text{ch}_{\text{Ar}}(E,D\rho)-\text{ch}_{\text{Ar}}(E,D\rho)\Big)
\cdot\text{td}_{\text{Ar}}(T_G(-{\log}\infty),\tau_G)\Big).\cr}$$}
\noindent
{\it Proof.} Surely, by (A.8), $i_t^*$ is
well-defined. So at least all terms make perfect sense as they
stand.
\vskip 0.20cm
The proof is rather formal but standard. The key points are  that
$E-E(-X\times {\Bbb P}^1)$ is supported only on $X\times {\Bbb
P}^1$; that $X\times {\Bbb P}^1$ does not intersect $B_XZ$ where
the hermitian exact sequence splits. As a direct consequence, in
the calculation, we may only pay our attention on $X\times
{\Bbb P}^1$, while pay no attention on $B_XZ$. With this in
mind, certainly,  we may also simply view 
$(T_{0},\tau_{g_0})$ (resp.
$(T_{\infty},\tau_{g_\infty})$) as the restriction of the
$(T_G(-{\log}\infty),\tau_G)$. (Recall that the induced metric
on $T_\infty$ from $\tau_G$ is singular, but the singularity
is concentrated on ${\Bbb P}\cap B_XZ$, which is away from
$X\times {\Bbb P}^1$ by the construction.) 
To be more precise, in practice, we need the following
preparation.
\vskip 0.20cm
First, following [F], we give a generalization of arithmetic
intersection  introduced in A.7. Instead of working over a single
regular projective variety over ${\Bbb C}$, we  are now working
over a triple
$(A;B;C)$. Here $A$ and $C$ are  regular
projective varieties over ${\Bbb C}$, $A\subset B\subset C$ with $B$ an open subset of
$C$, i.e., $B$ an open neighborhood of $A$ in $C$.
Define a relative arithmetic Chow group 
$\text{CH}_{\text{Ar}}^{A, B}(C)$ by setting it to be 
the quotient of the group generated by  arithmetic  cycles
$(S,g_S)$ of $C$ with 
$S\subset A$, $\text{Supp}(g_S)\subset B$, 
modulo the subgroup generated by arithmetic cycles defined by rational functions on cycles 
in $A$, together with the forms $\partial \alpha+\bar\partial\beta,$ where
$\alpha$ and $\beta$ are currents with support in $A$. We point out that even by
tensoring with ${\Bbb Q}$, the resulting space
$\text{CH}_{\text{Ar}}^{A, B}(C)_{\Bbb Q}$ only admits a group structure.
\vskip 0.20cm
Nevertheless, we may introduce a natural
$\text{CH}_{\text{Ar}}(C)_{\Bbb Q}$-module structure on 
$\text{CH}_{\text{Ar}}^{A, B}(C)_{\Bbb Q}$. Indeed, for an element $\omega\in \tilde A(C)$,
it is clear that for any element $(S,g_S)\in \text{CH}_{\text{Ar}}^{A, B}(C)_{\Bbb Q}$, the
standard arithmetic intersection on $C$ gives an arthmetic cycle $(0,\omega)\cdot(S,g_S)$
on $C$ which then turns out to be an element in $\text{CH}_{\text{Ar}}^{A, B}(C)_{\Bbb Q}$
as well. Similarly, for any hermitian line bundle $(L,\rho)$ on $C$,
 Theorem A.7, in which  the standard arithmetic
intersection is introduced, shows that
$c_{1,\text{Ar}}(L,\rho)\cdot (S,g_S):=(\text{div}(s),-\log\|s\|_\rho^2)\cdot(S,g_S)$ 
is again an element in $\text{CH}_{\text{Ar}}^{A, B}(C)_{\Bbb Q}$. Here $s$ is a non-zero
section of $L$. Easily, one checks that all the actions satisfies the module axioms.
Thus in particular, by Theorem A.10$'$ and splitting
principle, we have provd the first part of the following
\vskip 0.20cm
\noindent
{\bf Lemma.} ([F]) {\it With the same notation as above,
\vskip 0.20cm
\noindent
(i) $\text{CH}_{\text{Ar}}^{A, B}(C)_{\Bbb Q}$ admits a natural
$\text{CH}_{\text{Ar}}(C)_{\Bbb Q}$-module structure;
\vskip 0.20cm
\noindent
(ii) If $(E_1,\rho_1)$ and $(E_2,\rho_2)$ are two hermitian vector bundles on $C$ such that
over $B$, there is an isometry $(E_1,\rho_1)\Big|_B\simeq (E_2,\rho_2)\Big|_B$, then for
any arithmetic characteristic class $\phi_{\text{Ar}}$, and any element $(S,g_S)\in
\text{CH}_{\text{Ar}}^{A, B}(C)$, we have in $\text{CH}_{\text{Ar}}^{A, B}(C)_{\Bbb Q}$ and
hence in $\text{CH}_{\text{Ar}}(A)$, $$\phi_{\text{Ar}}(E_1,\rho_1)\cdot (S,g_S)=
\phi_{\text{Ar}}(E_2,\rho_2)\cdot (S,g_S).$$}
\vskip 0.20cm
\noindent
{\it Proof.} Let us start with the case for hermitian line bundles. In this case, 
 we have the following
\vskip 0.20cm
\noindent
{\bf Lemma.}$'$ {\it Let $(A_1;B_1;C_1)$ and $(A_2;B_2;C_2)$ 
be two triples as above such that
there is an isomorphism $a:B_1\simeq B_2$ which induces an isomorphism 
$A_1\simeq A_2$. Then
\vskip 0.20cm
\noindent
(i) there is a natural induced isomorphism
$$\text{CH}_{\text{Ar}}^{A_1, B_1}(C_1)\simeq \text{CH}_{\text{Ar}}^{A_2, B_2}(C_2);$$
\vskip 0.20cm
\noindent
(ii) If there are hermitian line bundles $(L_1,\rho_1)$ and $(L_2,\rho_2)$ on $Z_1$ and
$Z_2$ respectively such that there is an isometry
$(L_1,\rho_1)\Big|_{B_1}\simeq (L_2,\rho_2)\Big|_{B_2},$ then we have the following
commutative diagram;
$$\matrix \text{CH}_{\text{Ar}}^{A_1, B_1}(C_1)&\simeq&\text{CH}_{\text{Ar}}^{A_2,
B_2}(C_2)\\
&&\\
~_{c_{1,\text{Ar}}(L_1,\rho_1)\cdot(\cdot)}\downarrow&&\downarrow
~_{(\cdot)\cdot c_{1,\text{Ar}}(L_2,\rho_2)}\\
 &&\\ 
\text{CH}_{\text{Ar}}^{A_1,
B_1}(C_1)&\simeq&\text{CH}_{\text{Ar}}^{A_2,
B_2}(C_2).\endmatrix$$}
\vskip 0.20cm
\noindent
{\it Proof.} This comes  directly from the definition.
\vskip 0.3cm
Let us now go back to the proof of the lemma. For doing so, we use the flag varieties of
$E_1$ and $E_2$.  Denote by
$\text{Flag}_CE_1$ (resp. $\text{Flag}_CE_2$) the complete flag
variety of $E_1$ (resp. $E_2$), by
$\pi_i:\text{Flag}_CE_i\rightarrow Z$ the natural projections, and
$B_i:=\pi_i^{-1}(B)$, and $A_i:=\pi_i^{-1}(A)$. By our assumption
on $(E_i,\rho_i)$ and Lemma$'$(i), we have an isomorphism
$$\text{CH}_{\text{Ar}}^{A_1, B_1}(\text{Flag}_C{E_1})\simeq
\text{CH}_{\text{Ar}}^{A_2, B_2}(\text{Flag}_C{E_2})$$ with which
the pull-back of arithmetic cycles $\pi_i^*(S,g_S)$'s are
compactible.
\vskip 0.20cm
But for $i=1,2$, on $\text{Flag}_C{E_i}$, $\pi_i^*(E_i)$ has a
complete filtration by line bundles $L_{i,\alpha}$. Thus we can
let any polynomial in the
$c_{1,\text{Ar}}(L_{i,\alpha},\rho_{i,\alpha})$ act on the arithmetic cycles
in $\text{CH}_{\text{Ar}}^{A_i, B_i}(\text{Flag}_C{E_i})$. By
Lemma$'$ (ii), at this level, the actions for $i=1$ and for $i=2$
correspond to each other in a unique way. Easily, one sees that
so do the actions for elements in $\tilde A$, which in particular
includes all classical Bott-Chern secondary characteristic
classes $\text{ch}_{\text{BC}}$ measuring the change of
arithmetic characteristic classes for hermitian vector bundles
with respect to the extension of vector bundles and the change of
their metrics. Thus, by Theorem A.10, for
$\phi_{\text{Ar}}(E_i,\rho_i)$, there exists an operator
$Q(c_{1,\text{Ar}}(L_{i,\alpha},\rho_{i,\alpha}))+\Omega_{\text{BC}}$ on
$\text{CH}_{\text{Ar}}(\text{Flag}_C{E_i})_{\Bbb Q}$, where $Q$
is a power series in 
$c_{1,\text{Ar}}(L_{i,\alpha},\rho_{i,\alpha})$ and $\Omega_{\text{BC}}$ the multiplication
by a classical Bott-Chern secondary characteristic class, such
that for all $(S,g_S)\in
\text{CH}_{\text{Ar}}^{A,B}(C)$,
$$\phi_{\text{Ar}}(E_i,\rho_i)\,(S,g_S)
=\pi_{i*}\Big(\Bigl(\bigl(Q(c_{1,{\text{Ar}}}(L_{i,\alpha},\rho_{i,\alpha}))\bigr)
+\Omega_{\text{BC}}\Bigr)\,\pi^*_i(S,g_S)\Big).$$
But, by Lemma$'$(ii), we know that
$$\eqalign{~&\pi_{1*}\Big(\Bigl(\bigl(Q(c_{1,{\text{Ar}}}(L_{1,\alpha},\rho_{1,\alpha}))\bigr)
+\Omega_{\text{BC}}\Bigr)\,\pi^*_1(S,g_S)\Big)\cr
=&\pi_{2*}\Big(\Bigl(\bigl(Q(c_{1,{\text{Ar}}}(L_{2,\alpha},\rho_{2,\alpha}))\bigr)
+\Omega_{\text{BC}}\Bigr)\,\pi^*_2(S,g_S)\Big),\cr}$$
via the isomorphism 
$$\text{CH}_{\text{Ar}}^{A_1, B_1}(\text{Flag}_C{E_1})\simeq
\text{CH}_{\text{Ar}}^{A_2, B_2}(\text{Flag}_C{E_2}).$$
Thus in $\text{CH}_{\text{Ar}}^{A, B}(C)_{\Bbb Q}$, and hence in
$\text{CH}_{\text{Ar}}(C)_{\Bbb Q}$, we have
$$\phi_{\text{Ar}}(E_1,\rho_1)\,(S,g_S)=\phi_{\text{Ar}}(E_2,\rho_2)\,(S,g_S)$$
for all elements $(S,g_S)\in \text{CH}_{\text{Ar}}^{A, B}(C).$
This then completes the proof of the lemma.
\vskip 0.20cm
With this lemma, now the proof of the proposition can be easily deduced. Indeed, first, 
since $E-E(-X\times {\Bbb P}^1)$ is supported only on 
$X\times {\Bbb P}^1$ and $X\times {\Bbb P}^1$ does not intersect $B_XZ$ where
the hermitian exact sequence splits,
$\text{ch}_{\text{Ar}}(E-E(-X\times {\Bbb P}^1),D\rho-D\rho')\in
\text{CH}_{\text{Ar}}^{X\times{\Bbb P}^1;W\backslash\overline U}(W)_{\Bbb Q}$.  
Secondly, we easily have the  functorial
properties for the  relative arithmetic interesction defined
above. Thus, we may apply the lemma as follows: At the fiber over
0, we choose $(A;B;C)$ to be
$(X;W\backslash\overline U\cap Z;Z)$, 
$(S,g_S)$ to be $\text{ch}_{\text{Ar}}(E-E(-X\times {\Bbb P}^1),D\rho-D\rho')\Big|_Z$, 
$\phi_{\text{Ar}}$ to be $\text{td}_{\text{Ar}}$, and $(E_1,\rho_1)$ to be $(T_g,\tau_g)$,
$(E_2,\rho_2)$ to be 
$(T_G(-\log\infty),\tau_G)\Big|_Z$; while 
at the fiber over $\infty$ or better over ${\Bbb P}$, we choose $(A;B;C)$ to be
$(X;W\backslash\overline U\cap {\Bbb P};{\Bbb P})$, 
$(S,g_S)$ to be $(E-E(-X\times {\Bbb P}^1),D\rho-D\rho')\Big|_{\Bbb P}$,  $\phi_{\text{Ar}}$
to be $\text{td}_{\text{Ar}}$, and $(E_1,\rho_1)$ to be $(T_{\infty},\tau_{g_\infty})$,
$(E_2,\rho_2)$ to be  $(T_G(-\log\infty),\tau_G)\Big|_{\Bbb P}$. 
All this then  completes the proof of the proposition.
\vskip 0.20cm
\noindent
{\bf B. An Effective Construction of Relative Bott-Chern
Secondary Characteristic Classes}
\vskip 0.20cm
\noindent
(B.1) In (A.6), we see that there is an exact sequence
$$\tilde A^{p-1,p-1}(X)
\buildrel a\over\rightarrow \text{CH}^p_{\text{Ar}}(X)\buildrel\zeta\over\rightarrow \text{CH}^p(X)
\rightarrow 0.\eqno(7.5)$$ Therefore the $a$-image of $\tilde
A^{p-1,p-1}(X)$ in $\text{CH}^p_{\text{Ar}}(X)$ is a
well-defined space. Due to the fact that  
if $f:Z\to X$ is a morphism of
regular projective varieties over ${\Bbb C}$, then the maps
appeared in the above exact sequence are compactible with all
the pull-back $f^*$ for the corresponding spaces, and hence
$a(\tilde A^{p-1,p-1}(X))$, a quotient space of $A^{p-1,p-1}(X)$,
is as canonical as $A^{p-1,p-1}(X)$ itself. 
In particular, we may equally use this space, denoted by
$\overline A^{p-1,p-1}(X)$,  to develop the theory of
relative Bott-Chern secondary characteristic classes. That is to
say, in the axioms stated for relative Bott-Chernsecondary
characteristic classes, we may use $\overline A$,
instead of the original $\tilde A$, e.g., we consider the
classical and relative Bott-Chern secondary characteristic classes
as elements in $\overline A$ rather than as in $\tilde A$. One
checks easily that the previous uniqueness theorems and their
proofs work exactly the same way as before. 
\vskip 0.20cm
\noindent
(B.2) We now give an effective construction for relative
Bott-Chern secondary characteristic classes
$\text{ch}_{\text{BC}}$ at the level of
$\overline A:=\oplus_p \overline A^{p,p}$.
That  amounts to saying that for all smooth metrized morphisms
and all relative acyclic hermitian vector bundles, we should
effectively construct some classes in
$\overline A$, which satisfy the
(corresponding modified) six axioms for relative Bott-Chern
secondary characteristic classes.
\vskip 0.20cm
Let $(f:X\to Y;E,\rho;T_f,\tau_f)$ be a properly metrized
datum. Moreover, from now on, we always assume that all
manifolds are regular projective over ${\Bbb C}$. 
Then on the direct image, we have
a well-defined
$L^2$-metric
$L^2(\rho,\tau)$. (Please note that the K\"ahler conditions in
the definition of a properly metrized datum for $\tau_f$
guarantee that the
$L^2$-metric is well-defined via the Hodge decomposition.) 
\vskip 0.20cm
\noindent
{\bf Lemma.} {\it With the same notation as above, 
$$f_*\Big(\text{ch}_{\text{Ar}}(E,\rho)\cdot
\text{td}_{\text{Ar}}(T_f,\tau_f)\Big)
-\text{ch}_{\text{Ar}}\Big(f_*E,L^2(\rho,\tau_f)\Big)\in
\overline A(Y).$$}
\noindent
{\it Proof.}  In fact,
$$\eqalign{~&\zeta\Big(f_*\Big(\text{ch}_{\text{Ar}}(E,\rho)\cdot
\text{td}_{\text{Ar}}(T_f,\tau_f)\Big)
-\text{ch}_{\text{Ar}}\Big(f_*E,L^2(\rho,\tau_f)\Big)\Big)\cr
=&f_*\Big(\zeta\Big(\text{ch}_{\text{Ar}}(E,\rho)\cdot
\text{td}_{\text{Ar}}(T_f,\tau_f)\Big)\Big)
-\zeta\Big(\text{ch}_{\text{Ar}}\Big(f_*E,L^2(\rho,\tau_f)\Big)\Big)\cr
=&f_*\Big(\zeta\Big(\text{ch}_{\text{Ar}}(E,\rho)\Big)\cdot
\zeta\Big(\text{td}_{\text{Ar}}(T_f,\tau_f)\Big)
-\zeta\Big(\text{ch}_{\text{Ar}}\Big(f_*E,L^2(\rho,\tau_f)\Big)\Big)\cr
&\qquad(\text{by\ Theorem A.7(iii)})\cr
=&f_*\Big(\text{ch}(E)\cdot
\text{td}(T_f)\Big)
-\text{ch}\Big(f_*E)\Big)\qquad(\text{By\ Theorem\ A.10(iii)})\cr
=&0\qquad(\text{by\ the\ Grothendieck-Riemann-Roch\ Theorem\ in
Algebraic\ Geomretry}).\cr}$$
Hence, by the exact sequence (7.5), we complete the proof
of the lemma.
\vskip 0.20cm
Now set
$$\eqalign{~&\text{ch}_{\text{BC}}'(E,\rho;f,\tau_f)\cr
:=&
f_*\Big(\text{ch}_{\text{Ar}}(E,\rho)\cdot
\text{td}_{\text{Ar}}(T_f,\tau_f)\Big)
-\text{ch}_{\text{Ar}}\Big(f_*E,L^2(\rho,\tau_f)\Big)
\in \overline A(Y).\cr}$$
We claim that $\text{ch}_{\text{BC}}'(E,\rho;f,\tau_f)$ gives an
effective construction of the relative Bott-Chern secondary
characteristic classes. So it
remains to check that $\text{ch}_{\text{BC}}'(E,\rho;f,\tau_f)$
satisfies six axioms in Chapters 2 and 3.
\vskip 0.20cm
\noindent
(B.3) First, for Axiom 1,  the downstairs rule, we have the
following;
\vskip 0.20cm
\noindent
{\bf Lemma.} {\it With the same notation as above,
$$dd^c\Big(\text{ch}_{\text{BC}}'(E,\rho;f,\tau_f)\Big)=
f_*\Big(\text{ch}(E,\rho)\cdot
\text{td}(T_f,\tau_f)\Big)
-\text{ch}\Big(f_*E,L^2(\rho,\tau_f)\Big).$$}
\noindent
{\it Proof.} This may
be obtained by the following calculation. Indeed,
$$\eqalign{~&dd^c\Big(\text{ch}_{\text{BC}}'
(E,\rho;f,\tau_f)\Big)\cr
=&\omega\Big(\text{ch}_{\text{BC}}'
(E,\rho;f,\tau_f)\Big)\quad(\text{as}\ \text{ch}_{\text{BC}}
(E,\rho;f,\tau_f)\in \overline A(Y))\cr
=&\omega\Big(f_*\Big(\text{ch}_{\text{Ar}}(E,\rho)\cdot
\text{td}_{\text{Ar}}(T_f,\tau_f)\Big)
-\text{ch}_{\text{Ar}}\Big(f_*E,L^2(\rho,\tau_f)\Big)\Big)\cr
&\qquad(\text{by\ definition})\cr
=&f_*\Big(\omega\Big(\text{ch}_{\text{Ar}}(E,\rho)\cdot
\text{td}_{\text{Ar}}(T_f,\tau_f)\Big)\Big)
-\omega\Big(\text{ch}_{\text{Ar}}\Big(f_*E,L^2(\rho,\tau_f)\Big)\Big)\cr
=&f_*\Big(\omega\Big(\text{ch}_{\text{Ar}}(E,\rho)\Big)\cdot
\omega\Big(\text{td}_{\text{Ar}}(T_f,\tau_f)\Big)\Big)
-\omega\Big(\text{ch}_{\text{Ar}}\Big(f_*E,L^2(\rho,\tau_f)\Big)\Big)\cr
&\qquad(\text{by\ Theorem\
A.7(iii)})\cr
=&f_*\Big(\text{ch}(E,\rho)\cdot
\text{td}(T_f,\tau_f)\Big)
-\text{ch}\Big(f_*E,L^2(\rho,\tau_f)\Big)
\qquad(\text{by\ Theorem\ A.10(iii)}).\cr}$$ This completes the
proof.
\vskip 0.20cm
\noindent
(B.4) For Axiom 2, we should prove  the
following
\vskip 0.20cm
\noindent
{\bf Lemma.} {\it With the same notation as above, let
$(F,\rho')$ be a hermitian vector bundle on $Y$, then
$$\text{ch}_{\text{BC}}'(E\otimes f^*F,\rho\otimes
f^*\rho';f,\tau_f)=\text{ch}_{\text{BC}}'(E,\rho;f,\tau_f)
\wedge\text{ch}(F,\rho').$$}
\noindent
{\it Proof.} In fact,
$$\eqalign{~&\text{ch}_{\text{BC}}'(E\otimes f^*F,\rho\otimes
f^*\rho';f,\tau_f)\cr
=&f_*\Big(\text{ch}_{\text{Ar}}(E\otimes f^*F,\rho\otimes
f^*\rho')\cdot
\text{td}_{\text{Ar}}(T_f,\tau_f)\Big)
-\text{ch}_{\text{Ar}}\Big(f_*(E\otimes
f^*F),L^2(\rho\otimes f^*\rho',\tau_f)\Big)\Big)\cr
&\qquad(\text{by\ definition})\cr
=&f_*\Big(\text{ch}_{\text{Ar}}(E,\rho)\cdot
f^*(\text{ch}_{\text{Ar}}(F,\rho'))\cdot
\text{td}_{\text{Ar}}(T_f,\tau_f)\Big)
-\text{ch}_{\text{Ar}}\Big((f_*(E),L^2(\rho,\tau_f))
\otimes (F,\rho')\Big)\cr
&\qquad(\text{by\ Theorem\ A.10}'\ \text{and\ Theorem A.10(i)})\cr
=&f_*\Big(\text{ch}_{\text{Ar}}(E,\rho)
\cdot
\text{td}_{\text{Ar}}(T_f,\tau_f)\Big)\cdot
\text{ch}_{\text{Ar}}(F,\rho')
-\text{ch}_{\text{Ar}}\Big(f_*(E,L^2(\rho,\tau_f)
\Big)\cdot
\text{ch}_{\text{Ar}}(F,\rho')\cr
&\qquad(\text{by\ projection\ formula\ and\ Theorem\ A.10}')\cr
=&\Big(f_*\Big(\text{ch}_{\text{Ar}}(E,\rho)
\cdot
\text{td}_{\text{Ar}}(T_f,\tau_f)\Big)
-\text{ch}_{\text{Ar}}\Big(f_*(E,L^2(\rho,\tau_f)
\Big)\Big)\cdot
\text{ch}_{\text{Ar}}(F,\rho')\cr
=&\text{ch}_{\text{BC}}'(E,\rho;f,\tau_f)
\wedge\omega\Big(\text{ch}_{\text{Ar}}(F,\rho')\Big)
\qquad(\text{by\ Theorem A.7(iii)})\cr
=&\text{ch}_{\text{BC}}'(E,\rho;f,\tau_f)
\wedge\text{ch}(F,\rho')
\qquad(\text{by\ Theorem A.10(iii)}).\cr}$$ This completes the
proof.
\vskip 0.20cm
\noindent
(B.5) For Axiom 3, we should prove  the
following
\vskip 0.20cm
\noindent
{\bf Lemma.} {\it With the same notation as above, if $g:Y'\to
Y$ is a morphism between regular projective varieties over
${\Bbb C}$, then
$$g_f^*\Big(\text{ch}_{\text{BC}}'(E,\rho;f,\tau_f)\Big)=
\text{ch}_{\text{BC}}'\Big(g_f^*(E,\rho;f,\tau_f)\Big).$$
Here as in 2.C.2, we have the base change diagram $$\matrix
X\times_YY'&\buildrel g_f\over\to&X\\ &&\\
f_g\downarrow&&\downarrow f\\
&&\\
Y'&\buildrel g\over\to&Y,\endmatrix$$ and define the pull-back
of
$(E,\rho;f,\tau_f)$ via $g_f$ by
$$g_f^*(E,\rho;f,\tau_f)=(g_f^*(E,\rho);f_g,g_f^*\tau_f).$$}
\noindent
{\it Proof.} In fact, By Theorem A.10(i), the functorial property
for arithmetic characteristic classes and A.8, 
$$\eqalign{~&g^*f_*\Big(\text{ch}_{\text{Ar}}(E,\rho)\cdot
\text{td}_{\text{Ar}}(T_f,\tau_f)\Big)\cr
=&(f_g)_*(g_f)^*\Big(\text{ch}_{\text{Ar}}(E,\rho)\cdot
\text{td}_{\text{Ar}}(T_f,\tau_f)\Big)
\quad(\text{by\ A.8})\cr
=&(f_g)_*\Big(\text{ch}_{\text{Ar}}((g_f)^*(E,\rho))\cdot
\text{td}_{\text{Ar}}((g_f)^*(T_f,\tau_f))\Big)
\quad(\text{by\ Theorem A.10(i)})\cr
=&(f_g)_*\Big(\text{ch}_{\text{Ar}}((g_f)^*(E,\rho))\cdot
\text{td}_{\text{Ar}}(T_{f_g},\tau_{f_g}))\Big)
\quad(\text{by\ definition}).\cr}$$ So by definition and the
functorial property for arithmetic characteristic classes, it
suffices to show that $$g^*\Big(f_*E,L^2(\rho,\tau_f)\Big)\simeq
\Big((f_g)_*g_f^*E,L^2(g_f^*\rho,\tau_{f_g})\Big).$$
But this is  a direct consequence of our $f$-acyclic condition
on $E$. So Axiom 3 is checked.
\vskip 0.20cm
\noindent
(B.6) Now we check Axion 4. This amounts to proving the following
\vskip 0.20cm
\noindent
{\bf Lemma.} {\it With the same notation as above, let $0\to
E_1\to E_2\to E_3\to 0$ be an exact sequence of $f$-acyclic vector
bundles on
$X$. Then for hermitian metrics $\rho_i$ on $E_i$, $i=1,2,3$,
$$\eqalign{~&\text{ch}_{\text{BC}}'(E_2,\rho_2;f,\tau_f)
-\text{ch}_{\text{BC}}'(E_1,\rho_1;f,\tau_f)
-\text{ch}_{\text{BC}}'(E_3,\rho_3;f,\tau_f)\cr
=&f_*\Big(\text{ch}_{\text{BC}}(E.,\rho.)\cdot\text{td}(T_f,\tau_f)\Big)
-\text{ch}_{\text{BC}}\Big(f_*E.,L^2(\rho.,\tau_f)\Big).\cr}$$}
\vskip 0.30cm
\noindent
{\it Proof.} By definition, 
$$\eqalign{~&\text{ch}_{\text{BC}}'(E_2,\rho_2;f,\tau_f)
-\text{ch}_{\text{BC}}'(E_1,\rho_1;f,\tau_f)
-\text{ch}_{\text{BC}}'(E_3,\rho_3;f,\tau_f)\cr
=&f_*\Big(\text{ch}_{\text{Ar}}(E_2,\rho_2)\cdot\text{td}_{\text{Ar}}
(T_f,\tau_f)\Big)-\text{ch}_{\text{Ar}}\Big(f_*E_2,L^2(\rho_2,\tau_f)\Big)\cr
&-\Big(f_*\Big(\text{ch}_{\text{Ar}}(E_1,\rho_2)\cdot\text{td}_{\text{Ar}}
(T_f,\tau_f)\Big)-\text{ch}_{\text{Ar}}\Big(f_*E_1,L^2(\rho_1,\tau_f)\Big)\cr
&+f_*\Big(\text{ch}_{\text{Ar}}(E_3,\rho_3)\cdot\text{td}_{\text{Ar}}
(T_f,\tau_f)\Big)-\text{ch}_{\text{Ar}}\Big(f_*E_3,L^2(\rho_3,\tau_f)\Big)\Big)\cr
=&f_*\Big(\Big(\text{ch}_{\text{Ar}}(E_2,\rho_2)-\text{ch}_{\text{Ar}}(E_1,\rho_1)
-\text{ch}_{\text{Ar}}(E_3,\rho_3)\Big)\cdot\text{td}_{\text{Ar}}
(T_f,\tau_f)\Big)\cr
&-\Big(\text{ch}_{\text{Ar}}\Big(f_*E_2,L^2(\rho_2,\tau_f)\Big)
-\text{ch}_{\text{Ar}}\Big(f_*E_1,L^2(\rho_1,\tau_f)\Big)
-\text{ch}_{\text{Ar}}\Big(f_*E_3,L^2(\rho_3,\tau_f)\Big)\Big)\cr
=&f_*\Big(a\Big(\text{ch}_{\text{BC}}(E.,\rho.)\Big)
\cdot\text{td}_{\text{Ar}}
(T_f,\tau_f)\Big)-a\Big(\text{ch}_{\text{BC}}\Big(f_*E.,L^2(\rho.,\tau_f)\Big)\Big)\cr
&\qquad(\text{by\ Theorem\ 10(iv)})\cr
=&f_*\Big(a\Big(\text{ch}_{\text{BC}}(E.,\rho.)
\cdot\text{td}(T_f,\tau_f)\Big)\Big)-a\Big(\text{ch}_{\text{BC}}\Big(f_*E.,L^2(\rho.,\tau_f)\Big)\Big)\cr
&\qquad(\text{by\ Theorem\ A.7\ and\ Theorem\ 10(iii)})\cr
=&f_*\Big(\text{ch}_{\text{BC}}(E.,\rho.)
\cdot\text{td}(T_f,\tau_f)\Big)-\text{ch}_{\text{BC}}\Big(f_*E.,L^2(\rho.,\tau_f)\Big)\cr
&\qquad(\text{by\ Theorem\ A.7(iii)\ and\ (A.8)}).\cr}$$
This completes the proof of the lemma.
\vskip 0.30cm
\noindent
(B.7) Similarly, for Axiom 5, we need to show the following
\vskip 0.20cm
\noindent
{\bf Lemma.} {\it With the same notation as in section 2.E, then
$$\eqalign{~
&\text{ch}_{\text{BC}}'(E,\rho;g\circ f,\tau_{g\circ
f})-g_*\Big(\text{ch}_{\text{BC}}'(E,\rho;f,\tau_f)\cdot
\text{td}(T_g,\tau_g)\Big)
-\text{ch}_{\text{BC}}'\Big(f_*E,L^2(\rho,\tau_f);g,\tau_g\Big)\cr
&=(g\circ f)_*\Big(\text{ch}(E,\rho)\cdot
\text{td}_{\text{BC}}(T.,\tau.)\Big)-
\text{ch}_{\text{BC}}\Big((g\circ f)_*E;L^2(\rho,\tau_{g\circ
f}),L^2(L^2(\rho,\tau_f),\tau_g)\Big).\cr}$$}
\vskip 0.30cm
\noindent
{\it Proof.} By definition,
$$\eqalign{~
&\text{ch}_{\text{BC}}'(E,\rho;g\circ f,\tau_{g\circ
f})-g_*\Big(\text{ch}_{\text{BC}}'(E,\rho;f,\tau_f)\cdot
\text{td}(T_g,\tau_g)\Big)
-\text{ch}_{\text{BC}}'\Big(f_*E,L^2(\rho,\tau_f);g,\tau_g\Big)\cr
=&(g\circ f)_*\Big(\text{ch}_{\text{Ar}}(E,\rho)
\cdot\text{td}_{\text{Ar}}(T_{g\circ f},\tau_{g\circ
f})\Big)-\text{ch}_{\text{Ar}}\Big((g\circ
f)_*E,L^2(\rho,\tau_{g\circ f})\Big)\cr
&-g_*\Big(\Big(f_*\Big(\text{ch}_{\text{Ar}}(E,\rho)\cdot\text{td}_{\text{Ar}}
(T_f,\tau_f)\Big)-\text{ch}_{\text{Ar}}\Big(f_*E,L^2(\rho,\tau_f)\Big)\Big)\cdot
\text{td}(T_g,\tau_g)\Big)\cr
&-\Big(g_*\Big(\text{ch}_{\text{Ar}}(f_*E,L^2(\rho,\tau_f))
\cdot\text{td}_{\text{Ar}}
(T_g,\tau_g)\Big)-\text{ch}_{\text{Ar}}\Big(g_*(f_*E),L^2(L^2(\rho,\tau_f),\tau_g))\Big)\Big)\cr
=&(g\circ f)_*\Big(\text{ch}_{\text{Ar}}(E,\rho)
\cdot\text{td}_{\text{Ar}}(T_{g\circ f},\tau_{g\circ
f})\Big)-\text{ch}_{\text{Ar}}\Big((g\circ
f)_*E,L^2(\rho,\tau_{g\circ f})\Big)\cr
&-g_*\Big(\Big(f_*\Big(\text{ch}_{\text{Ar}}(E,\rho)\cdot\text{td}_{\text{Ar}}
(T_f,\tau_f)\Big)-\text{ch}_{\text{Ar}}\Big(f_*E,L^2(\rho,\tau_f)\Big)\Big)\cdot
\text{td}_{\text{Ar}}(T_g,\tau_g)\Big)\cr
&-\Big(g_*\Big(\text{ch}_{\text{Ar}}(f_*E,L^2(\rho,\tau_f))
\cdot\text{td}_{\text{Ar}}
(T_g,\tau_g)\Big)-\text{ch}_{\text{Ar}}\Big(g_*(f_*E),L^2(L^2(\rho,\tau_f),\tau_g))\Big)\Big)\cr
&\qquad(\text{by Theorem A.7(iii)\ and\ Theorem A.10(iii)}\cr
&\qquad\text{as}\
f_*\Big(\text{ch}_{\text{Ar}}(E,\rho)\cdot\text{td}_{\text{Ar}}
(T_f,\tau_f)\Big)-\text{ch}_{\text{Ar}}\Big(f_*E,L^2(\rho,\tau_f)\Big)\in
\bar A(Y))\cr
=&(g\circ f)_*\Big(\text{ch}_{\text{Ar}}(E,\rho)
\cdot\Big(\text{td}_{\text{Ar}}(T_{g\circ f},\tau_{g\circ
f})-\text{td}_{\text{Ar}}
(T_f,\tau_f)\cdot f^*\Big(
\text{td}_{\text{Ar}}(T_g,\tau_g)\Big)\Big)\Big)\cr
&-\Big(\text{ch}_{\text{Ar}}\Big((g\circ
f)_*E,L^2(\rho,\tau_{g\circ
f})\Big)
-\text{ch}_{\text{Ar}}\Big(g_*(f_*E),
L^2(L^2(\rho,\tau_f),\tau_g)\Big)\Big)\cr
&\qquad(\text{by projection formula})\cr
=&(g\circ
f)_*\Big(\text{ch}(E,\rho)\cdot
\text{td}_{\text{BC}}(T.,\tau.)\Big)-
\text{ch}_{\text{BC}}\Big((g\circ f)_*E;L^2(\rho,\tau_{g\circ
f}),L^2(L^2(\rho,\tau_f),\tau_g)\Big)\cr
&\qquad(\text{by Theorem A.10(iii)}.\cr}$$
This completes the proof of the lemma.
\vskip 0.30cm
\noindent
(B.8) Finally, we are left with checking that 
$\text{ch}_{\text{BC}}'$ satisfies  Axiom 6 for relative
Bott-Chern secondary characteristic classes.
\vskip 0.20cm
\noindent 
{\bf Lemma.} {\it With the same notation as in B.2 and as in 
3.C.3, we have $$\eqalign{~&\lim_{t\to
\infty}\Big(\big(\text{ch}_{\text{BC}}'(E_t,\rho_t;g_t,\tau_{g_t})
-\text{ch}_{\text{BC}}'(E_t(-X),\rho_t';g_t,\tau_{g_t})\big)\cr
&+\big(\text{ch}_{\text{BC}}((g_t)_*(E_t);
L^2(\rho_t,\tau_t),\gamma_t)-\text{ch}_{\text{BC}}((g_t)_*(E_t(-X));
L^2(\rho_t',\tau_t),\gamma_t')\big)\Big)\cr
=&\big(\text{ch}_{\text{BC}}'(E_\infty,\rho_\infty;g_\infty,\tau_{g_\infty})-
\text{ch}_{\text{BC}}'(E_\infty(-X),\rho_\infty';g_\infty,\tau_{g_\infty})\big)+\cr
&\big(\text{ch}_{\text{BC}}((g_\infty)_*(E_\infty);
L^2(\rho_\infty,\tau_\infty),\gamma_\infty)
-\text{ch}_{\text{BC}}((g_\infty)_*(E_\infty(-X));
L^2(\rho_\infty',\tau_\infty),\gamma_\infty')\big).\cr}$$}
\noindent
{\it Proof.} By definition, Theorem A.10(v), and the fact that
the $L^2$-metrics are now replaced by a continuous family of
metrics,  it is sufficient to show that
$$\eqalign{~&\lim_{t\to\infty}\Big(\text{ch}_{\text{Ar}}(E_t,\rho_t)-
\text{ch}_{\text{Ar}}(E_t(-X),\rho_t')\Big)
\cdot
\text{td}_{\text{Ar}}(T_t,\tau_{g_t})\cr
=&\Big(\text{ch}_{\text{Ar}}(E_\infty,\rho_\infty)-
\text{ch}_{\text{Ar}}(E_\infty(-X),\rho_\infty')\Big)
\cdot
\text{td}_{\text{Ar}}(T_\infty,\tau_{g_\infty}).\cr}$$
But this is a direct consequence of Proposition A.11. In fact,
by Proposition A.11, for all $t\in {\Bbb P}^1$,
$$\eqalign{~&\Big(\text{ch}_{\text{Ar}}(E_t,\rho_t)-
\text{ch}_{\text{Ar}}(E_t(-X),\rho_t')\Big)
\cdot
\text{td}_{\text{Ar}}(T_t,\tau_{g_t})\cr
=&i_t^*\Big(\Big(\text{ch}_{\text{Ar}}(E,D\rho)-\text{ch}_{\text{Ar}}(E,D\rho)\Big)
\cdot\text{td}_{\text{Ar}}(T_G(-{\log}\infty),\tau_G)\Big).\cr}$$
So, by the fact that there exists a well-defined pull-back
morphisms  $\text{CH}_{\text{Ar}}({\Bbb
P}^1)\to\text{CH}_{\text{Ar}}(W)$ and that on ${\Bbb P}^1$, the
arithmetic cycle
$\text{div}_{\text{Ar}}(z):=(0-\infty,-\log|z|^2)$ is rational
equivalent to zero, we have
$$\eqalign{~&\Big(\Big(\text{ch}_{\text{Ar}}(E(B_XZ),D\rho)-
\text{ch}_{\text{Ar}}(E(B_XZ-X\times {\Bbb P}^1),D\rho_0)\Big)\cdot
\text{td}_{\text{Ar}}(T_G(-\log\infty),\tau_G)\Big)\cr
&\qquad\cdot(W_0-W_\infty,
-\text{Log}|z|^2)=0.\cr}$$ Here $\text{Log}|z|$ denotes the
pull-back of $\log|z|$ on $W$. Note that $W_\infty$ has two
components which intersects properly, the restriction of $0\to
E(B_XZ-X\times {\Bbb P}^1)\to E(B_XZ)$ together with metrics
$D\rho'$ and $D\rho$ results an isometry
$(E_\infty',\rho_\infty'')\simeq (E_\infty'',\rho_\infty''')$,
$X\times {\Bbb P}^1$ is away from
$B_XZ$, and $E_\infty-E_\infty(-X)$ is supported on $X\subset
{\Bbb P}$,
$$\eqalign{~&(\Big(\text{ch}_{\text{Ar}}(E(B_XZ),D\rho)-
\text{ch}_{\text{Ar}}(E(B_XZ-X\times {\Bbb
P}^1),D\rho')\Big)\Big|_{W_\infty}\cr
=&(\Big(\text{ch}_{\text{Ar}}(E(B_XZ),D\rho)-
\text{ch}_{\text{Ar}}(E(B_XZ-X\times {\Bbb
P}^1),D\rho')\Big)\Big|_{\Bbb P}.\cr}$$ Hence, in particular,
$$\eqalign{~&\Big(\text{ch}_{\text{Ar}}(E_\infty,\rho_\infty)-
\text{ch}_{\text{Ar}}(E_\infty(-X),\rho_\infty')\Big)
\cdot
\text{td}_{\text{Ar}}(T_\infty,\tau_{g_\infty})\cr
&-
\Big(\text{ch}_{\text{Ar}}(E_0,\rho_0)-
\text{ch}_{\text{Ar}}(E_0(-X),\rho_0')\Big)
\cdot
\text{td}_{\text{Ar}}(T_0,\tau_{g_0})\cr
=&\int_{{\Bbb
P}^1}\Big(\Big(\text{ch}_{\text{Ar}}(E(B_XZ),D\rho)-
\text{ch}_{\text{Ar}}(E(B_XZ-X\times {\Bbb P}^1),D\rho')\Big)
\cr
&\qquad\cdot
\text{td}_{\text{Ar}}(T_G(-\log\infty),\tau_G)\Big)
\cdot[\text{Log}|z|^2].\cr}$$ Similarly, 
$$\eqalign{~&\Big(\text{ch}_{\text{Ar}}(E_t,\rho_t)-
\text{ch}_{\text{Ar}}(E_t(-X),\rho_t')\Big)
\cdot
\text{td}_{\text{Ar}}(T_t,\tau_{g_t})\cr
&-
\Big(\text{ch}_{\text{Ar}}(E_0,\rho_0)-
\text{ch}_{\text{Ar}}(E_0(-X),\rho_0')\Big)
\cdot
\text{td}_{\text{Ar}}(T_0,\tau_{g_0})\cr
=&\int_{{\Bbb
P}^1}\Big(\Big(\text{ch}_{\text{Ar}}(E(B_XZ),D\rho)-
\text{ch}_{\text{Ar}}(E(B_XZ-X\times {\Bbb P}^1),D\rho')\Big)\cr
&\qquad\cdot
\text{td}_{\text{Ar}}(T_G(-\log\infty),\tau_G)\Big)
\cdot[\text{Log}|{{zt}\over {z-t}}|^2].\cr}$$
Easily, one sees that $\lim_{t\to\infty}|{{zt}\over {z-t}}|^2
=|z|^2.$ This completes the proof of the lemma.
\vskip 0.20cm
\noindent
(B.9) Recall that $(f:X\to Y; E,\rho;T_f,\tau_f)$ is called a
properly metrized datum 
if $f:X\to Y$ is a smooth morphism of compact
K\"ahler manifolds,
$(E,\rho)$ is an $f$-acyclic hermitian vector bundle, and  
$\tau_f$ is a hermitian metric on the relative tangent bundle
$T_f$ of $f$ such that the induced metrics on all fibers of $f$
are K\"ahler. Then all in all, with the space $\bar A$ defined in
B.1,  what we have already proved is the following
\vskip 0.20cm
\noindent
{\bf Theorem.} ({\bf The Weak Existence and Strong Uniqueness 
for Relative Bott-Chern Secondary Characteristic Classes})
{\it  There exists a construction $\text{ch}_{\text{BC}}$
satisfies the six axioms for relative Bott-Chern secondary
characteristic classes with values in $\overline A$.
Moreover, if there  are two constructions
$\text{ch}_{\text{BC}}$ and
$\text{ch}_{\text{BC}}'$ which satisfy these six axioms, 
then, there exists an additive characteristic classes $R$
such that,  for all properly metrized data
$(f:X\to Y;E,\rho;T_f,\tau_f)$, in $\overline A(Y)$,
$$\text{ch}_{\text{BC}}'(E,\rho;f,\tau_f)=
\text{ch}_{\text{BC}}(E,\rho;f,\tau_f)+a\Big(f_*\Big(\text{ch}(E)\cdot
\text{td}(T_f)\cdot R(T_f)\Big)\Big).$$}
\vskip 0.20cm
\noindent
(B.10) We end this paper by the following remark. There is in
fact another effective construction for  relative
Bott-Chern secondary characteristic classes by using heat
kernels and Mellin transform. Even though such an alternative
construction involves only analysis, yet it offers
us a strong existence of the relative Bott-Chern secondary
characteristic classes, i.e., the relative Bott-Chern classes
at the level of $\tilde A$. Sure, once this is done, then the
arithmetic Grothendieck-Riemann-Roch theorem may be viewed as
a direct consequence of the uniqueness for relative
Bott-Chern classes. We will study it carefully
in another occasion.
\vskip 1.00cm
\noindent
{\bf Acknowledgement}. This work was started when I was at MPI
f\"ur Mathematik at Bonn and has been completed while I am at
Mathematics Department of Osaka University. I would like thank
both institutes for their hospitality and support. I also would
like to thank Fujiki, Lang, and Ueno for their
interests and encourgement. Obviously, the works done by
Bott-Chern, Faltings and  Gillet-Soul\'e provide the essential
ingredients without which it is impossible to develop a theory of
relative Bott-Chern classes.
\vfill\eject
\vskip 0.20cm
\noindent
{\bf REFERENCES}
\vskip 0.20cm
\noindent
\item {[Ar1]} ARAKELOV, S.J.: Intersection theory of divisors
on an arithmetic surface, Math. USSR Izvestija {\bf 8}, 1974, pp.
1167-1180
\vskip 0.20cm
\noindent
\item {[Ar2]} ARAKELOV, S.J.: Theory of intersections on an
arithmetic surface, {\it Proc.  Intrnl. Cong. of Mathematicians},
Vancouver, 1975, Vol. {\bf 1}, pp. 405-408
\vskip 0.20cm
\noindent
\item {[BC]} BOTT, R., CHERN, S.S.: Hermitian vector bundles and 
the equidistribution of the zeroes of their holomorphic
sections. Acta Math. {\bf 114} (1965), pp. 71-112
\vskip 0.20cm
\noindent
\item {[BT]} BOTT, R., TU, L. W.: {\it Differential forms in
algebraic topology}, Graduate Texts in Mathematics, {\bf 82},
Springer-Verlag, New York-Berlin, 1982. xiv+331 pp.
\vskip 0.20cm
\noindent
\item {[De1]} DELIGNE, P.: {\it Equations diff\'erentielles a
points  singuliers r\'eguliers}, Springer Lecture Notes in
Mathematics, no. {\bf 163}, 1970
\vskip 0.20cm
\noindent
\item {[De2]} DELIGNE, P.: Le d\'eterminant de la cohomologie,
in {\it Current trands in Arithmetic Algebraic Geometry},
Contemporary Math., {\bf 67} (1987), pp. 93-178
\vskip 0.20cm
\noindent
\item {[Fa]} FALTINGS, G.: {\it Lectures on the arithmetic
Riemann-Roch theorem}, noted by S. W. Zhang, Ann. of Math. Study
{\bf 127} (1992)
\vskip 0.20cm
\noindent
\item {[GS1]} GILLET, H.,  SOUL\'E, CH.: Arithmetic Intersection
Theory, Publ. Math. IHES {\bf 72} (1990), pp. 93-174
\vskip 0.20cm
\noindent
\item {[GS2]} GILLET, H., SOUL\'E, CH.: Characteristic classes
for algebraic vector bundles with hermitian metric, Annals of
Math., {\bf 131}, pp. 163-203, pp. 205-238 (1990)
\vskip 0.20cm
\noindent
\item {[We1]} WENG L.: Arithmetic Riemann-Roch Theorem for
smooth morphism: An approach with relative Bott-Chern secondary
characteristic forms, MPI preprint, 1991
\vskip 0.20cm
\noindent
\item {[We2]} WENG, L. {\it Relative Bott-Chern secondary
characteristic objects and an Arithmetic Riemann-Roch Theorem},
monograph, 1992 (= MPI preprints, 1994)
\end